\documentclass[11pt]{article}
\usepackage{amsfonts}
\usepackage{amssymb}
\usepackage{graphicx}
\usepackage{amsmath}
\usepackage{amsthm}
\usepackage{color}
\usepackage{multirow}
\usepackage{url}
\usepackage{float}
\usepackage{booktabs} 
\usepackage{makecell}
\usepackage{listings}
\usepackage{subcaption}
\usepackage{algorithm}
\usepackage{enumerate}
\usepackage{tikz}

\usepackage{algpseudocode}
\lstnewenvironment{Rcode}{\lstset{language=R}}{}
\newtheorem{theorem}{Theorem}[section]

\newtheorem{corollary}[theorem]{Corollary}

\newtheorem{lemma}[theorem]{Lemma}

\newtheorem{remark}[theorem]{Remark}

\definecolor{battleshipgrey}{rgb}{0.52, 0.52, 0.51}
\definecolor{navyblue}{rgb}{0.0, 0.0, 0.5}
\definecolor{arsenic}{rgb}{0.23, 0.27, 0.29}
\definecolor{oldmauve}{rgb}{0.4, 0.19, 0.28}
\usepackage[colorlinks=TRUE]{hyperref}
\hypersetup{
	colorlinks=true,
	linkcolor=navyblue,
	filecolor=blue,      
	urlcolor=oldmauve,
	citecolor=navyblue,
	pdftitle={Overleaf Example},
	pdfpagemode=FullScreen,
}
\usepackage{natbib}
\setcitestyle{authoryear}
\begin{document}

	%\title{ \bf  Variable selection in the semi-functional partial linear single-index model\vspace{10pt} }
	\title{ \bf Sparse semiparametric regression when predictors are mixture of functional and high-dimensional variables}
		\author{Silvia Novo{$^a$}\footnote{Corresponding author email address: \href{s.novo@udc.es}{s.novo@udc.es}} \hspace{2pt} Germ\'{a}n Aneiros{$^b$} \hspace{2pt} Philippe Vieu{$^c$} \\		
		{\normalsize $^a$ Department of Mathematics, MODES, CITIC, Universidade da Coruña, A Coruña, Spain}\\
		{\normalsize $^b$ Department of Mathematics, MODES, CITIC, ITMATI, Universidade da Coruña, A Coruña, Spain}\\
		{\normalsize $^c$ Institut de Math\'{e}matiques, Universit\'e  Paul Sabatier, Toulouse, France}
	}
	
	\date{}
	\maketitle
	
	\begin{abstract} 
	This paper aims to front with dimensionality reduction in regression setting when the predictors are a mixture of functional variable and high-dimensional vector. A flexible model, combining both sparse linear ideas together with semiparametrics, is proposed. A wide scope of asymptotic results is provided: this covers as well rates of convergence of the estimators as asymptotic behaviour of the variable selection procedure. Practical issues are analysed through finite sample simulated experiments while an application to Tecator's data illustrates the usefulness of our methodology.
	\end{abstract}
	
	\vspace{1cm}
	\noindent \textit{Keywords: } Functional data analysis; Big data analysis; Variable selection; Sparse model; Dimension reduction; Functional single-index model; Semiparametrics

	\section{Introduction}
	\label{intro}
	During the two last decades, Functional Data Analysis (FDA) rose from a rather confidential field to one of the major topics in Statistics. As of today, as attested for instance by various recent surveys on FDA (see eg \citealt*{cue14}, \citealt*{goiv16} and \citealt*{acfgv19}), most of multivariate data analysis methodologies have been tentatively adapted to functional data. Because of the infinite-dimensionality feature of the random variables  involved in FDA, one of the main issues for insuring the good behaviour of any functional statistical procedure is to control, in one way or another, the dimensionality of the model (see eg \citealt*{vie18}). As a matter of fact, this dimensionality challenge is not so far from what exists in the related field of Big Data analysis in which traditionally the statistical variable is a high-dimensional vector. In the recent past, the necessary links between the two fields have been highlighted as well by the Big Data community (see eg \citealt*{san18}) as by the FDA one (see eg \citealt*{goiv16} and \citealt*{acfgv19}).
	
	In many fields of applications, one could have to deal with data consisting of mixtures of functional and high-dimensional variables and the statistical methodologies to be constructed have to cross both fields of FDA and Big Data. Our paper is part of this category, since the purpose is to develop models for regression problems involving some scalar response (let's say $Y$) and predictors composed of some functional variable, $\mathcal{X}$, and some high-dimensional vector, $(X_1,\ldots ,X_{p_n})$. The literature on functional and/or high-dimensional regression is usually balancing between parametric and nonparametric modelling. On the one hand, a parametric approach (see eg \citealt*{cars11} for general presentation) is  fewly impacted by dimensionality effect but suffers in counterpart from high lack of flexibility. 
	On the other hand, the high degree of flexibility of nonparametric modelling (see eg \citealt*{gee11} or \citealt*{linv18} for general discussion) is unfortunately going together with dramatically bad dimensionality effects. 
	
	In the multivariate literature  popular ways for dealing with dimensionality effects include principal component and partial least square ideas (mainly in linear modelling) or semiparametric ones (in non-linear modelling). These ideas have been recently adapted to the functional setting: the readers interested in techniques for  principal components or  partial least squares in functional linear regression may have a look at  \citealt*{pres05}, \citealt*{reio07}, \citealt*{delh12}, \citealt*{aguetal16} and \citealt*{kraetal08}, while they can find advances in functional semiparametrics in the survey by \citealt*{goiv14}. Functional semiparametric ideas  will be one of the main features of the model to be presented in this paper.
	%Semiparametrics are nice trade-off (see \citealt*{goiv14} for a survey on functional semiparametrics) and are one .
	
	Based on these considerations, our model has to take into account three important features of our problem: i) firstly, additive ideas are needed to separate the effects of the functional predictor, $\mathcal{X}$, from those of  multivariate predictor, $(X_1,\ldots ,X_{p_n})$; ii) secondly, sparse ideas are needed in order to control the high number of variables, $p_n$ (which is allowed to go to infinity as $n$ does), involved in the multivariate predictor; iii) finally, functional semiparametric ideas are required for modelling the effect of the infinite-dimensional predictor, $\mathcal{X}$. This leads to the so-called sparse semi-functional partial linear single-index model (SSFPLSIM) that will be presented in Section \ref{model}. A variable selection method and estimators of the components of the model will be constructed along Section \ref{estimates}, while a wide set of asymptotics will be provided in Section \ref{asymptotics}. Finite sample behaviour of the method will be assessed through Monte Carlo experiments in Section \ref{MC}. In addition, Section \ref{real-data} provides an application to Tecator's data. Technical proofs and lemmas are gathered in a supplementary file.    
	
	\section{The model}\label{model}
	The SSFPLSIM is defined by the relationship
	\begin{equation}
		\label{modelo}
		Y_i=X_{i1}\beta_{01}+\dots+X_{ip_n}\beta_{0p_n}+m\left(\left<\theta_0,\mathcal{X}_i\right>\right)+\varepsilon_i, \ \forall i=1,\dots,n,
	\end{equation}
	where $Y_i$ denotes a scalar response, $X_{i1},\dots,X_{ip_n}$ are random covariates taking values in $\mathbb{R}$ and $\mathcal{X}_i$ is a functional random covariate valued in a separable Hilbert space $\mathcal{H}$ with inner product $\left\langle \cdot, \cdot \right\rangle$. In this equation,
	$\pmb{\beta}_0=(\beta_{01},\dots,\beta_{0p_n})^{\top} \in \mathbb{R}^{p_n}$, $\theta_0\in\mathcal{H}$ and $m(\cdot)$ are a vector of unknown real parameters, an unknown functional direction and an unknown smooth real-valued function, respectively. Finally, $\varepsilon_i$ is the random error, which verifies
	\begin{equation}
		\mathbb{E}\left(\varepsilon_i|X_{i1},\dots,X_{ip_n},\mathcal{X}_i\right)=0.
		\label{centred_error}
	\end{equation} 
	To ensure identifiability, it is usual to assume that either $\left<\Gamma \theta_0,\theta_0\right>=1$, where $\Gamma$ denotes the covariance function of the functional variable $\mathcal{X}$ and $(\Gamma \theta_0)(t)=\left<\Gamma(\cdot,t),\theta_0\right>$, or $\left<\theta_0,\theta_0\right>=1$, and that, for some arbitrary $t_0$ in the domain of $\theta_0$, one has $\theta_0(t_0)>0$. Note that these conditions are common in the literature on this kind of models (see eg \citealt*{wang_2016} for model (\ref{modelo}), or \citealt*{ding_2017}, \citealt*{liang10} and \citealt*{wangzhu_2017} for related models).
	
	Finally, the model needs to incorporate the situation of high number of  covariates, that is when  %$p_n$ depends on the sample size, in such a way that 
	$p_n \rightarrow \infty$ as $n \rightarrow \infty$, and it has therefore to include some sparsity parameter. This is the role of the additional integer parameter $s_n$ which is the number of relevant covariates (that is, the covariates associated with $\beta_{0j} \neq 0$) that will be supposed to be much smaller than $p_n$ (see technical assumptions later). With such a sparse regression model, the aim is to construct a procedure to select the relevant covariates and, simultaneously, to estimate their corresponding effects, $\beta_{0j}$. From a theoretical point of view, the challenge is double: i) obtain consistency of the model selection procedure; ii) get  the same rate of convergence for the estimator of $\pmb{\beta}_0$ in (\ref{modelo}) as those obtained in the standard literature for the linear model
	\begin{equation}
		Y_i=X_{i1}\beta_{01}+\dots+X_{ip_n}\beta_{0p_n}+\varepsilon_i, \ \forall i=1,\dots,n,
		\label{modelo.lineal}
	\end{equation}
	or for the semi-functional partial linear model:
	\begin{equation}
		Y_i=X_{i1}\beta_{01}+\dots+X_{ip_n}\beta_{0p_n}+r\left(\mathcal{X}_i\right)+\varepsilon_i, \ \forall i=1,\dots,n.
		\label{modelo.sfpl}
	\end{equation}
	\noindent Recall that such a rate can be expressed as $O_p(\sqrt{s_n}n^{-1/2})$ (see eg \citealt*{fanlv2011} and \citealt*{aneiros_2015} for models (\ref{modelo.lineal}) and (\ref{modelo.sfpl}), respectively). The variable selection procedure is presented in Section \ref{estimates} and their asymptotics are stated in Section \ref{asymptotics}. Once the linear part of the model is dealt with, the single-functional component $m\left(\left<\theta_0,.\right>\right)$ will be estimated with univariate nonparametric rate of convergence (see Section \ref{asymptotics}).

	\section{The penalized least-squares estimators}\label{estimates}
	\subsection{Some initial notation}
	\label{some-notation}
	Throughout this paper we will use the notation
	$$\pmb{X}_i=\left(X_{i1},X_{i2},\dots,X_{ip_n}\right)^{\top}, \ \pmb{X}=\left(\pmb{X}_1,\dots,\pmb{X}_n\right)^{\top} \mbox{ and } \pmb{Y}=\left(Y_1,\dots,Y_n\right)^{\top}.$$
	In addition, for any $(n\times q)$-matrix $\pmb{A}$ $(q\geq 1)$ and $\theta\in\mathcal{H}$, we denote
	\begin{equation*}
		\widetilde{\pmb{A}}_{\theta}=\left(\pmb{I}-\pmb{W}_{h,\theta}\right)\pmb{A}, \mbox{ where } \pmb{W}_{h,\theta}=\left(w_{n,h,\theta}(\mathcal{X}_i,\mathcal{X}_j)\right)_{i,j}
	\end{equation*}
	with $w_{n,h,\theta}(\cdot,\cdot)$ being the weight function
	\begin{equation*}
		w_{n,h,\theta}(\chi,\mathcal{X}_i)=\frac{K\left(d_{\theta}\left(\chi,\mathcal{X}_i\right)/h\right)}{\sum_{j=1}^n K\left(d_{\theta}\left(\chi,\mathcal{X}_j\right)/h\right)},
	\end{equation*}
	where $K:\mathbb{R}^+\rightarrow\mathbb{R}^+$ is a kernel function, $h>0$ is a smoothing parameter and $d_{\theta}(\cdot,\cdot)$ is the semimetric defined as
	$$d_{\theta}\left(\chi,\chi'\right)=\left|\left<\theta,\chi-\chi'\right>\right|, \ \forall \chi,\chi'\in\mathcal{H}.$$
	
	\subsection{The estimators}
	\label{estimators}
	For simultaneously estimating the linear $\beta$-parameters and selecting the relevant $X$-covariates in the model (\ref{modelo}), we will use a penalized least-squares approach. The first idea is to transform the SSFPLSIM into a linear one by extracting from $Y_i$ and $X_{ij}$ ($j=1,\ldots,p_n$) the effect of the functional covariate $\mathcal{X}_i$ when projected on the direction $\theta_0$. Specifically, the fact that
	\begin{equation}
		Y_i-\mathbb{E}\left(Y_i|\left<\theta_0,\mathcal{X}_i\right>\right)=\left(\pmb{X}_i-\mathbb{E}\left(\pmb{X}_i|\left<\theta_0,\mathcal{X}_i\right>\right)\right)^{\top}\pmb{\beta}_0+\varepsilon_i, \ \forall i=1,\dots,n,\label{mod_lineal}
	\end{equation}
	allows to consider the following approximate linear model:
	\begin{equation}
		\widetilde{\pmb{Y}}_{\theta_0}\approx\widetilde{\pmb{X}}_{\theta_0}\pmb{\beta}_0+\pmb{\varepsilon}, \label{linear}
	\end{equation}
	where $\pmb{\varepsilon}=\left(\varepsilon_1,\dots,\varepsilon_n\right)^{\top}$ (note that, to obtain (\ref{linear}), the conditioned expectations in (\ref{mod_lineal}) were estimated by means of functional nonparametric techniques). Then, in a second step, the penalized least-squares approach is applied to model (\ref{linear}). In this way, the considered penalized profile least-squares function is defined as
	%Specifically, $\pmb{\beta}_0$ and $\theta_0$ are estimated by considering a  (local) minimizer, $(\widehat{\pmb{\beta}}_0, \widehat{\theta}_0)$, of the penalized profile least-squares function
	\begin{equation}
		\mathcal{Q}\left(\pmb{\beta},\theta\right)=\frac{1}{2}\left(\widetilde{\pmb{Y}}_{\theta}-\widetilde{\pmb{X}}_{\theta}\pmb{\beta}\right)^{\top}\left(\widetilde{\pmb{Y}}_{\theta}-\widetilde{\pmb{X}}_{\theta}\pmb{\beta}\right)+n\sum_{j=1}^{p_n}\mathcal{P}_{\lambda_{j_n}}\left(|\beta_j|\right),
		\label{func_minimizar}
	\end{equation}
	where $\pmb{\beta}=(\beta_1,\ldots,\beta_{p_n})^{\top}$, $\mathcal{P}_{\lambda_{j_n}}\left(\cdot\right)$ is a penalty function and $\lambda_{j_n} > 0$ is a tuning parameter (note that the role of the summatory in (\ref{func_minimizar}) is to penalize the presence of non zero $\beta$-parameters; in fact, under suitable conditions on $\mathcal{P}_{\lambda}$ (see eg \citealt*{fanli_2001}), the penalized least-squares estimators produce sparse solutions (many estimated coefficients are zero)). At this moment it is noteworthy, on the one hand, the fact that the objective function $\mathcal{Q}$ in (\ref{func_minimizar}) is not necessarily convex. This is the reason why, as usual in the related literature (see eg \citealt*{fanli_2001}, \citealt*{fanpeng_2004} and \citealt*{wangzhu_2017}), our asymptotic results in next section are focused on a local minimizer, $(\widehat{\pmb{\beta}}_0, \widehat{\theta}_0)$, of $\mathcal{Q}$ (in particular, existence of such local minimizer will be established in Theorem \ref{teorema1}). On the other hand, this parameter estimation procedure can be used also as a variable selection method in a simple way: if $\widehat{\pmb{\beta}}_{0j}$ is a non-null component of $\widehat{\pmb{\beta}}_0$, then $X_j$ is selected as an influential variable.

	Finally, after estimating $\pmb{\beta}_0$ and $\theta_0,$ we can deal with the estimation of the nonlinear function $m_{\theta_0}(\cdot)\equiv m\left(\left<\theta_0,\cdot\right>\right)$ in (\ref{modelo}). A natural way is employing again nonparametric procedures and smoothing the partial residuals $Y_i-\pmb{X}_i^{\top}\widehat{\pmb{\beta}}_0$. For that, we introduce the following statistic: %estimator
	\begin{equation}
		\label{m_statistic}
		\widehat{m}_{\theta}\left(\chi\right)\equiv\widehat{m}\left(\left<\theta,\chi\right>\right)=\sum_{i=1}^n w_{n,h,\theta}(\chi,\mathcal{X}_i)\left(Y_i-\pmb{X}_i^{\top}\widehat{\pmb{\beta}}_0\right),
	\end{equation}
	the estimator of $m_{\theta_0}\left(\chi\right)$ being $\widehat{m}_{\widehat{\theta}_0}\left(\chi\right)$. Note that the same bandwidth, $h$, is used to estimate both the functional index $\theta_0$ and the parameter ${\pmb{\beta}}_0$ from (\ref{func_minimizar}), as well as to estimate the smooth real-valued function $m(\cdot)$ from (\ref{m_statistic}). Despite the partial residuals in (\ref{m_statistic}) could be smoothed by considering a different bandwidth, we have adopted the more usual procedure of using the same bandwidth twice (see eg \citealt*{{liang10}} for the case of a non functional partial linear single-index model). Extension to the case with different bandwidths does not involve any extra difficulties.

	\section{Asymptotic theory}\label{asymptotics}
	
	\subsection{Some additional notation}
	\label{notation}
	Let us first introduce some additional notation to be used in the results that will be presented in this paper, as well as in their proofs:
	\begin{itemize}
		\item We will denote  $J_n=\{1,\dots,p_n\}$ and $S_n=\{j \in J_n; \ \beta_{0j} \not=0\}$.  
		In addition $s_n$ will mean card$(S_n)$.
		
		\item Given any vector $\pmb{v}\in\mathbb{R}^{p_n}$ and any $p_n\times p_n$ matrix $\pmb{M}$, $\pmb{v}_{S_n}$ and $\pmb{M}_{S_n\times S_n}$ denote the vector and the matrix obtained from $\pmb{v}$ and $\pmb{M}$ by retaining only the components corresponding to the index sets $S_n$ and $S_n\times S_n$, respectively.
		
		\item For any $\theta \in \mathcal{H}$, and for $1\leq i\leq n$, $1\leq j\leq p_n$, we denote
		\begin{equation*}
			g_{0,\theta}(\mathcal{X}_i)=\mathbb{E}\left(Y_i|\left<\theta,\mathcal{X}_i\right>\right), \ g_{j,\theta}(\mathcal{X}_i)=\mathbb{E}\left(X_{ij}|\left<\theta,\mathcal{X}_i\right>\right),
		\end{equation*}
		and 
		\begin{equation*}
			\pmb{\eta}_{i,\theta_0}=\left(\eta_{i1,\theta_0},\dots,\eta_{ip_n,\theta_0}\right)^{\top}, \textrm{ where } \eta_{ij,\theta_0}=X_{ij}-g_{j,\theta_0}(\mathcal{X}_i).
		\end{equation*} 
		
		\item  $\Delta_{min}(\pmb{M})$ and $\Delta_{max}(\pmb{M})$ denote the smallest and the largest eigenvalues of the matrix $\pmb{M}$, respectively.
		\item The symbol $||\cdot||$ is used for denoting the $L_2$ norm of vectors and matrices. % (its norm induced by such vector norm).
		The same symbol is also employed for denoting the norm induced by the inner product $\left<\cdot,\cdot\right>$. Specifically:
		\begin{align*}
			&||\pmb{a}||=\left(a_1^2+\dots+a_q^2\right)^{1/2} \ \textrm{for} \ \pmb{a}=\left(a_1,\dots,a_q\right)^{\top}\in\mathbb{R}^q,
			\\
			&||\pmb{A}||=\max_{0\not=\pmb{x}\in\mathbb{R}^q}\frac{||\pmb{Ax}||}{||\pmb{x}||} \textrm{ for any }r\times q \textrm{ matrix } \pmb{A}
			\\
			& \textrm{and} \\
			&||\chi||=\left<\chi,\chi\right>^{1/2} \textrm{ for any }\chi \in\mathcal{H}.
		\end{align*}	
		
		\item  $\forall \ \chi, \theta \in \mathcal{H}$ and $\forall \epsilon>0$, we will use the notation:
		$$
		B_\theta(\chi,\epsilon)=\{\chi' \in \mathcal{H}; \ d_\theta(\chi,\chi')<\epsilon\}, \
		\phi_{\chi,\theta}(\epsilon)=\mathbb{P}\left(\mathcal{X} \in B_\theta(\chi,\epsilon)\right)
		$$
		and
		\begin{equation*}
			B(\theta,\epsilon)=\{\theta' \in \mathcal{H}; \ d(\theta,\theta')<\epsilon\},
		\end{equation*}
		where,  $\forall \ \chi,\chi'\in\mathcal{H}$,
		%$$d\left(\chi,\chi'\right)=\left<\chi-\chi',\chi-\chi'\right>^{1/2} .$$
		$$d\left(\chi,\chi'\right)=||\chi-\chi'|| .$$
		
	\end{itemize}
	
	\subsection{Assumptions}
	\label{assumptions}
	In order to state rates of convergence of the proposed estimators and model selection consistency, we will use a large number of assumptions (some of them very technical). Such number is directly linked to the complexity of the model and the results to be obtained. These assumptions, that will be justified in next Remark \ref{r_assumptions}, are the following:
	
	\begin{description}
		\item[\textit{Conditions on the set of values of $\mathcal{X}$ and the topologies induced by $d_{\theta}(\cdot,\cdot)$.}] The functional variable $\mathcal{X}$ is valued in some subset $\mathcal{C}$ of $%
		\mathcal{H}$ such that
		\begin{equation}
			\mathcal{C} \subset \bigcup_{k=1}^{N_{\mathcal{C},\epsilon}^{\theta}}B_\theta(\chi_{\epsilon,k}^{\theta},\epsilon), \ \forall \theta \in \Theta_n,
			\label{cover-C}
		\end{equation}
		where
		\begin{equation}
			\Theta_n=\{\theta \in \mathcal{H}; \ d(\theta,\theta_0)\leq v_n\} \mbox{ with } v_n\rightarrow 0 \mbox{ as } n\rightarrow \infty ,
			\label{Theta}
		\end{equation}
		and
		$N_{\mathcal{C},\epsilon}^{\theta}$ is the minimal number of open balls in $(\mathcal{H},d_\theta(\cdot,\cdot))$ of radius $\epsilon$ which are necessary to cover $\mathcal{C}$.
		
		\item[\textit{Conditions on the entropies and the balls in (\ref{cover-C}).}] Let us denote
		\begin{equation}
			N_{\mathcal{C},\epsilon}=\sup_{\theta \in \Theta_n} N_{\mathcal{C},\epsilon}^{\theta}, \ \psi_\mathcal{C}(\epsilon)=\log(N_{\mathcal{C},\epsilon}), \ 
			k_{(\theta,k,\epsilon)}^{0}=\arg\min_{k' \in \left\{1,\ldots,N_{\mathcal{C},\epsilon}^{\theta_0}\right\}}d(\chi_{\epsilon,k}^{\theta},\chi_{\epsilon,k'}^{\theta_{0}}) \label{varias}
		\end{equation}
		and, in the sake of brevity,
		$$
		\chi_{k}^{\theta}=\chi_{1/n,k}^{\theta} \mbox{ and } k^0=k_{(\theta,k,1/n)}^{0}.
		$$
		
		It is assumed that:
		\begin{equation}
			\exists \beta>1 \textrm{ such that } p_n\exp\left\{\left(1-\beta\log p_n\right)\psi_{\mathcal{C}}\left(\frac{1}{n}\right)\right\}\rightarrow 0 \textrm{ as } n\rightarrow \infty ,
			\label{fun_spaces_0}	
		\end{equation}
		and
		\begin{equation}
			\sup_{\theta \in \Theta_n} \max_{k \in \{1,\ldots,N_{\mathcal{C},1/n}^\theta\}}d(\chi_k^{\theta},\chi_{k^0}^{\theta_0})=O(1/n).
			\label{chi-chi_0}
		\end{equation}
		
		\item[\textit{Conditions on the small-ball probabilities.}]  There exist constants $C_1>0$,  $0<C_2\leq C_{3}<\infty$
		and a function $f:\mathbb{R}\longrightarrow(0,\infty)$ such that 
		%$\int_0^1f(sh)ds/f(h)>\rho_0>0$.
		\begin{equation}
			\int_{0}^{1}f \left( hs\right) ds>C _{1}f\left( h\right)%
			\label{h33}
		\end{equation}
		and
		\begin{equation}
			\forall \chi \in \mathcal{C} \mbox{ and } \forall \theta \in \Theta_n, \ C_2f(h)\leq\phi_{\chi,\theta}(h)\leq C_{3}f(h).
			\label{small_ball}
		\end{equation}
		
		\item[\textit{Conditions linking the entropies and the small-ball probabilities.}] There exists a constant $C_4>0$ such that, for $n$ large enough,
		\begin{equation}
			\psi_{\mathcal{C}}\left(\frac{1}{n}\right)\leq\frac{C_4 n f(h)}{\alpha_n\log p_n}, \mbox{ where } \alpha_n\rightarrow\infty \mbox{ as } n\rightarrow\infty.
			\label{entropy_0}
		\end{equation}
		
		\item[\textit{Conditions on the kernel $K$.}] 
		\begin{align}
			&K \textrm{ is Lipschitz continuous on its support $[0,1)$ and}, \textrm{
				if $K(1)=0$},\nonumber\\ 
			&K \textrm{ also satisfies } -\infty<C_5<K'(\cdot)<C_6<0,\nonumber\\
			&\mbox{where } C_5 \mbox{ and } C_6 \mbox{ denote constants.}
			\label{kernel}
		\end{align}
		
		\item[\textit{Conditions on the smoothness.}] For some constants $0\leq C_7<\infty$ and $\alpha>0$, $\forall (\chi,\chi')\in \mathcal{C}\times \mathcal{C}$, and $\forall z\in \{g_{0,\theta_0},g_{1,\theta_0},\dots,g_{p_n,\theta_0}\}$, it verifies that
		%\begin{eqnarray}
		\begin{equation}
			\left|z(\chi)-z(\chi')\right|\leq C_7d_{\theta_0}(\chi,\chi')^{\alpha}.%\nonumber\\
			\label{smooth}
		\end{equation}
		%\end{eqnarray}
		\item[\textit{Conditions on the random variables.}]
		
		\begin{equation}
			\left\{(Y_i,X_{i1},\dots,X_{ip_n},\mathcal{X}_i)\right\} \mbox{ are random vectors iid verifying model (\ref{modelo})}.
			\label{var_indep_1}
		\end{equation}
		\begin{equation}
			\{\pmb{\eta}_{i,\theta_0}\} \mbox{ and } \{\varepsilon_i\} \mbox{ are independents.}
			\label{var_indep_2}
		\end{equation}
		\begin{equation}
			\left<\mathcal{X},\mathcal{X}\right>^{1/2}<C_8, \mbox{where } C_8 \mbox{ denotes a positive constant.} \label{chi}
		\end{equation}
		
		\item[\textit{Conditions on the moments.}] Let $C_{9}$, $C_{\eta_{\theta_0}}$ and $C_{r_{\varepsilon}}$ be positive constants. $\forall r\geq 2$, there exists a continuous operator in $\mathcal{C}$, $\sigma_r(\cdot)$, such that $\forall \chi\in \mathcal{C},$
		\begin{equation}
			\max_{j \in \{1,\dots,p_n\}}\left\{\mathbb{E}\left(|Y_1|^r|\left<\theta_0,\mathcal{X}_1\right>=\left<\theta_0,\chi\right>\right),\mathbb{E}\left(|X_{1j}|^r|\left<\theta_0,\mathcal{X}_1\right>=\left<\theta_0,\chi\right>\right)\right\}<\sigma_r(\chi)<C_{9}.\label{mom_1}
		\end{equation}
		\begin{equation}
			\forall \ r\geq 2 \ \textrm{and} \ \forall 1\leq j\leq p_n,\ \mathbb{E}|\eta_{1j,\theta_0}|^r\leq C_{\eta_{\theta_0}}\left(\frac{r!}{2}\right). \label{mom_2}
		\end{equation}
		\begin{equation}
			\exists r_{\varepsilon}>4  \textrm{ such that }  \mathbb{E}\left|\varepsilon_1\right|^{r_{\varepsilon}}\leq C_{r_{\varepsilon}}.
			\label{mom_3}
		\end{equation}
		In addition, there exists a constant $C_{10}$ such that
		\begin{eqnarray}
			0<C_{10}<\Delta_{min}\left(\pmb{B}_{\theta_0S_n\times S_n}\right),
			\label{mom_4}
		\end{eqnarray}
		where $\pmb{B}_{\theta_0}=\mathbb{E}\left(\pmb{\eta}_{1\theta_0}\pmb{\eta}_{1\theta_0}^{\top}\right)$.
		
		\item[\textit{Conditions on the non null parameters and the penalty functions.}] Let $C_{11}$ and $C_{12}$ be positive constants. 
		\begin{equation}
			\mathcal{P}_{\lambda_{jn}}(\cdot) \mbox{ is a continuous and nonnegative function verifying } \mathcal{P}_{\lambda_{jn}}(0)=0. \label{u1}
		\end{equation}
		\begin{equation}
			\mathcal{P}_{\lambda_{jn}}(\cdot) \mbox{ is differentiable excepted perhaps at } 0. \label{u2}
		\end{equation}
		\begin{equation}
			\left|\mathcal{P}_{\lambda_{j_n}}''(a)-\mathcal{P}_{\lambda_{jn}}''(b)\right|\leq C_{11}|a-b|,\ \forall a,b>C_{12}\lambda_{jn}. \label{u3}
		\end{equation}
		\begin{equation}
			\lim\inf_{n\rightarrow\infty}\min_{j\in S_n^c}\left\{\lim \inf_{d\rightarrow 0^+}\frac{\mathcal{P}'_{\lambda_{jn}}(d)}{\lambda_{jn}}\right\}>0.
			\label{null_par1}
		\end{equation}
		
		Finally,
		\begin{equation}
			\min_{j\in S_n}\left\{\frac{\left|\beta_{0j}\right|}{\lambda_{jn}}\right\}\rightarrow\infty \textrm{ as } n\rightarrow \infty\label{null_par2}
		\end{equation}
		and
		\begin{equation}
			\max_{j\in S_n}\left\{\left|\beta_{0j}\right|\right\}=O(1).\label{null_par3}
		\end{equation}
	\end{description}
	
	\begin{remark}
		\label{r_assumptions}
		The hypotheses listed above are, in general, usual (or natural extensions of those) in the related literature. For instance, conditions (\ref{cover-C}), (\ref{fun_spaces_0}) and (\ref{entropy_0}) are related to the topology of $(\mathcal{C}, d_{\theta})$ and, in the particular case of known $\theta_0$, they are common when one needs to obtain uniform orders over $\mathcal{C}$ (see eg \citealt*{ferraty_2010} or \citealt*{aneiros_2015}). In the general case dealt here, where $\theta_0$ is unknown and one needs to control the behaviour of the profile function $\mathcal{Q}(\cdot,\cdot)$ (\ref{func_minimizar}) around $\theta_0$, conditions (\ref{cover-C}), (\ref{fun_spaces_0}) and (\ref{entropy_0}) are the natural extension of the corresponding to such particular case. In the same way, conditions (\ref{Theta}), (\ref{h33}) and (\ref{small_ball}) also allow to control the effect of $\theta$; more specifically, condition (\ref{Theta}) establishes the set of values of $\theta$ where the profile function $\mathcal{Q}(\cdot,\cdot)$ achieves a local minimum (see \citealt*{ma16}), while conditions (\ref{h33}) and (\ref{small_ball}) are natural extensions of usual assumptions (related to the concentration properties of the probability measure of the functional variable $\mathcal{X}$) from the case where $\theta_0$ is known (see eg \citealt*{ferraty_2010}) to the one of unknown $\theta_0$ (see eg \citealt*{ait} and \citealt*{novo}). In addition, conditions (\ref{kernel})-(\ref{mom_4}) are standard ones in nonparametric and semiparametric estimation of the regression function when functional covariates are present (see eg \citealt*{ferraty_2010}, \citealt*{aneiros_2015}, \citealt*{wang_2016}). Basically, they are mild conditions on the kernel, on the smoothness of the nonparametric components related to both the response variable and the scalar covariates, and on both the dependence within and the moments of the variables in the model. Focusing now on the conditions directly linked to the penalty procedure (conditions (\ref{u1})-(\ref{null_par3})), they are usual assumptions in the topic of variable selection using nonconcave penalized functions (see eg \citealt*{fanli_2001}, \citealt*{fanpeng_2004}, \citealt*{aneiros_2015}). It is noteworthy that a main role of these conditions is to produce sparse solutions; that is, automatically to set small estimated coefficients to
		zero to reduce model complexity. Note that, under some specific condition (see eg \citealt*{aneiros_2015}), the smoothly clipped absolute deviation (SCAD) penalty function (proposed in \citealt*{Fan97}) verifies our assumptions. Finally, the condition (\ref{chi-chi_0}) is really specific to the functional setting addressed here and, therefore, requires a deeper discussion. It will be discussed in a more general setting in Section B.1 (see Remark 5 in the supplementary file).
	\end{remark}
	
	\subsection{Results}
	Our first result focuses on both the existence and rate of convergence of a local minimizer of the penalized least-squares objective function $\mathcal{Q}\left(\pmb{\beta},\theta\right)$ (see (\ref{func_minimizar})). Let us denote
	\begin{equation}
		\label{not-4.3}
		\delta_n=\max_{j\in S_n}\left\{\left|\mathcal{P}_{\lambda_{j_n}}'\left(|\beta_{0j}|\right)\right|\right\}, \ \rho_n=\max_{j\in S_n}\left\{\left|\mathcal{P}_{\lambda_{j_n}}''\left(|\beta_{0j}|\right)\right|\right\}  \mbox{ and } u_n=\sqrt{s_n}\left(n^{-1/2}+\delta_n\right).
	\end{equation}
	
	\begin{theorem}
		\label{teorema1}
		Assume that the assumptions (\ref{centred_error}), (\ref{cover-C}), (\ref{Theta}) and (\ref{fun_spaces_0})-(\ref{null_par3}) hold. Assume, in addition, that $p_n\rightarrow\infty$ as $n\rightarrow\infty$, $p_n=o\left(n^{1/2}\right)$ and
		\begin{align}
			&\max\left\{ns_n^2h^{4\alpha},s_n h^{\alpha}\log n\right\}=O(1),\nonumber\\
			&s_n^2\log p_n \log^2n=O\left(\frac{nf(h)}{\psi_{\mathcal{C}}\left(1/n\right)}\right),\nonumber\\
			&s_n^2\log^2 n=O\left(n\left(\frac{f(h)}{\psi_{\mathcal{C}}\left(1/n\right)}\right)^2\right),\nonumber\\
			&ns_nv_n=O(hf(h)) \nonumber \\
			&and \nonumber \\
			&\max\left\{\rho_n,\frac{u_n}{\min_{j\in S_n^c}\left\{\lambda_{jn}\right\}},\frac{u_n\Delta_{max}^{1/2}(\pmb{B}_{\theta_0})}{\min_{j\in S_n^c}\left\{\lambda_{jn}\right\}},\frac{n^{-1/2+1/r_{\epsilon}}\log n}{\min_{j\in S_n^c}\left\{\lambda_{jn}\right\}}\right\}=o(1). \nonumber
		\end{align}
		Then, there exists a local minimizer $\left(\widehat{\pmb{\beta}}_0,\widehat{\theta}_0\right)$ of $\mathcal{Q}\left(\pmb{\beta},\theta\right)$ such that
		\begin{eqnarray}
			\left\|\widehat{\pmb{\beta}}_0-\pmb{\beta}_0\right\|=O_p(u_n) \mbox{ and }
			\left\|\widehat{\theta}_0-\theta_0\right\|=O_p(v_n). \nonumber
		\end{eqnarray}
		(Note that $v_n$ was defined in (\ref{Theta}))
	\end{theorem}
	
	\begin{remark}
		\label{r_theorem1}
		Theorem \ref{teorema1} can be seen, in a certain sense, as an extension of Theorem 3.1 in \cite{aneiros_2015} from the case $\Theta_n=\{\theta_0\}$ (i.e., $v_n=0$ in (\ref{Theta}); equivalently, $\theta_0$ known) to the case where $\{\theta_0\} \subset \Theta_n$ (i.e., $\theta_0$ unknown). For that, all that one must do is to consider the results in \cite{aneiros_2015} when the semimetric $d_{\theta_0}(\cdot,\cdot)$ is used. From Theorem \ref{teorema1} we have that the rate of convergence ($u_n$) achieved by the local minimizer $\widehat{\pmb{\beta}}_0$ is the same as that obtained in the least complex scenario studied in \cite{aneiros_2015} (as well as in the linear model considered in \citealt*{fanlv2011}, where $\delta_n=0$), this being one of our main aims. Naturally, for this to be possible, it is necessary to have a very good estimator of the parameter $\theta_0$. Such estimator is obtained by means of the local minimizer $\widehat{\theta}_0$ (note that the local feature of the minimizers plays a main role to obtain fast rates of convergence).
	\end{remark}
	
	Our second result states the model selection consistency. Let us denote
	\begin{equation*}
		\widehat{S}_n=\left\{ j\in J_n; \widehat{\beta}_{0j}\not = 0\right\},
	\end{equation*}
	where $\widehat{\pmb{\beta}}_0=(\widehat{\beta}_{01},\dots,\widehat{\beta}_{0p_n})^{\top}$ is the estimator in Theorem \ref{teorema1}.
	\begin{theorem}[Model selection consistency]\label{teorema2}
		Under assumptions in Theorem \ref{teorema1}, we have that 
		\begin{equation*}
			\mathbb{P}\left(\widehat{S}_n=S_n\right)\rightarrow 1 \textrm{ as } n\rightarrow\infty.
		\end{equation*}
	\end{theorem}
	
	The following result focuses on the asymptotic distribution of certain projections of $\widehat{\pmb{\beta}}_0$. First, let us introduce some additional notation. We will denote
	$$\pmb{c}=\left(c_1,\dots,c_{p_n}\right)^{\top}, \textrm{ being } c_j=\mathcal{P}_{\lambda_{jn}}'\left(\left|\beta_{0j}\right|\right)\textrm{sgn}(\beta_{0j})1_{\{j\in S_n\}},$$
	and
	$$\pmb{V}=\textrm{diag}\left\{V_1,\dots,V_{p_n}\right\},\textrm{ where } V_j=\mathcal{P}_{\lambda_{jn}}''\left(\left|\beta_{0j}\right|\right)1_{\{j\in S_n\}}.$$
	In addition, we will denote $\sigma_{\varepsilon}^2=\mathbb{E}(\varepsilon_i^2)$, while
	$\pmb{A}_n$ will be any $q\times s_n$ matrix such that $\pmb{A}_n\pmb{A}_n^{\top}\rightarrow \pmb{A}$ as $n\rightarrow\infty$, where $\pmb{A}$ is a $q\times q$ definite positive matrix.
	
	\begin{theorem}[Asymptotic normality]\label{teorema3} 
		Adding the following conditions to assumptions in Theorem \ref{teorema1} (where, if $s_n=1$, $\log s_n$ must be interpreted as $1$):
		\begin{align}
			&\exists \beta'>1 \textrm{ such that } s_n\exp\left\{\left(1-\beta'\log s_n\right)\psi_{\mathcal{C}}\left(\frac{1}{n}\right)\right\}\rightarrow 0 \textrm{ as } n\rightarrow \infty,\nonumber\\	
			&\max\left\{ns_n^3h^{4\alpha},n^{2/r_{\varepsilon}}s_n h^{2\alpha}\log^2n,s_n^3h^{2\alpha}\log^2n,n^{-1}s_n^3,ns_n^3\delta_n^4\right\}=o(1)\nonumber\\
			&and \nonumber\\
			&\max\left\{n^{2/r_{\varepsilon}}s_n\log s_n\log^2 n,s_n^3\log s_n\log^2n\right\}=o\left(\frac{nf(h)}{\psi_{\mathcal{C}}\left(1/n\right)}\right),\nonumber
		\end{align}
		the following result can be established:
		\begin{equation*}
			n^{1/2}\pmb{A}_n\pmb{C}_{\theta_0,S_n}\left(\widehat{\pmb{\beta}}_{0S_n}-\pmb{\beta}_{0S_n}+\left(\pmb{B}_{\theta_0S_n\times S_n}+\pmb{V}_{S_n\times S_n}\right)^{-1}\pmb{c}_{S_n}\right)\overset{\textnormal{d}}{\longrightarrow}\textrm{N}(\pmb{0},\pmb{A}),
		\end{equation*}	
		where we have denoted $\pmb{C}_{\theta_0,S_n}=\sigma_{\varepsilon}^{-1}\pmb{B}_{\theta_0S_n\times S_n}^{-1/2}\left(\pmb{B}_{\theta_0S_n\times S_n}+\pmb{V}_{S_n\times S_n}\right)$.
	\end{theorem}
	
	\begin{remark}
		\label{r_theorems3y3}
		Theorems \ref{teorema2} and \ref{teorema3} show that $\widehat{\pmb{\beta}}_0$ enjoys the oracle property in the meaning given, for instance, in \cite{xie_2009}: ``the estimator can correctly select the nonzero coefficients with probability converging to one, and that the estimators of the nonzero coefficients are asymptotically normal with the same means and covariances that they would have if the zero coefficients were known in advance". In the setting of multivariate regression (not functional), the interested reader can find estimators verifying such property in \cite{fanpeng_2004} (linear regression), \cite{xie_2009} (partially linear regression) or \cite{wangzhu_2017} (partial linear single-index regression), among others. See also \cite{aneiros_2015} for the case of the semi-functional partial linear regression. 
	\end{remark}
	
	Finally, the next theorem states the uniform rate of convergence of the statistic $\widehat{m}_{\theta_0}\left(\chi\right)$ in (\ref{m_statistic}).
	\begin{theorem}\label{teorema4}
		Under assumptions of Theorem \ref{teorema1}, if in addition the following conditions are verified:
		\begin{enumerate}
			
			\item[A)] $\forall (\chi,\chi')\in\mathcal{C}\times\mathcal{C}, \textrm{ } \left|m_{\theta_0}(\chi)-m_{\theta_0}(\chi')\right|\leq C_{13}d_{\theta_0}\left(\chi,\chi'\right)^{\alpha}$, where $\alpha$ was defined in (\ref{smooth}),
			
			\item[B)] $\sup_{\chi\in\mathcal{C},j\in S_n}\left|g_{j,\theta_0}(\chi)\right|=O(1)$
			
			and
			
			\item[C)] $\psi_{\mathcal{C}}\left(1/n\right)\rightarrow\infty$ as $n\rightarrow \infty$,
		\end{enumerate}
		then we have that
		\begin{equation*}
			\sup_{\theta\in\Theta_n}\sup_{\chi\in\mathcal{C}}\left|\widehat{m}_{\theta}\left(\chi\right)-m_{\theta_0}\left(\chi\right)\right|=O_p\left(h^{\alpha}+\sqrt{\frac{\psi_{\mathcal{C}}\left(1/n\right)}{nf(h)}}\right)+O_p\left(\sqrt{s_n}u_n\right).
		\end{equation*}
	\end{theorem}
	
	\begin{corollary}\label{cor}
		Under assumptions of Theorem \ref{teorema4}, we have that
		\begin{equation*}
			\sup_{\chi\in\mathcal{C}}\left|\widehat{m}_{\widehat{\theta}_0}\left(\chi\right)-m_{\theta_0}\left(\chi\right)\right|=O_p\left(h^{\alpha}+\sqrt{\frac{\psi_{\mathcal{C}}\left(1/n\right)}{nf(h)}}\right)+O_p\left(\sqrt{s_n}u_n\right).
		\end{equation*}
	\end{corollary}

	\begin{corollary}\label{final-cor}
		Under assumptions of Theorem \ref{teorema4}, if in addition the following conditions hold:
		\begin{enumerate}
			
			\item[A)] $\forall \theta \in \Theta_n$, the random variables $\left<\theta,\mathcal{X}\right>$ are valued in the same compact subset, $\mathcal{R}$, of $\mathbb{R}$, and are absolutely continuous with respect to the Lebesgue measure, with density $f_\theta$ satisfying
			\begin{equation*}
				0<\inf_{\theta\in\Theta_n, \ u \in \mathcal{R}} f_\theta(u) \leq \sup_{\theta\in\Theta_n, \ u \in \mathcal{R}}f_\theta(u) < \infty,
			\end{equation*}
			\item[B)] $h \approx C(\log n / n)^{1/(2\alpha+1)}$,
			\item[C)] $s_n\approx cn^\gamma$ with $0<2\gamma<1-2\alpha/(2\alpha+1)$
			
			and
			
			\item[D)] $\delta_n=O(n^{-1/2})$ ($\delta_n$ was defined in (\ref{not-4.3})),
		\end{enumerate}
		then we have that
		\begin{equation*}
			\sup_{\theta\in\Theta_n}\sup_{\chi\in\mathcal{C}}\left|\widehat{m}_{\theta}\left(\chi\right)-m_{\theta_0}\left(\chi\right)\right|=O_p\left( \left(\frac{\log n}{n}\right)^{\alpha/(2\alpha + 1)}\right)
		\end{equation*}
		and
		\begin{equation*}
			\sup_{\chi\in\mathcal{C}}\left|\widehat{m}_{\widehat{\theta}_0}\left(\chi\right)-m_{\theta_0}\left(\chi\right)\right|=O_p\left(\left(\frac{\log n}{n}\right)^{\alpha/(2\alpha + 1)}\right).
		\end{equation*}
	\end{corollary}
	
	\begin{remark}
		\label{r_theorem4}
		Theorem \ref{teorema4} extends, in the same sense as in Remark \ref{r_theorem1}, Theorem 3.3 in \cite{aneiros_2015} from the case $\Theta_n=\{\theta_0\}$ (i.e., $v_n=0$ in (\ref{Theta}); equivalently, $\theta_0$ known) to the case where $\{\theta_0\} \subset \Theta_n$ (i.e., $\theta_0$ unknown). Corollary \ref{final-cor} shows a nice property of dimensionality reduction: the semi-functional nonparametric component $m_{\theta_0}(\cdot)\equiv m(\left<\theta_0,\cdot \right>)$ is estimated with univariate nonparametric rate (note that the Condition A imposed in Corollary \ref{final-cor} was used, for instance, in \citealt*{ferraty_2013}, while the Condition D is satisfied, for instance, for the SCAD penalty function; finally, Condition B considers a bandwidth with optimal rate for univariate nonparametric regression while Condition C is a non-restrictive technical condition). 
	\end{remark}

	\section{Simulation study}\label{MC}
	\label{sim-study}
	The aim of this section is to show the finite sample behaviour of the two statistical procedures presented before for the  SSFPLSIM (\ref{modelo}); that is, firstly the  penalized least-squares procedure (for both variable selection and estimation of the linear parameters $\pmb{\beta}_0$) and secondly the single-index approach for estimating the functional semiparametric component of the model.
	
	\subsection{Design}
	\label{design-sim}
\begin{figure}[h]
		\centering
		\includegraphics[width=6cm]{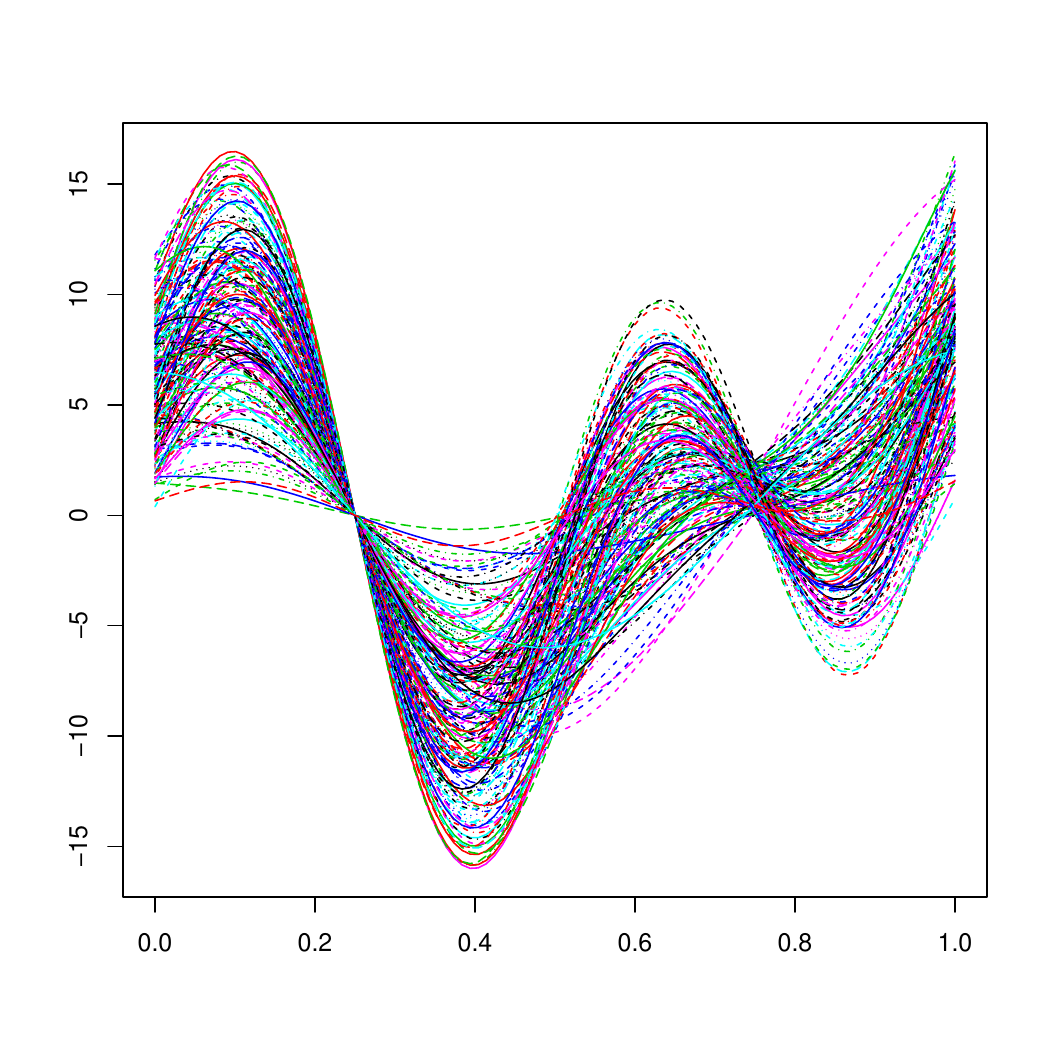}%\hspace{1.5cm}
		\includegraphics[width=6cm]{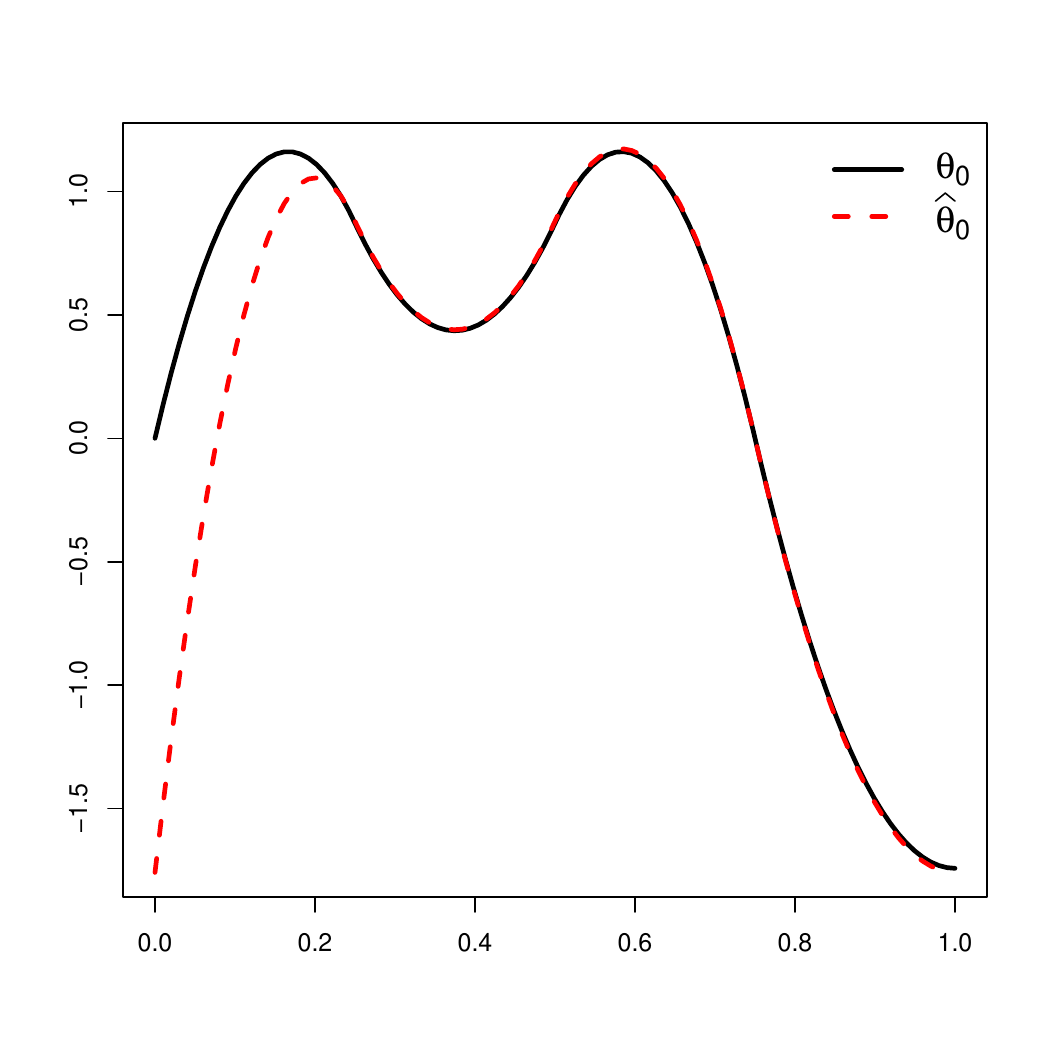}
		\caption{Sample of 200 curves generated from (\ref{X-sim}) (left panel) and  functional direction $\theta_0$ (right panel). In addition, in right panel, it is displayed the estimation, $\widehat{\theta}_0$, of $\theta_0$ obtained from a particular sample in the scenario $(n,p_n,\rho,c)=(100,50,0,0.05)$.}
		\label{fig1}
	\end{figure}

	For $(n,p_n) \in \{(100,50),(200,100)\}$, samples of iid data $\{(X_{i1},\dots,X_{ip_n},\mathcal{X}_i,Y_i)\}_{i=1}^{n}$ were constructed according to the following model:
	\begin{equation}
		%\label{modelo}
		Y_i=X_{i1}\beta_{01}+\dots+X_{ip_n}\beta_{0p_n}+m\left(\left<\theta_0,\mathcal{X}_i\right>\right)+\varepsilon_i, \ \forall i=1,\dots,n,\label{mod_sim}
	\end{equation}
	where: \begin{itemize}
		\item The vectors of real covariates, $(X_{i1},\dots,X_{ip_n})^{\top}$ ($i=1,\ldots,n$), were generated from a multivariate normal distribution with zero mean and covariance matrix given by $(\rho^{|j-k|})_{jk}$ ($j,k=1,\dots,p_n$). Two values for $\rho$ (namely $\rho=0$ and $\rho=0.5$) were considered.
		\item The functional covariate, $\mathcal{X}_i$ ($i=1,\ldots,n$), was generated in the following way:
		\begin{equation}
			\mathcal{X}_i(t)=a_i\cos(2\pi t) + b_i\sin(4\pi t) + 2c_i(t-0.25)(t-0.5) \ \forall t \in [0,1],
			\label{X-sim}
		\end{equation}
		where the random variables $a_i,b_i$ and $c_i$ ($i=1,\ldots,n$) were independent and uniformly distributed on the interval $[0,10]$ (note that we refer to independence both between and within vectors $(a_i,b_i,c_i)^\top$). These curves were discretized on the same grid of 100 equispaced points in $[0,1]$.  
		\item The iid random errors, $\varepsilon_i$ ($i=1,\ldots,n$), were simulated from a $N(0,\sigma_{\varepsilon})$ distribution, where $\sigma_{\varepsilon}^2=c\sigma_r^2$ with $\sigma_r^2$ denoting the empirical variance of the regression $X_{i1}\beta_{01}+\dots+X_{ip_n}\beta_{0p_n}+m\left(\left<\theta_0,\mathcal{X}_i\right>\right)$. Note that $c$ is the signal-to-noise ratio, and two values (namely $c=0.01$ and $c=0.05$) were considered.
		\item The true vector of linear coefficients was $$\pmb{\beta}_0=(\beta_{01},\dots,\beta_{0p_n})^{\top}=(3,1.5,0,0,2,0,\dots,0)^{\top}.$$ 
		\item The true direction of projection was 
		\begin{equation}
			\theta_0(\cdot)=\sum_{j=1}^{d_n}\alpha_{0j}e_j(\cdot), 
			\label{theta_0_base}
		\end{equation}
		where $\{e_1(\cdot),\ldots,e_{d_n}(\cdot)\}$ is a set of B-spline basis functions and $d_n=l+m_n$ ($l$ denotes the order of the splines and $m_n$ is the number of regularly interior knots). Values $l=3$ and $m_n=3$ were considered, while the vector of coefficients of $\theta_0$ in expression (\ref{theta_0_base}) was $$(\alpha_{01},\dots,\alpha_{0d_n})^\top=(0,1.741539,0,1.741539,-1.741539,-1.741539)^\top$$ (note that this vector is obtained by calibrating the one $(0,1,0,1,-1,-1)^\top$ to ensure identifiability; for details, see \citealt*{novo}).
		\item  The inner product and the link function were $\left<f,g\right>=\int_{0}^1 f(t)g(t)dt$ and $m(\left<\theta_0,\chi\right>)=\left<\theta_0,\chi\right>^3$, respectively.
	\end{itemize}
	Figure \ref{fig1} (left panel) shows a sample of $200$ curves generated from (\ref{X-sim}), while in its right panel, in black colour and solid line, the functional direction $\theta_0$ is displayed.
	
	For each simulation case $((n,p_n),\rho,c)\in\{(100,50),(200,100)\}\times\{0,0.5\}\times\{0.01,0.05\}$, $M=100$ independent samples were generated from model (\ref{mod_sim}). For each sample, we obtained an estimator of the pair $(\pmb{\beta}_0,\theta_0)$ by minimizing the penalized profile least-squares function $\mathcal{Q}\left(\pmb{\beta},\theta\right)$ (see (\ref{func_minimizar})). For that, we considered as eligible functional directions
	\begin{equation*}
		\theta(\cdot)=\sum_{j=1}^{d_n}\alpha_je_j(\cdot)
		%\label{theta_base}
	\end{equation*}
	for a wide set of vectors of coefficients, $(\alpha_1,\dots,\alpha_{d_n})^\top$, constructed  following the procedure described in \cite{novo}. Epanechnikov kernel was used while the penalty function considered was the SCAD one with parameter $a=3.7$ (see \citealt*{Fan97}, \citealt*{fanli_2001} or \citealt*{aneiros_2015}). To reduce  the quantity of tuning parameters, $\lambda_j$, to be selected for each sample, we consider $\lambda_j = \lambda \widehat{\sigma}_{\beta_{0,j,OLS}}$, where $\beta_{0,j,OLS}$ denotes the OLS estimate of $\beta_{0,j}$ in (\ref{mod_sim}) and $\widehat{\sigma}_{\beta_{0,j,OLS}}$ is the estimated standard deviation. This tuning parameter, $\lambda$, as well as the bandwidth, $h$, were selected by means of the BIC procedure. More specifically, the BIC value corresponding to $(\widehat{\pmb{\beta}}_{0,h,\lambda},\widehat{{\theta}}_{0,h,\lambda})$ (the estimate of the parameter $(\pmb{\beta}_0,\theta_0)$ in the linear model (\ref{linear}) obtained by minimizing the profile least-squares function (\ref{func_minimizar})) was computed from the routine {\tt select} of the R package \textit{grpreg}. The main reason why we have used this selector is its low computational cost compared to cross-validation-based selectors (which are time consuming procedures). Taking into account that this BIC selector shows good behaviour both in this simulation study and in the real data application reported in Section \ref{real-data}, it is noteworthy that the feature of its low computational cost takes a main relevance in a so complex model as SSFPLSIM. 
	
	\subsection{Results}
\begin{figure}[h]
		\centering
		\includegraphics[width=5.8cm]{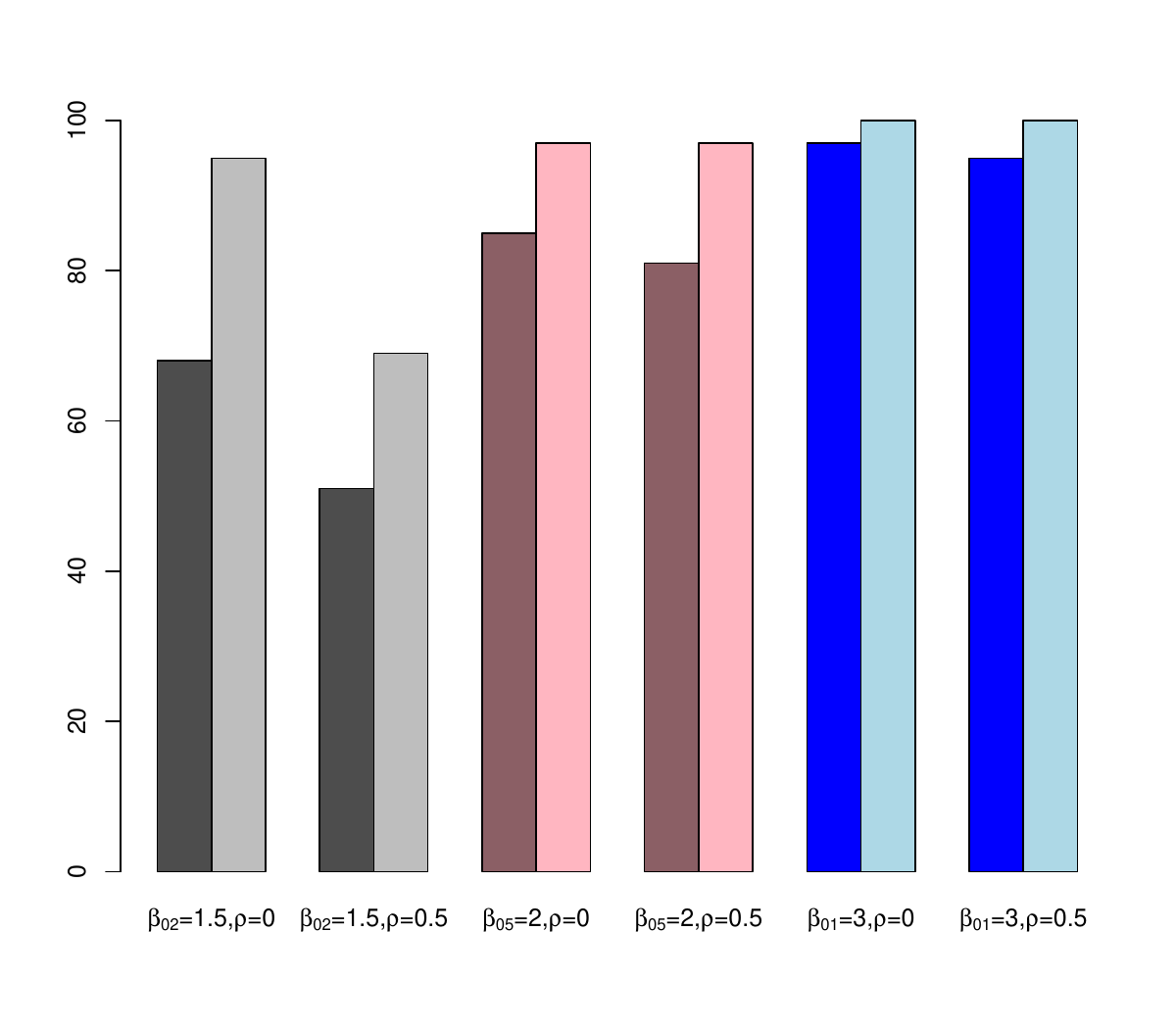}
		\includegraphics[width=5.8cm]{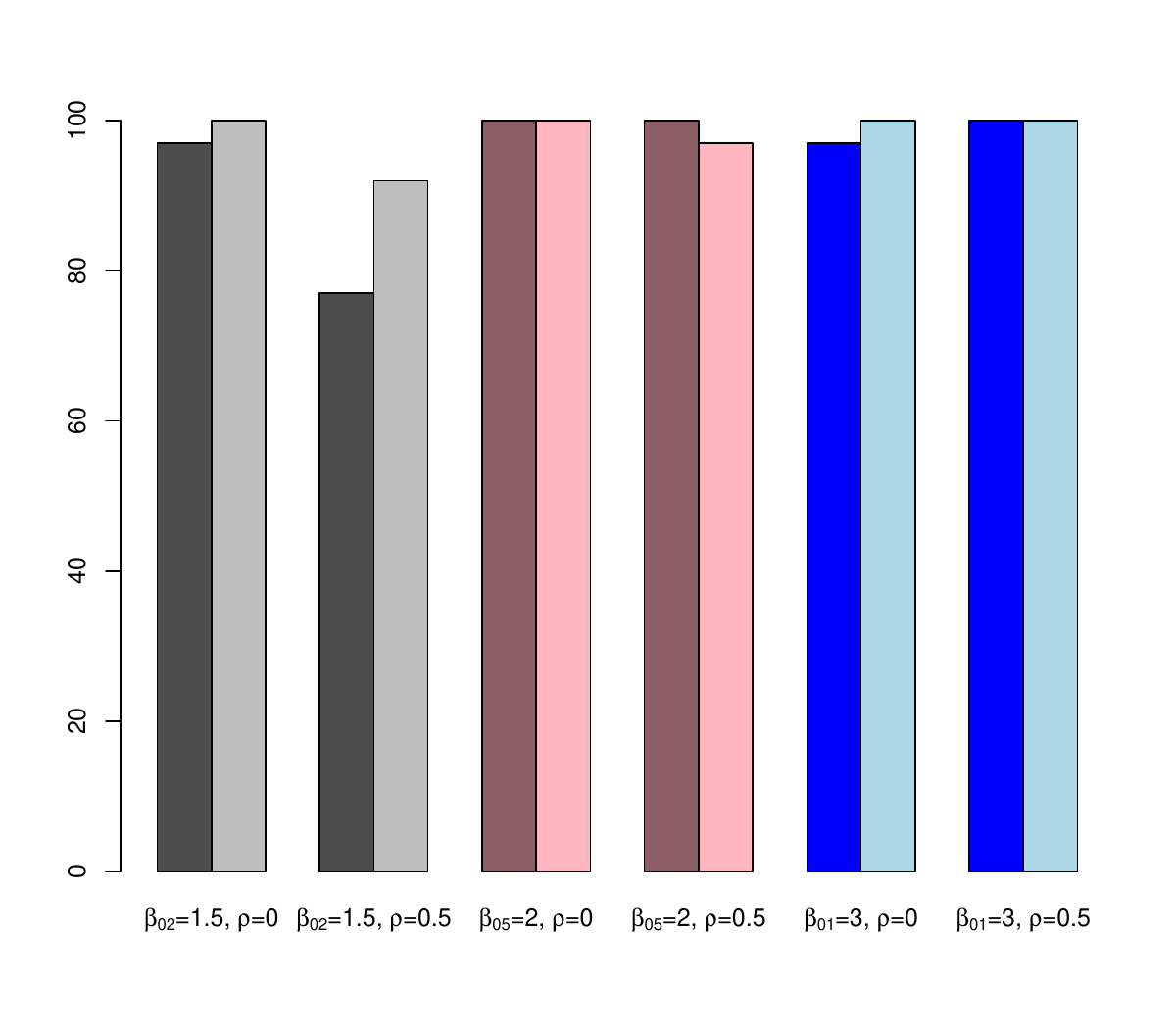}
		\caption{Percentage of times that each non-zero coefficient of $\pmb{\beta}_0$ is not set to zero (left panel: $c=0.05$; right panel: $c=0.01$).  We use grey for $\beta_{02}=1.5$, pink for $\beta_{05}=2$ and blue for $\beta_{01}=3$. Dark colours correspond to $n=100$ while light colours match $n=200$.  Values $\rho=0$ and $\rho=0.5$ are considered.}
		\label{fig2}
	\end{figure}

	First results of the simulation study are presented in Table \ref{tabla_ic_bic} and Figure \ref{fig2} (variable selection) and Table \ref{suma_errores_2_BIC} ($\pmb{\beta}_0$ estimation).
	
	Table \ref{tabla_ic_bic} shows both the average percentage (restricted only to the true zero coefficients) of coefficients correctly set to zero and the average percentage (restricted only to the true non-zero coefficients) of coefficients erroneously set to zero. 
	
	From Table \ref{tabla_ic_bic} we can observe that, as sample size increases, our procedure is able to detect a greater percentage of non-significant variables. In addition,  the percentage of erroneously set to zero significant variables decreases. It is noteworthy that positive dependence between variables gives some advantage in detecting non-significant variables, but it is detrimental to detection of the significant ones (similar conclusions were obtained in nonfunctional both linear (\citealt*{HuMaZh08}) and partial linear (\citealt*{xie_2009} models), as well in the semi-functional partial linear model (\citealt*{aneiros_2015})). We can also observe that results are better for $c=0.01$ than $c=0.05$, especially for finding the true relevant variables. 
	
	\begin{table}[h]%[ht]
		\centering
		%\scalebox{0.85}{
		\begin{tabular}{|ccr|rr|rr|}
			\hline
			&&&\multicolumn{2}{c|}{$\rho=0$}&\multicolumn{2}{c|}{$\rho=0.5$}\\
			\hline
			$n$& $p_n$&	$c$	& Correct & Incorrect & Correct & Incorrect \\ 
			\hline	
			\multirow{2}{1cm}{$100$} &	\multirow{2}{1cm}{$50$}	& 0.05  & 77.404 & 16.667 & 84.447 & 24.333 \\ 
			
			&	& 0.01 & 92.830 & 1.000 & 96.319 & 7.667 \\ 
			
			\hline
			\multirow{2}{1cm}{$200$ }&\multirow{2}{1cm}{$100$}	&0.05 & 85.052 & 2.667 & 91.072 & 11.333 \\ 
			
			&	&0.01 &  98.619 & 0.000 & 99.732 & 2.667 \\ 
			
			\hline
		\end{tabular}%}
		\caption{Column ``Correct'': Average percentage (restricted only to the true zero coefficients) of coefficients correctly set to zero. Column ``Incorrect'': Average percentage (restricted only to the true non-zero coefficients) of coefficients erroneously set to zero.}
		\label{tabla_ic_bic}
	\end{table}
	
	Figure \ref{fig2} shows barplots with the percentages of time that each non-zero coefficient ($\beta_{01}=3$, $\beta_{02}=1.5$  and $\beta_{05}=2$)  is not set to zero. Therefore, since the linear covariates are identically distributed, Figure \ref{fig2} allows to analyse the influence of the size of each $\beta_{0j}$ non-zero coefficient in the detection of the $j^{th}$ variable as influential one. A first conclusion is that, as intuition says, as bigger is the value of  $\beta_{0j}$, greater is the percentage of success. In general, results also improve if we increase the sample size or if we reduce $c$ (and then, $\sigma^2_{\varepsilon}$). In addition, positive dependence between variables makes more difficult the detection of the significant variables, especially for smaller values of  $\beta_{0j}$.
	
	Table \ref{suma_errores_2_BIC} reports information about the performance of the penalized least-squares (PLS) estimator of $\pmb{\beta}_0$ in the SSFPLSIM (\ref{mod_sim}). More specifically, on the one hand it shows both the mean and standard deviation of the squared errors, 
	\begin{equation}
		||\pmb{\widehat{\beta}}_0-\pmb{\beta}_0||^2=\sum_{j=1}^{p_n}(\widehat{\beta}_{0j} - \beta_{0j})^2,
		\label{error-beta0}
	\end{equation}
	obtained from the $M$ replicates when both the proposed PLS approach and the ordinary least-squares (OLS) estimator are applied to the SSFPLSIM. On the other hand, it reports the corresponding both mean and standard deviation of the squared errors obtained from the OLS approach assuming that one knows in advance what are the non-null coefficients (oracle estimator); note that the oracle estimator can't be used in practice, but can be seen as a benchmark method in simulation. As expected, the oracle estimator performs the best and the OLS one performs the worst. The proposed PLS estimator presents a good behaviour, being its performance much closer to that of the oracle estimator than to the OLS one. In addition, its performance significantly improves when the sample size increases (note that the estimation of $\theta_0$ is needed, which, given the complexity of the SSFPLSIM, requires a sufficient sample size) or the signal-to-noise ratio ($c$) decreases. Note also that the dependence in the linear covariates has effects on both the estimators (oracle and OLS) and the variable selection procedure (PLS): in general, small values of $\rho$ provide better results.
	
	\begin{table}[h]%[ht]
		\centering
		%\scalebox{0.85}{
		\begin{tabular}{|ccc|rr|rr|rr|}
			\hline
			&	&	&\multicolumn{2}{c|}{ORACLE}&\multicolumn{2}{c|}{PLS}&\multicolumn{2}{c|}{OLS}\\
			
			\hline
			$c$& $\rho$	&$n$	& Mean & SD & Mean & SD  & Mean & SD \\ 
			\hline
			\multirow{4}{1cm}{$0.05$}	&\multirow{2}{1cm}{$0$}	& 100 & 1.879 & 1.855 & 12.742 & 13.241 & 89.209 & 51.619 \\

			&& 200	 & 0.796 & 0.628 & 4.055 & 2.807 & 61.023 & 26.082 \\ 
			
			&\multirow{2}{1cm}{$0.5$}		&100& 2.347 & 2.882 & 11.014 & 9.610 & 148.431 & 91.162 \\ 
			&	& 200	 & 1.051 & 0.956 & 4.804 & 3.821 & 99.872 & 41.521 \\

			\hline
			\multirow{4}{1cm}{$0.01$}	&\multirow{2}{1cm}{$0$}& 100	& 0.440 & 0.502 & 1.359 & 1.063 & 22.425 & 17.874 \\ 
			
			&	& 200	 & 0.175 & 0.140 & 0.656 & 0.519 & 12.847 & 4.738 \\
			
			&\multirow{2}{1cm}{$0.5$}	& 100	 & 0.540 & 0.733 & 2.645 & 1.720 & 37.210 & 30.573 \\
			&	&	200 & 0.231 & 0.207 & 1.825 & 1.080  & 21.067 & 7.868 \\ 
			\hline
		\end{tabular}%}
		\caption{Mean and standard deviation (SD) of the squared errors (\ref{error-beta0}) obtained from ORACLE, (the proposed) PLS and OLS procedures.}
		\label{suma_errores_2_BIC}
	\end{table}
	
	The remainder of this section focuses on the estimation of the semiparametric component of the model; that is, on both the single index, $\theta_0$, and the nonparametric component, $m(\cdot)$. Table \ref{suma_errores_theta0_BIC} shows the performance of the proposed estimator of $\theta_0$ in the SSFPLSIM (\ref{mod_sim}). Such performance is measured by means of the mean of the squared errors,
	\begin{equation}
		||\widehat{\theta}_0-\theta_0||^2=\int_0^1\left(\widehat{\theta}_{0}(t)-\theta_0(t)\right)^2dt ,
		\label{error-theta0}
	\end{equation}
	obtained from the $M$ replicates.
	
	\begin{table}[ht]
		\centering
		\begin{tabular}{|ccc|rr|}
			\hline
			c& $\rho$	&$n$	& Mean & SD \\ 
			\hline
			\multirow{4}{1cm}{$0.05$}	&\multirow{2}{1cm}{$0$}	& 100 & 0.010 & 0.011 \\ 
			&& 200 & 0.003 & 0.007 \\ 
			&\multirow{2}{1cm}{$0.5$}		&100& 0.009 & 0.011 \\ 
			&&200& 0.004 & 0.004 \\ 
			
			\hline
			
			\multirow{4}{1cm}{$0.01$}	&\multirow{2}{1cm}{$0$}& 100	& 0.004 & 0.008 \\ 
			&	& 200	& 0.001 & 0.004 \\ 
			&\multirow{2}{1cm}{$0.5$}	& 100	 & 0.004 & 0.007 \\ 
			&	&	200 & 0.001 & 0.004 \\ 
			\hline
		\end{tabular}
		\caption{Mean and standard deviation (SD) of the squared errors (\ref{error-theta0}) obtained from the proposed procedure.}
		\label{suma_errores_theta0_BIC}
	\end{table}
	\noindent From Table \ref{suma_errores_theta0_BIC} we can conclude that the performance of the proposed estimator for the single-index (infinite-dimensional) parameter, $\theta_0$, clearly improves when the sample size increases or the signal-to-noise ratio ($c$) decreases (in a similar way as happened for estimation of the linear (finite-dimensional) parameter $\pmb{\beta}_0$). In addition, it appears that, once the variable selection is performed, the estimation errors are not really affected by higher values of $\rho$. 
	
\begin{figure}[h]
		\centering
		\includegraphics[width=5.8cm]{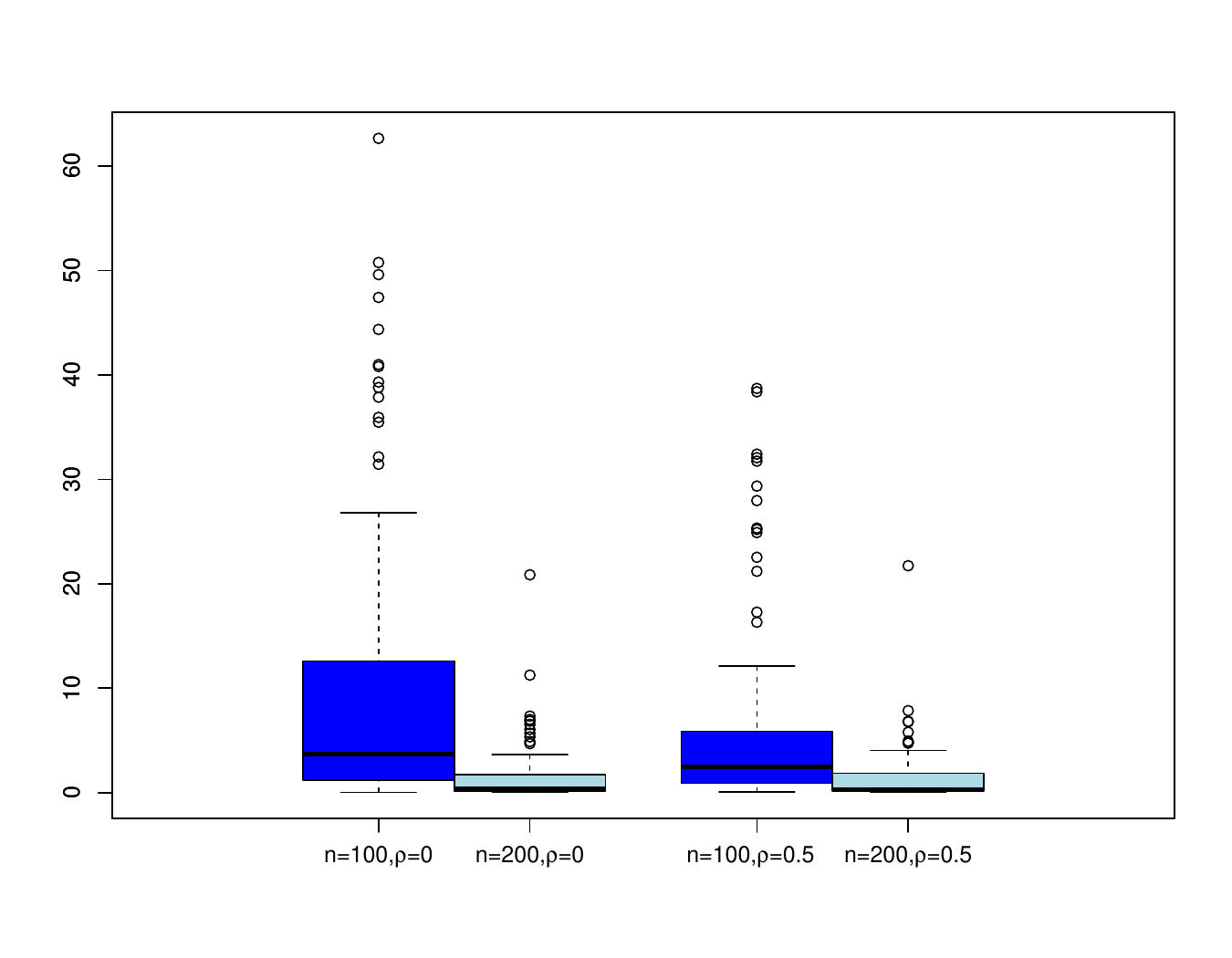}
		\includegraphics[width=5.8cm]{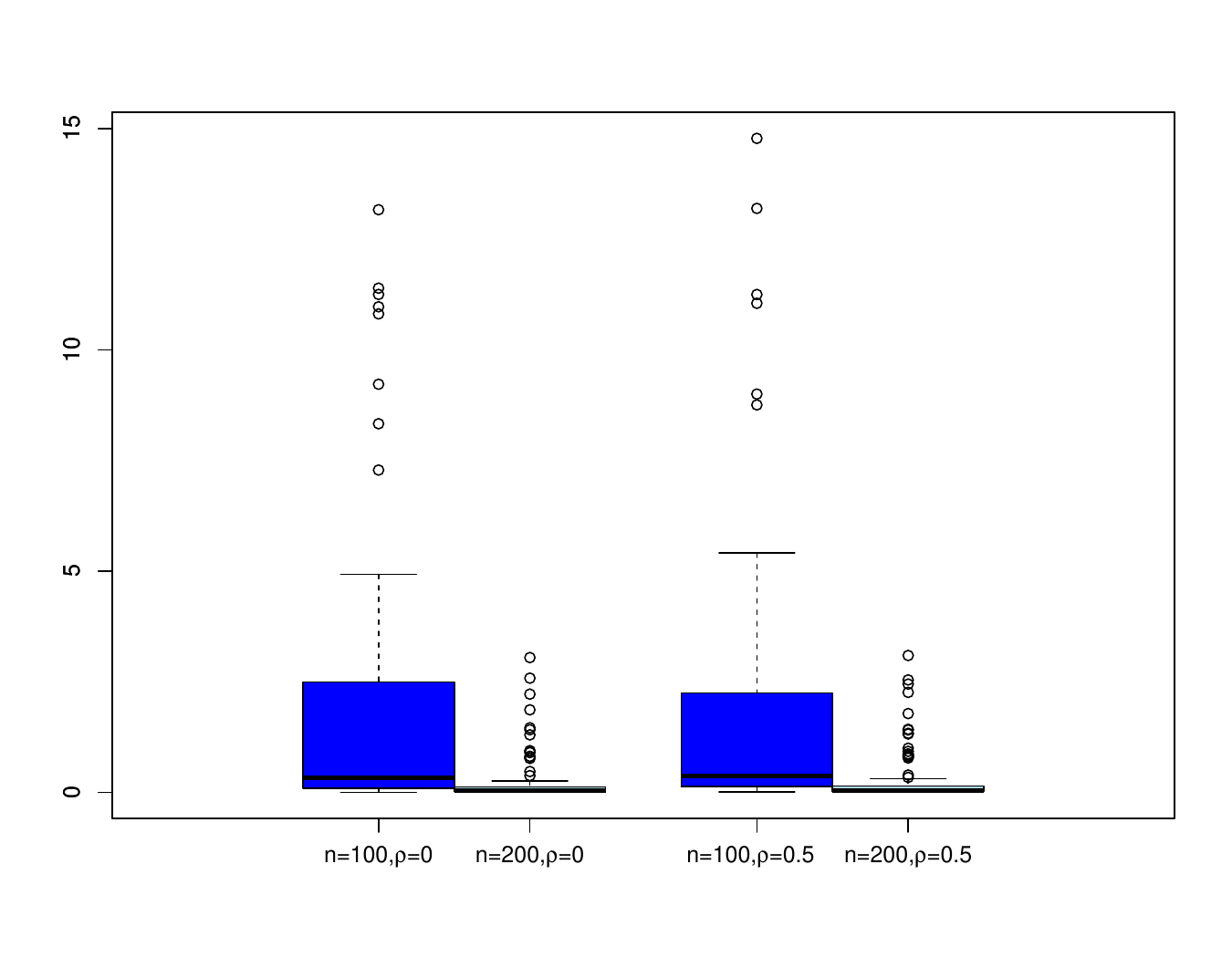}
		\caption{Boxplots of the squared errors (\ref{error-m}) obtained from the proposed procedure for the several considered scenarios. Left panel: $c=0.05$; right panel: $c=0.01$.}
		\label{fig3}
	\end{figure}

	In order to measure the performance of the proposed estimator ($\widehat{m}(\cdot)$) for the nonparametric component ($m(\cdot)$), $M$ independent test samples with sample size $N=100$,
	%$$\left\{(X_{j1}^{(k)},\dots,X_{jp_n}^{(k)},\mathcal{X}_j^{(k)},Y_j^{(k)})\right\}_{j=1}^{N}, \ k=1,\ldots,M,$$
	$$\left\{\mathcal{X}_j^{(k)}\right\}_{j=1}^{N}, \ k=1,\ldots,M,$$
	were constructed is a similar way as in Section \ref{design-sim} (note that these $M$ test samples were also independent of the $M$ (training) samples considered until now). Then, the performance of the estimate for $m(\cdot)$ constructed from the $k$-th training sample was measured by means of the Mean Squared Error of Prediction (MSEP),
	\begin{equation}
		MSEP_k=\frac{1}{N}\sum_{j=1}^N\left(\widehat{m}_k\left(\left<\widehat{\theta}_{0k},\mathcal{X}_j^{(k)}\right>\right)-m\left(\left<\theta_0,\mathcal{X}_j^{(k)}\right>\right)\right)^2,
		\label{error-m}
	\end{equation}
	where $\widehat{\theta}_{0k}$ and $\widehat{m}_k(\cdot)$ denote estimators for $\theta_0$ and $m(\cdot)$, respectively, constructed from information in the $k$-th training sample. Figure \ref{fig3} displays, for each considered scenario, boxplots with the corresponding $MSEP_k$ values.
	
	\noindent In a similar way as for estimation of both the linear (finite-dimensional) parameter $\pmb{\beta}_0$ (see Table \ref{suma_errores_2_BIC}) and the single-index (infinite-dimensional) parameter $\theta_0$ (see Table \ref{suma_errores_theta0_BIC}), Figure \ref{fig3} shows that the performance of the estimate of the nonparametric component, $\widehat{m}(\cdot)$, clearly improves when the sample size increases or the signal-to-noise ratio ($c$) decreases.
	
\begin{figure}[h]
		\centering
		\includegraphics[width=6cm]{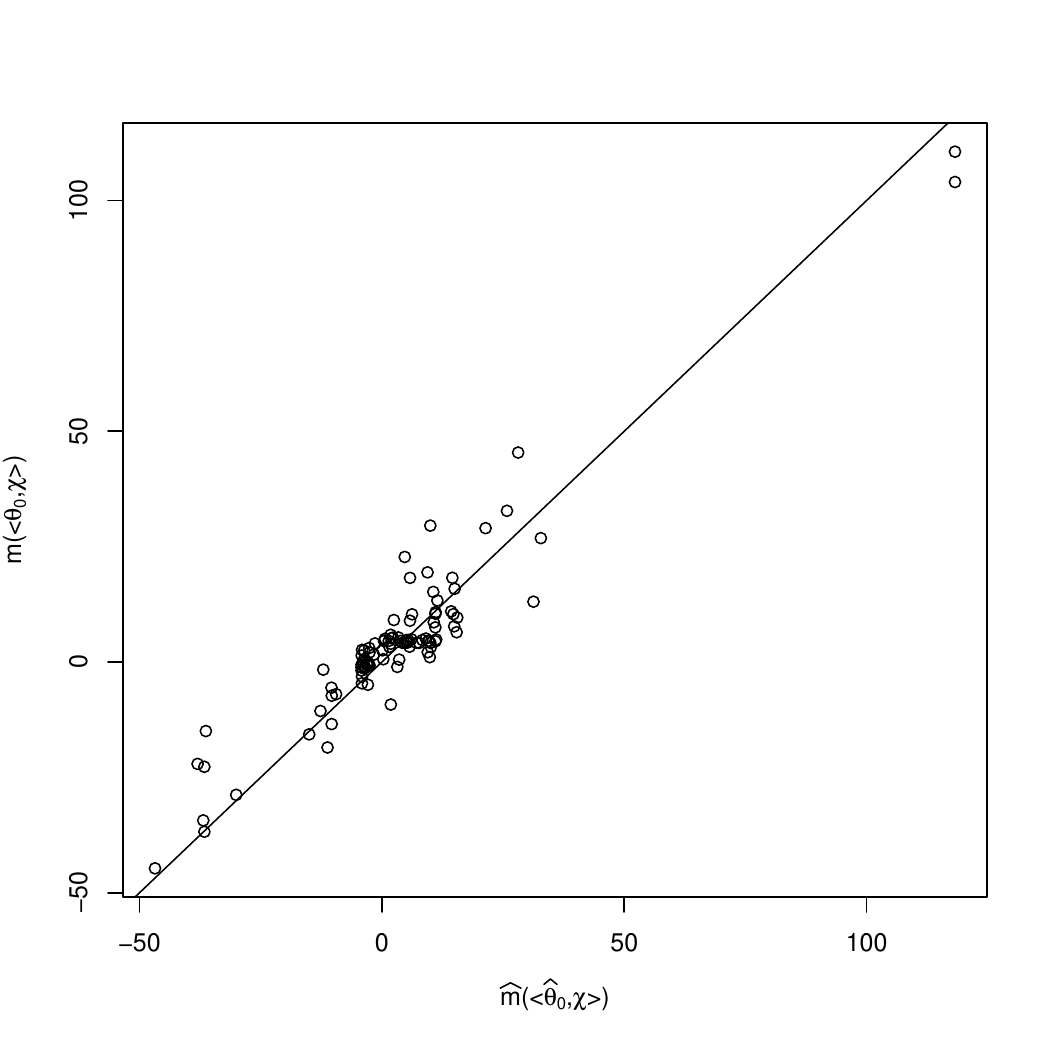}%\hspace{1.5cm}
		\includegraphics[width=6cm]{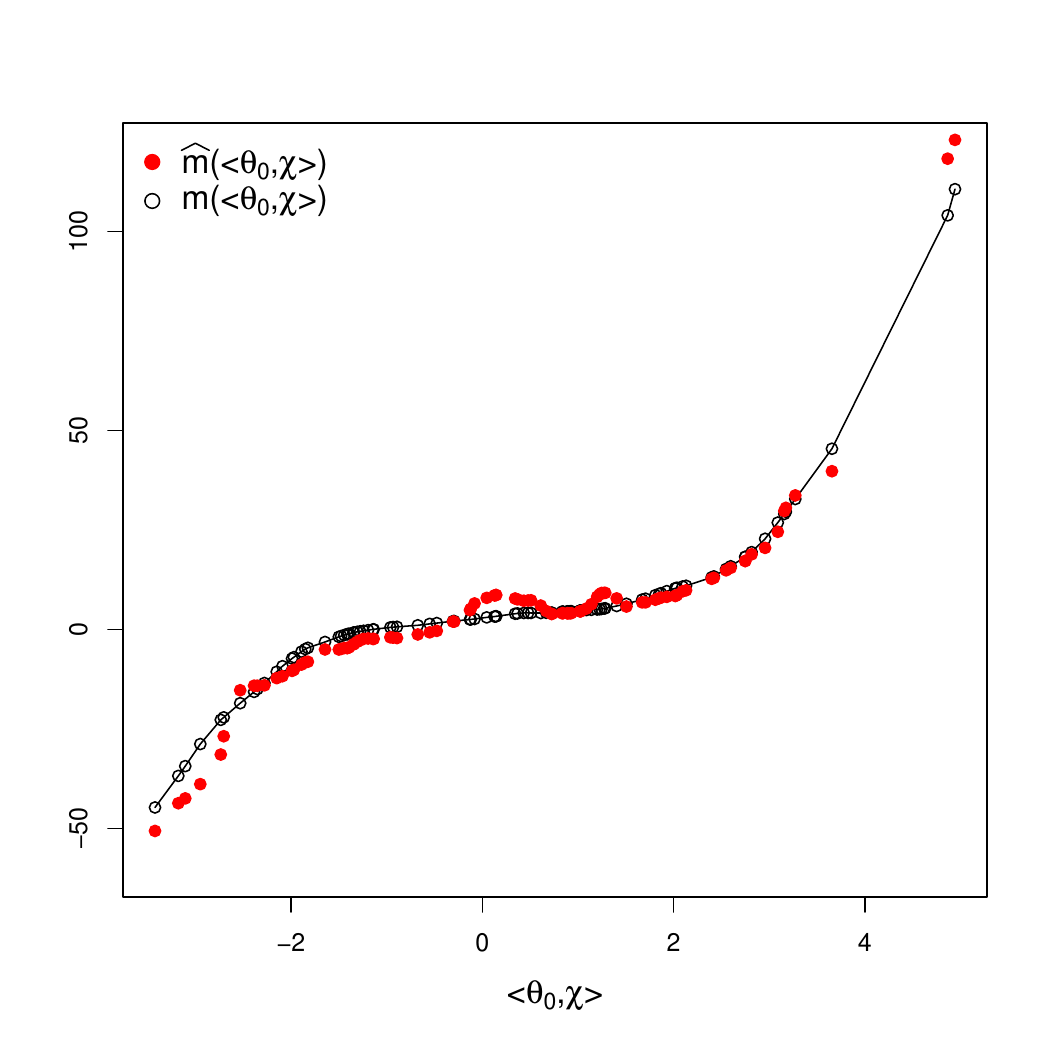}
		\caption{Real and estimated values, from a particular sample in the scenario $(n,p_n,\rho,c)=(100,50,0,0.05)$, related to the semiparametric component, $m\left(\left< \theta_0,\cdot\right>\right)$, of the SSFPLSIM (\ref{mod_sim}). The curve in the right panel is the true $m$.}
		\label{fig4}
	\end{figure}

	Finally, Figure \ref{fig4} displays, for a particular replicate, values of $m\left(\left< \theta_0,\cdot\right>\right)$ vs. $\widehat{m}\left(\left< \widehat{\theta}_0,\cdot\right>\right)$, as well as both $m\left(\left< \theta_0,\cdot\right>\right)$ and $\widehat{m}\left(\left< \theta_0,\cdot\right>\right)$ vs. $\left< \theta_0,\cdot\right>$. For a graphic showing the estimate of $\theta_0$ obtained from such particular replicate, see right panel in Figure \ref{fig1} (red color and dashed line).

	\section{Real data application}
	\label{real-data}
	In this section, a benchmark data set in the nonparametric functional context is modelled thought different functional regression models, including the SSFPLSIM (\ref{modelo}) proposed in this paper. The results obtained show the usefulness of both the SSFPLSIM and the proposed PLS estimation procedure.
	
	Before beginning the next sections dedicated to present the data set, modelling, variable selection and prediction, we indicate that, in the estimation of the three models that require variable selection (see models SLM, SFPLM and SSFPLSIM in Table \ref{table5}, Section \ref{mod-pred}), both the tuning parameter, $\lambda$, and the bandwidth, $h$, were selected by means of the BIC procedure, and the Epanechnikov kernel and the penalty function SCAD (with parameter $a = 3.7$) were used (in a similar way as in the simulation study in Section \ref{sim-study}). In addition, in the SSFPLSIM the order of the splines was $l=3$ while the number of regularly interior knots, $m_n$, was selected by means of the BIC procedure (resulting $\widehat{m}_n=4$); for details on the role of the splines, see (\ref{theta_0_base}) in Section \ref{design-sim}.
	
	\subsection{Tecator's data}
\begin{figure}[h]
		\vspace{-30pt}
		\centering
		\includegraphics[width=\textwidth]{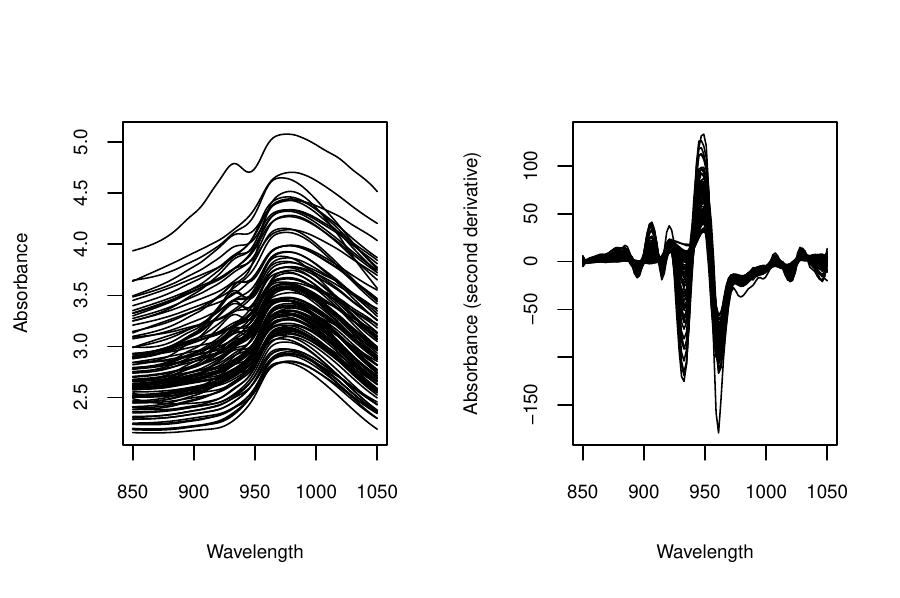}
		\vspace{-20pt}
		\caption{Sample of 100 absorbance curves $\mathcal{X}$ (left panel) together with their second derivatives $\mathcal{X}^{(2)}$ (right panel).}
		\label{fig5}
	\end{figure}

	The real data application will be focused on the well-known Tecator's data, which is possibly one of the data sets most used in the statistical literature to illustrate the usefulness of different functional nonparametric or semiparametric procedures (despite the known fact that Tecator's data contains several duplicated curves). Tecators's data include the percentages of fat, protein and moisture contents, and the near-infrared absorbance spectra of $215$ finely chopped pieces of meat. For each piece of meat, the percentages of fat, protein and moisture contents ($Y_i$, $X_{1i}$ and $X_{2i}$, respectively) are scalar, while the corresponding near-infrared absorbance spectra observations were collected on $100$ equally spaced wavelengths ($t_j$, $j=1,\dots,100$) in the range $850$--$1050$ $nm$; so each subject can be considered as a continuous curve, $\mathcal{X}_i$. As usual when one deals with Tecator's data set, %(see e.g. \citealt{aneiros_2016}), 
	we will use the second derivatives of the absorbance curves, $\mathcal{X}_i^{(2)}$, as functional covariate instead of the original curve (see eg \citealt*{ferratyvieu_2006} for details). Figure \ref{fig5} displays samples of both the absorbance curves and their second derivatives.
	
	Our purpose is modelling the relationship between the fat content (response), the protein and moisture contents (scalar covariates), and the absorbance spectra (functional covariate) and then, use the model to predict the fat content. In addition, we are interested in whether there are any interaction effects, quadratic effects and/or cubic effects between these scalar covariates.
	
	In order to compare the behaviour of each considered model and estimation procedure, we will split the original sample into two subsamples:
	a training sample, $$\mathcal{T}_1=\{(X_{i1},X_{i2},\mathcal{X}_i^{(2)},Y_i)\}_{i=1}^{160},$$ and a testing one, $$\mathcal{T}_2=\{(X_{i1},X_{i2},\mathcal{X}_i^{(2)},Y_i)\}_{i=161}^{215}.$$ In this way, all the estimation task is made only by means of the training sample, while the testing sample is used to measure the quality of the predictions. To quantify the error in the prediction task, the mean square error of prediction (MSEP) will be used:
	\begin{equation}
		\textrm{MSEP}=\frac{1}{55}\sum_{i=161}^{215}\left(Y_i-\widehat{Y}_i\right)^2,
		\label{MSEP}
	\end{equation}
	where $\widehat{Y}_i$ is the predicted value for $Y_i$ obtained from each considered model and estimation procedure.
	
	\subsection{Modelling, variable selection and prediction}
	\label{mod-pred}
	In literature, several models have been used to describe the relation between the fat content and the absorbance spectra (see eg \citealt*{ferratyvieu_2006} for a functional nonparametric model, and \citealt*{chen_2011} for a multiple index functional model). \cite{novo} modelled this data set using the functional single-index  model (FSIM) and compared the performance of the obtained predictions  with that provided by the functional linear model (FLM) and the pure functional nonparametric model (FNM). %; see Table \ref{table3}. 
	Such three models as well as the corresponding MSEPs obtained from kernel-based estimation procedures (in the case of FNM and FSIM) and functional principal components regression (in the case of FLM; see eg \citealt*{aguetal10} for partial least squares regression, including an application to Tecator's data, and \citealt*{febetal17} for a comparative study between these two dimensionality reduction techniques) are summarized in Table \ref{table4} (for details on the results corresponding to the FLM, see  \citealt*{ferraty_2013}; see also  \citealt*{novo} for details related to the FNM and FSIM).
	
	\begin{table}[h]
		\centering
		\begin{tabular}{|llc|}
			\hline
			& Model & MSEP \\ 
			\hline
			FLM: & $Y=\alpha _{0}+\int_{850}^{1050}\mathcal{X}^{(2)}(t)\alpha (t)dt+\varepsilon $ & $7.17$ \\
			& & \\
			FNM: & $Y=r_1(\mathcal{X}^{(2)})+\varepsilon$ & $4.06$  \\ 	
			& & \\
			FSIM: & $Y=r_2\left( \left\langle \theta_{0},\mathcal{X}^{(2)}\right\rangle \right)
			+\varepsilon $ & $3.49$\\ 
			%	& & \\
			%	LM: & $Y=\sum_{j=1}^{7}X_j\beta _{j}++\varepsilon$ & 1.95 ($X_1,X_2,X_7$)
			%	\\
			%	& & \\
			%	AM: & $Y=\mu +f_1(X_1)+f_2(X_2)+f_7(X_7)+\varepsilon$ & 1.96 ($X_1,X_2,X_7$ nodos opt 2)
			%\\
			\hline
		\end{tabular}
		\caption{Values of the MSEPs from some functional models.}
		\label{table4}
	\end{table}
	To improve the performance of the FNM and FSIM, \cite{aneiros_2006} and \cite{wang_2016} included   in such models, respectively, information from the scalar covariates $X_1$ and $X_2$. Nevertheless, in those two papers only linear effects of $X_1$ and $X_2$ were considered: no interaction effects, and neither quadratic nor cubic etc effects. In order to take into account such potential effects, one can extend the case studies of \cite{aneiros_2006} and \cite{wang_2016} by considering as linear covariates $X_{2j-1}=X_1^j$ and $X_{2j}=X_2^j$ ($j=1,\cdots,q_n$), and $X_{p_n}=X_1X_2$ (we have denoted $p_n=2q_n+1$). The corresponding semi-functional partial linear model (SFPLM) and SSFPLSIM  for the particular case of $p_n=7$ (equivalently, for models allowing linear, quadratic and cubic effects, as well as interaction between the covariates $X_1$ and $X_2$) are shown in Table \ref{table5}. The sparse linear model (SLM) is also included in such table. Table \ref{table5} also reports the selected variables when the PLS procedures in \cite{fanpeng_2004} (SLM), \cite{aneiros_2015} (SFPLM) and our proposal (SPLSIM) are applied, as well as the corresponding MSEPs.
	
	\begin{table}[h]
		
		\centering
		\begin{tabular}{|lllc|}
			\hline
			& Model & Selected variables & MSEP \\ 
			\hline
			SLM: & $Y=\sum_{j=1}^{7}X_j\beta _{j}+\varepsilon $ & $X_1,X_2, X_7$ & $1.95$  \\ 	
			& & &\\
			
			SFPLM: & $Y=\sum_{j=1}^{7}X_j\beta _{j}+m_{1}(%
			\mathcal{X}^{(2)})+\varepsilon $ & $X_1,X_2$ & $1.48$  \\ 	
			& & &\\
			SSFPLSIM: & $Y=\sum_{j=1}^{7}X_j\beta_j+m\left(\left\langle \theta_{0},\mathcal{X}^{(2)}\right\rangle\right)+\varepsilon $ & $X_1,X_2,X_4,X_5$ & $1.29$ \\ 
			\hline
		\end{tabular}
		\caption{Values of the MSEPs from some scalar parametric and functional semiparametric models when PLS variable selection methods are used. The selected variables are also shown.}
		\label{table5}
	\end{table}
	
	Several conclusions can be drawn from Tables \ref{table4} and \ref{table5}. First, Table \ref{table4} shows that the functional semiparametric model (FSIM) improves both the functional linear (FLM) and nonparametric (FNM) ones. Second, Table \ref{table5} indicates that to add scalar linear effects in the FNM and FSIM (or, equivalently, to add functional nonparametric or semiparametric effects in the SLM) improves the predictive power of these simpler models. In addition, the percentages of protein ($X_1$) and moisture ($X_2$) contents linearly influence on the percentage of fat content ($Y$). $X_1$ and $X_2$ also present cubic and quadratic influence on $Y$, respectively, when the SSFPLSIM is considered, while interaction effects (the covariate $X_7$ is selected) only are detected from the SLM. Finally, Table \ref{table5} also shows that our proposed model (SSFPLSIM), which is a mix of all these ideas (semiparametric and partial linear ideas), presents the better performance.
	
\begin{figure}[h]
		\centering
		\includegraphics[width=6cm]{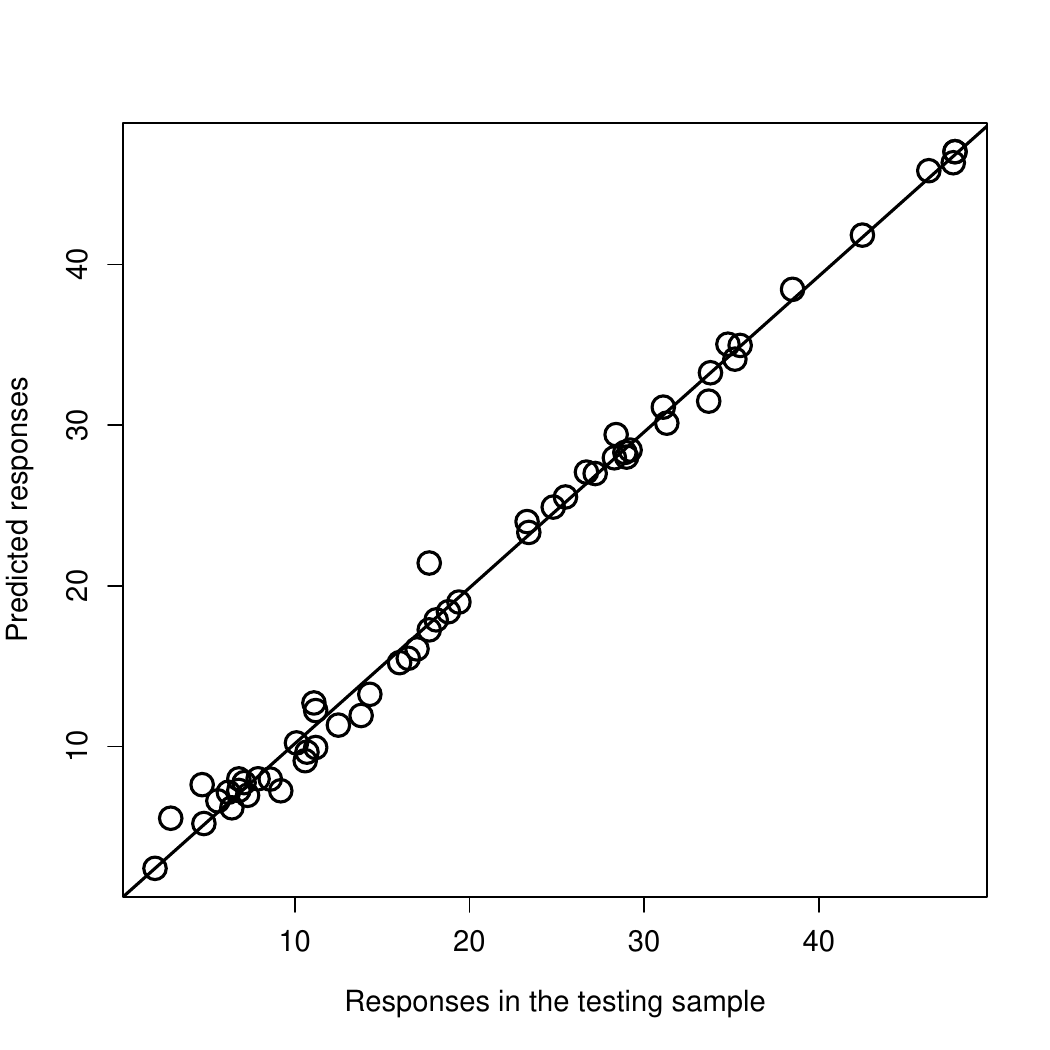}%\hspace{1.5cm}
		\caption{Predicted values from the SSFPLSIM vs Observed values.}
		\label{fig6}
	\end{figure}

\begin{figure}[h]
		\centering
		\includegraphics[width=6cm]{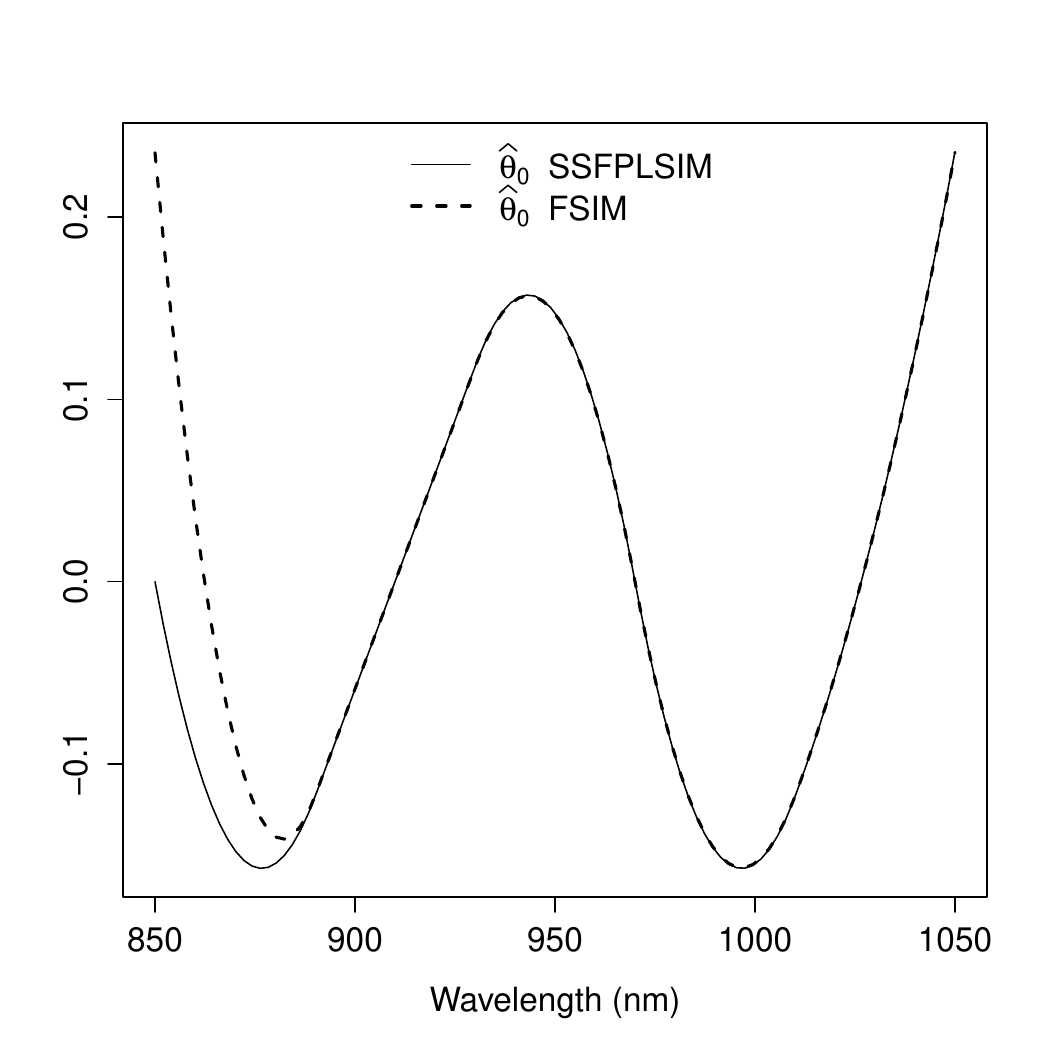}%\hspace{1.5cm}
		\includegraphics[width=6cm]{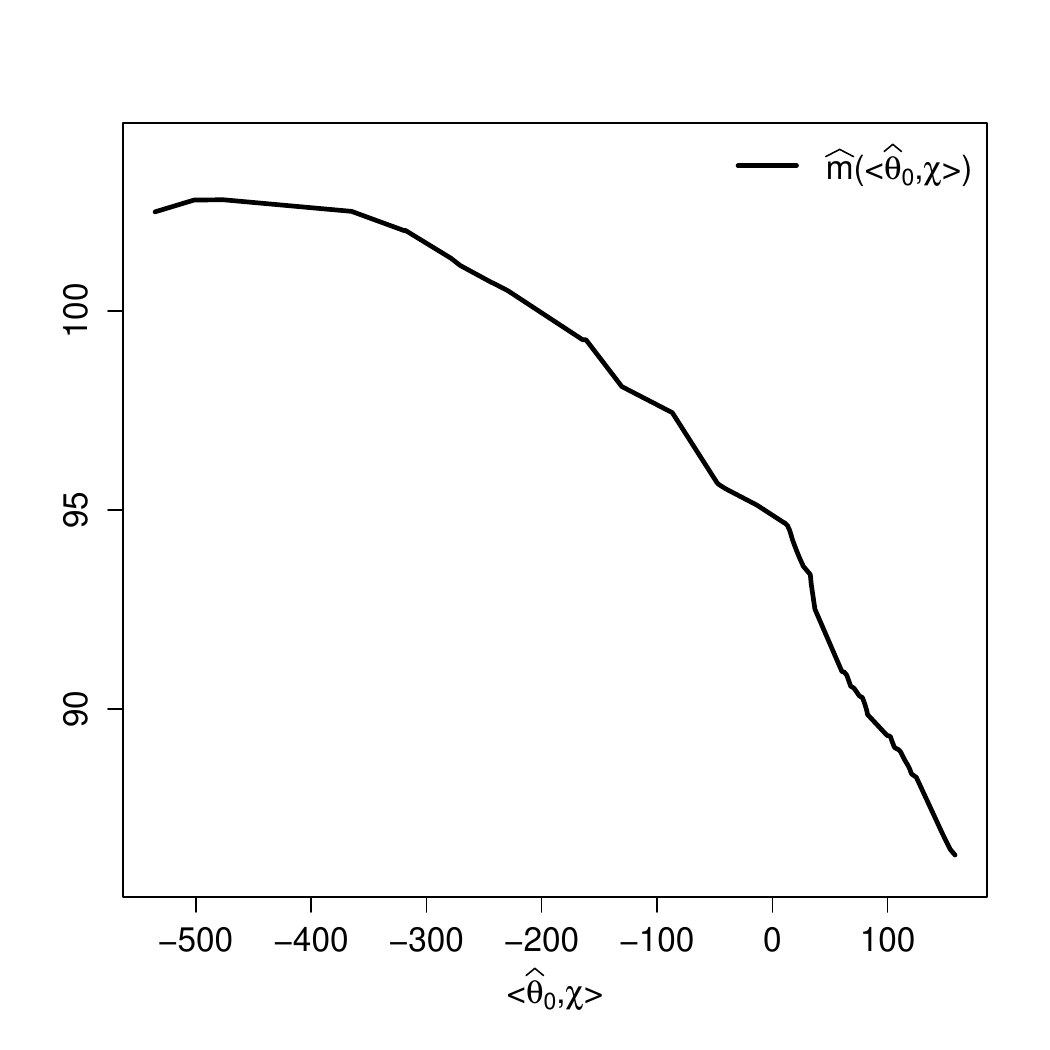}
		\caption{Left panel: Estimates of the functional directions ($\theta_0$) in the FSIM and SSFPLSIM. Right panel: Estimate of the function $m$ in the SSFPLSIM.}
		\label{fig7}
	\end{figure}

	Figure \ref{fig6} displays the predicted values ($\widehat{Y}_i, \ i=161,\ldots,215$) from the SSFPLSIM versus the observed ones ($Y_i, \ i=161,\ldots,215$). The high predictive power of the SSFPLSIM is evident. The estimates of the functional directions, $\theta_0$, in the FSIM and SSFPLSIM are displayed in Figure \ref{fig7} (left panel). It is worth being noted that both graphics of $\widehat{\theta}_0$ suggests that the two bumps
	around wavelengths 880 and 1000, as well as the peak around wavelength 940, could be
	important indicators of the fat content (note that this suggestion is compatible with the
	findings in \citealt*{novo}). Finally, Figure \ref{fig7} (right panel) shows the estimate of the smooth real-valued function, $m$, in the SSFPLSIM.  
	
	\subsection{Summary}
	%To sum up, 
	Our real data application evidences the advantages of using the SSFPLSIM together with the proposed PLS procedure in terms of accuracy of predictions. 
	% and number of selected linear covariates. 
	In addition, as in the case of FSIM, the SSFPLSIM presents the advantage of the interpretation of the estimated direction of projection, $\widehat{\theta}_0$, which could also complement the information about how the (second derivative of the) spectrometric curves affect to the fat content.
	
	\bigskip
	\noindent{\large \bf{Acknowledgements}}
	
	\noindent The authors wish to thank two anonymous referees for their helpful comments and suggestions, which greatly improved the quality of this paper. This work was supported in part by the Spanish Ministerio de Econom\'ia y Competitividad under Grant MTM2014-52876-R and
	Grant MTM2017-82724-R, in part by the Xunta de Galicia through Centro Singular de Investigaci\'on de Galicia accreditation
	under Grant ED431G/01 2016-2019 and through the Grupos de Referencia Competitiva under Grant ED431C2016-015,
	and in part by the European Union (European Regional Development Fund - ERDF). The first author also thanks the financial support from the Xunta de Galicia and the European Union (European Social Fund - ESF), the reference of which is ED481A-2018/191.

\newpage

\appendix

\begin{center}
{\large{\textbf{Sparse semiparametric regression when predictors are mixture of functional and high-dimensional variables}}}

Silvia Novo$^a$\footnote{Corresponding author email address: \href{s.novo@udc.es}{s.novo@udc.es} } \hspace{2pt} Germ\'an Aneiros$^b$ \hspace{2pt} Philippe Vieu$^c$
\end{center}

\noindent$^a$ Department of Mathematics, MODES, CITIC, Universidade da Coruña, A Coruña, Spain \\
\noindent$^b$ Department of Mathematics, MODES, CITIC, ITMATI, Universidade da Coruña, A Coruña, Spain\\
\noindent$^c$ Institut de Math\'{e}matiques, Universit\'e  Paul Sabatier,  Toulouse, France			

\bigskip

\bigskip

\noindent \textbf{SUPPLEMENTARY MATERIAL:} These pages contain the proofs of the asymptotic results presented in our paper. In addition, they show some known lemmas used to prove such results. Novel lemmas, as well as their proofs, are also presented here. \\

The used assumptions and the enunciates of the theorems, as well as the used references, can be found in the paper. All the enumeration corresponding to the paper is maintained here (this includes enumeration related to equations, assumptions, Theorems, Remarks, Sections, etc).

\section{Proofs}
This section presents the proofs of our main results. For that, a main role is played by the technical lemmas provided in Section \ref{results-lemas}. Note that the Remark \ref{rem-app} 
in Section \ref{general-assump} justifies that such lemmas can be applied under the conditions of our theorems.

Without loss of generality, we will assume that $S_n=\{1,\dots,s_n\}$.
\subsection{Proof of Theorem \ref{teorema1}}
Before starting the proof, let us complete in the following way the notations introduced in Sections \ref{some-notation} and \ref{notation} of the paper: 
\begin{equation*}
\pmb{g}_{j,\theta_0}=\left(g_{j,\theta_0}(\mathcal{X}_1),\dots,g_{j,\theta_0}(\mathcal{X}_n)\right)^{\top} \ (0\leq j\leq p_n),\ \pmb{G}_{\theta_0}=\left(\pmb{g}_{1,\theta_0},\dots,\pmb{g}_{p_n,\theta_0}\right),
\end{equation*} 	
and for each $\theta\in\Theta_n$:
\begin{align*}
&\widehat{g}_{j,\theta}(\chi)=\sum_{i=1}^n w_{n,h,\theta}(\chi,\mathcal{X}_i)Z_{ij} \mbox{  with  }
Z_{i0}=Y_i \mbox{  and  } Z_{ij}=X_{ij} \ (1 \leq j \leq p_n),\nonumber\\ &\widehat{\pmb{g}}_{j,\theta}=\left(\widehat{g}_{j,\theta}(\mathcal{X}_1),\dots,\widehat{g}_{j,\theta}(\mathcal{X}_n)\right)^{\top} \ (0\leq j\leq p_n)\mbox{ and }
\widehat{\pmb{G}}_{j,\theta}=\left(\widehat{\pmb{g}}_{1,\theta},\dots,\widehat{\pmb{g}}_{p_n,\theta}\right).\nonumber
\end{align*}
In addition, we denote
\begin{equation*}
\pmb{\eta}_{\theta_0}=\left(\pmb{\eta}_{1,\theta_0},\dots,\pmb{\eta}_{n,\theta_0}\right)^{\top}%, \ \pmb{\beta}_{S_n}=\left(\beta_1,\dots,\beta_{s_n}\right)^{\top}\in\mathbb{R}^{s_n}
\end{equation*}
and
\begin{equation*}
\mathcal{Q}^*\left(\pmb{\beta}_{S_n},\theta\right)=\frac{1}{2}\left(\widetilde{\pmb{Y}}_{\theta}-\widetilde{\pmb{X}}_{\theta S_n}\pmb{\beta}_{S_n}\right)^{\top}\left(\widetilde{\pmb{Y}}_{\theta}-\widetilde{\pmb{X}}_{\theta S_n}\pmb{\beta}_{S_n}\right)+n\sum_{j=1}^{s_n}\mathcal{P}_{\lambda_{j_n}}\left(|\beta_j|\right).
\end{equation*}

To obtain the desired result, it suffices to prove that there exists a local minimizer $\left(\widehat{\pmb{\beta}}_{0S_n},\widehat{\theta}_0\right)$ of $\mathcal{Q}^*(\pmb{\beta}_{S_n},\theta)$ such that
\begin{equation}
\left\|\widehat{\pmb{\beta}}_{0S_n}-\pmb{\beta}_{0S_n}\right\|=O_p(u_n), \ \left\| \widehat{\theta}_0-\theta_0\right\|=O_p(v_n)
\label{rates}
\end{equation}
and
\begin{equation}
\left(\widehat{\pmb{\beta}}_0,\widehat{\theta}_0\right)=\left(\left(\widehat{\pmb{\beta}}_{0S_n}^{\top},\pmb{0}_{p_n-s_n}^{\top}\right)^{\top},\widehat{\theta}_0\right) \textrm{ is a local minimizer of $\mathcal{Q}(\cdot,\cdot)$},
\label{minimo}
\end{equation}
where $\pmb{0}_{p_n-s_n}$ is a vector of zero components with dimension $p_n-s_n$.

First, we will obtain the results in (\ref{rates}). For that, it suffices to show that, for any given $\gamma>0$, there exists a constant $C$ such that, for $n$ large enough,
\begin{equation*}
\mathbb{P}\left(\inf_{\lvert\lvert\pmb{u}\rvert\rvert=C, \theta \in \Theta_n^\ast}\mathcal{Q}^*\left(\pmb{\beta}_{0S_n}+u_n\pmb{u},\theta\right)>\mathcal{Q}^*\left(\pmb{\beta}_{0S_n},\theta_0\right)\right)\geq 1-\gamma,
\end{equation*}
where $\pmb{u}=\left(u_1,\dots,u_{s_n}\right)^{\top}\in\mathbb{R}^{s_n}$ and $\Theta_n^\ast=\{\theta \in \Theta_n; d(\theta,\theta_0)=v_n\}$.

Let us denote 
\begin{equation}\mathcal{Q}^*(\pmb{\beta}_{S_n},\theta)=\mathcal{L}^*(\pmb{\beta}_{S_n},\theta)+\mathcal{P}^*(\pmb{\beta}_{S_n}),\label{def_Q*}\end{equation}
where
\[
\mathcal{L}^*(\pmb{\beta}_{S_n},\theta)=\frac{1}{2}\left(\widetilde{\pmb{Y}}_{\theta}-\widetilde{\pmb{X}}_{\theta S_n}\pmb{\beta}_{S_n}\right)^{\top}\left(\widetilde{\pmb{Y}}_{\theta}-\widetilde{\pmb{X}}_{\theta S_n}\pmb{\beta}_{S_n}\right)
\mbox{ and }
\mathcal{P}^*(\pmb{\beta}_{S_n})=n\sum_{j \in S_n}\mathcal{P}_{\lambda_{j_n}}\left(|\beta_j|\right).
\]
We have that
\begin{equation}
\mathcal{Q}^*(\pmb{\beta}_{0S_n},\theta_0)-\mathcal{Q}^*\left(\pmb{\beta}_{0S_n}+u_n\pmb{u},\theta\right)=A_1+A_2,
\label{Q*_inicial}
\end{equation}
where
\begin{equation}
A_1=\mathcal{L}^*(\pmb{\beta}_{0S_n},\theta_0)-\mathcal{L}^*\left(\pmb{\beta}_{0S_n}+u_n\pmb{u},\theta\right)
\mbox{ and } A_2=\mathcal{P}^*(\pmb{\beta}_{0S_n})-\mathcal{P}^*(\pmb{\beta}_{0S_n}+u_n\pmb{u}).
\label{a}
\end{equation}
Focusing on $A_1$, we can write
\begin{eqnarray}
2A_1&=&\left(\widetilde{\pmb{Y}}_{\theta_0}^{\top}\widetilde{\pmb{Y}}_{\theta_0}-2\widetilde{\pmb{Y}}_{\theta_0}^{\top}\widetilde{\pmb{X}}_{\theta_0S_n}\pmb{\beta}_{0S_n}\right)+\left(\widetilde{\pmb{X}}_{\theta_0S_n}\pmb{\beta}_{0S_n}\right)^{\top}\widetilde{\pmb{X}}_{\theta_0S_n}\pmb{\beta}_{0S_n}
\nonumber\\
&-&\left(\widetilde{\pmb{Y}}_{\theta}^{\top}\widetilde{\pmb{Y}}_{\theta}-2\widetilde{\pmb{Y}}_{\theta}^{\top}\widetilde{\pmb{X}}_{\theta S_n}\pmb{\beta}_{0S_n}\right)
\nonumber\\
&-&\left(\widetilde{\pmb{X}}_{\theta S_n}\pmb{\beta}_{0S_n}\right)^{\top}\widetilde{\pmb{X}}_{\theta S_n}\pmb{\beta}_{0S_n}+ 2u_n\left(\widetilde{\pmb{Y}}_{\theta}^{\top}\widetilde{\pmb{X}}_{\theta S_n}- \left(\widetilde{\pmb{X}}_{\theta S_n}\pmb{\beta}_{0S_n}\right)^{\top}\widetilde{\pmb{X}}_{\theta S_n} \right)\pmb{u}
\nonumber\\
&-&u_n^2\pmb{u}^{\top}\widetilde{\pmb{X}}_{\theta S_n}^{\top}\widetilde{\pmb{X}}_{\theta S_n}\pmb{u}\equiv A_{11}+A_{12}-A_{13}-A_{14}+2u_nA_{15}-A_{16}.\label{descomposicion_L}
\end{eqnarray}
Taking into account that
\begin{equation*}
\widetilde{\pmb{Y}}_{\theta_0}=\pmb{g}_{0,\theta_0}-\widehat{\pmb{g}}_{0,\theta_0}+\pmb{\eta}_{\theta_0S_n}\pmb{\beta}_{0S_n}+\pmb{\varepsilon} \mbox{ and }
\widetilde{\pmb{Y}}_{\theta}=\pmb{g}_{0,\theta_0}-\widehat{\pmb{g}}_{0,\theta}+\pmb{\eta}_{\theta_0S_n}\pmb{\beta}_{0S_n}+\pmb{\varepsilon},
\end{equation*}
as well as that
\begin{equation*}
\widetilde{\pmb{X}}_{\theta_0 S_n}=\left(\pmb{G}_{\theta_0}-\widehat{\pmb{G}}_{\theta_0}\right)_{S_n}+\pmb{\eta}_{\theta_0S_n} \mbox{ and }
\widetilde{\pmb{X}}_{\theta S_n}=\left(\pmb{G}_{\theta_0}-\widehat{\pmb{G}}_{\theta}\right)_{S_n}+\pmb{\eta}_{\theta_0S_n},
\end{equation*}
we obtain that
\begin{eqnarray}
A_{11}
&=&\left(\pmb{g}_{0,\theta_0}-\widehat{\pmb{g}}_{0,\theta_0}\right)^{\top}\left(\pmb{g}_{0,\theta_0}-\widehat{\pmb{g}}_{0,\theta_0}\right)+2\left(\pmb{g}_{0,\theta_0}-\widehat{\pmb{g}}_{0,\theta_0}\right)^{\top}\pmb{\varepsilon}\nonumber\\
&-&2\left(\pmb{g}_{0,\theta_0}-\widehat{\pmb{g}}_{0,\theta_0}\right)^{\top}\left(\pmb{G}_{\theta_0}-\widehat{\pmb{G}}_{\theta_0}\right)_{S_n}\pmb{\beta}_{0S_n}-2\pmb{\varepsilon}^{\top}\left(\pmb{G}_{\theta_0}-\widehat{\pmb{G}}_{\theta_0}\right)_{S_n}\pmb{\beta}_{0S_n}\nonumber\\
&+&\pmb{\varepsilon}^{\top}\pmb{\varepsilon}-\left(\pmb{\eta}_{\theta_0S_n}\pmb{\beta}_{0S_n}\right)^{\top}\pmb{\eta}_{\theta_0S_n}\pmb{\beta}_{0S_n}-2\left(\pmb{\eta}_{\theta_0S_n}\pmb{\beta}_{0S_n}\right)^{\top}\left(\pmb{G}_{\theta_0}-\widehat{\pmb{G}}_{\theta_0}\right)_{S_n}\pmb{\beta}_{0S_n},\label{expr1}
\end{eqnarray}
\begin{eqnarray}
A_{12}&=&
\pmb{\beta}_{0S_n}^{\top}\left(\pmb{G}_{\theta_0}-\widehat{\pmb{G}}_{\theta_0}\right)_{S_n}^{\top}\left(\pmb{G}_{\theta_0}-\widehat{\pmb{G}}_{\theta_0}\right)_{S_n}\pmb{\beta}_{0S_n}\nonumber\\&+&2\pmb{\beta}_{0S_n}^{\top}\left(\pmb{G}_{\theta_0}-\widehat{\pmb{G}}_{\theta_0}\right)_{S_n}^{\top}\pmb{\eta}_{\theta_0 S_n}\pmb{\beta}_{0S_n}+\left(\pmb{\eta}_{\theta_0 S_n}\pmb{\beta}_{0S_n}\right)^{\top}\pmb{\eta}_{\theta_0 S_n}\pmb{\beta}_{0S_n},\label{expr3}
\end{eqnarray}
\begin{eqnarray}
A_{13}&=&\left(\pmb{g}_{0,\theta_0}-\widehat{\pmb{g}}_{0,\theta}\right)^{\top}\left(\pmb{g}_{0,\theta_0}-\widehat{\pmb{g}}_{0,\theta}\right)+2\left(\pmb{g}_{0,\theta_0}-\widehat{\pmb{g}}_{0,\theta}\right)^{\top}\pmb{\varepsilon}\nonumber\\
&-&2\left(\pmb{g}_{0,\theta_0}-\widehat{\pmb{g}}_{0,\theta}\right)^{\top}\left(\pmb{G}_{\theta_0}-\widehat{\pmb{G}}_{\theta}\right)_{S_n}\pmb{\beta}_{0S_n}-2\pmb{\varepsilon}^{\top}\left(\pmb{G}_{\theta_0}-\widehat{\pmb{G}}_{\theta}\right)_{S_n}\pmb{\beta}_{0S_n}\nonumber\\
&+&\pmb{\varepsilon}^{\top}\pmb{\varepsilon}-\left(\pmb{\eta}_{\theta_0 S_n}\pmb{\beta}_{0S_n}\right)^{\top}\pmb{\eta}_{\theta_0 S_n}\pmb{\beta}_{0S_n}-2\left(\pmb{\eta}_{\theta_0S_n}\pmb{\beta}_{0S_n}\right)^{\top}\left(\pmb{G}_{\theta_0}-\widehat{\pmb{G}}_{\theta}\right)_{S_n}\pmb{\beta}_{0S_n},
\label{expr2}
\end{eqnarray}
\begin{eqnarray}
A_{14}&=&
\pmb{\beta}_{0S_n}^{\top}\left(\pmb{G}_{\theta_0}-\widehat{\pmb{G}}_{\theta}\right)_{S_n}^{\top}\left(\pmb{G}_{\theta_0}-\widehat{\pmb{G}}_{\theta}\right)_{S_n}\pmb{\beta}_{0S_n}\nonumber\\&+&2\pmb{\beta}_{0S_n}^{\top}\left(\pmb{G}_{\theta_0}-\widehat{\pmb{G}}_{\theta}\right)_{S_n}^{\top}\pmb{\eta}_{\theta_0 S_n}\pmb{\beta}_{0S_n}+\left(\pmb{\eta}_{\theta_0 S_n}\pmb{\beta}_{0S_n}\right)^{\top}\pmb{\eta}_{\theta_0 S_n}\pmb{\beta}_{0S_n}\label{expr4}
\end{eqnarray}
and
\begin{eqnarray}
A_{15}
&=&\left(\pmb{g}_{0,\theta_0}-\widehat{\pmb{g}}_{0,\theta}\right)^{\top}\left(\pmb{G}_{\theta_0}-\widehat{\pmb{G}}_{\theta}\right)_{S_n}\pmb{u}+\left(\pmb{g}_{0,\theta_0}-\widehat{\pmb{g}}_{0,\theta}\right)^{\top}\pmb{\eta}_{\theta_0S_n}\pmb{u}\nonumber\\
&-&
\pmb{\beta}_{0S_n}^{\top}\left(\left(\pmb{G}_{\theta_0}-\widehat{\pmb{G}}_{\theta}\right)_{S_n}^{\top}\left(\pmb{G}_{\theta_0}-\widehat{\pmb{G}}_{\theta}\right)_{S_n}+\left(\pmb{G}_{\theta_0}-\widehat{\pmb{G}}_{\theta}\right)_{S_n}^{\top}\pmb{\eta}_{\theta_0S_n}\right)\pmb{u}\nonumber\\
&+&\pmb{\varepsilon}^{\top}\left(\pmb{G}_{\theta_0}-\widehat{\pmb{G}}_{\theta}\right)_{S_n}\pmb{u}+\pmb{\varepsilon}^{\top}\pmb{\eta}_{\theta_0S_n}\pmb{u}.
\label{expr5}
\end{eqnarray}
Let us denote
\begin{equation} 
B=A_{11}+A_{12}-A_{13}-A_{14}. \label{b}
\end{equation}
From decompositions (\ref{expr1})-(\ref{expr4}), it is easy to obtain that
\begin{eqnarray}
B &=&
\left(\widehat{\pmb{g}}_{0,\theta_0}-\widehat{\pmb{g}}_{0,\theta}\right)^{\top}\left(\widehat{\pmb{g}}_{0,\theta_0}+\widehat{\pmb{g}}_{0,\theta}\right)+2\pmb{g}_{0,\theta_0}^{\top}\left(\widehat{\pmb{g}}_{0,\theta}-\widehat{\pmb{g}}_{0,\theta_0}\right)+2\pmb{\varepsilon}^{\top}\left(\widehat{\pmb{g}}_{0,\theta}-\widehat{\pmb{g}}_{0,\theta_0}\right)\nonumber\\
&+&2\pmb{g}_{0,\theta_0}^{\top}\left(\widehat{\pmb{G}}_{\theta_0}-\widehat{\pmb{G}}_{\theta}\right)_{S_n}\pmb{\beta}_{0S_n}+2\left(\widehat{\pmb{g}}_{0,\theta_0}-\widehat{\pmb{g}}_{0,\theta}\right)^{\top}\pmb{G}_{\theta_0S_n}\pmb{\beta}_{0S_n}\nonumber\\
&+&2\left(\widehat{\pmb{g}}_{0,\theta}^{\top}\widehat{\pmb{G}}_{\theta S_n}-\widehat{\pmb{g}}_{0,\theta_0}^{\top}\widehat{\pmb{G}}_{\theta_0 S_n}\right)^{\top}\pmb{\beta}_{0S_n}+2\pmb{\varepsilon}^{\top}\left(\widehat{\pmb{G}}_{\theta_0}-\widehat{\pmb{G}}_{\theta}\right)_{S_n}\pmb{\beta}_{0S_n}\nonumber\\
&+&2\pmb{\beta}_{0S_n}^{\top}\widehat{\pmb{G}}_{\theta_0 S_n}^{\top}\left(\widehat{\pmb{G}}_{\theta}-\widehat{\pmb{G}}_{\theta_0 }\right)_{S_n}\pmb{\beta}_{0S_n}+\pmb{\beta}_{0S_n}^{\top}\left(\widehat{\pmb{G}}_{\theta_0}-\widehat{\pmb{G}}_{\theta }\right)_{S_n}^{\top}\left(\widehat{\pmb{G}}_{\theta_0}+\widehat{\pmb{G}}_{\theta}\right)_{S_n}\pmb{\beta}_{0S_n}\nonumber\\
&=&B_1+B_2+B_3+B_4+B_5+B_6+B_7+B_8+B_9.
\label{expr6}
\end{eqnarray}
Now, we are going to obtain bounds (in probability) for each term, $B_k$ ($k=1,\ldots,9$), in (\ref{expr6}). Let us denote, for $0\leq j \leq p_n$, $$\left(\widehat{\pmb{g}}_{j,\theta_0}-\widehat{\pmb{g}}_{j,\theta}\right)=\left(d_{j1}',\dots,d_{jn}'\right)^{\top} \mbox{ and } \left(\widehat{\pmb{g}}_{j,\theta_0}+\widehat{\pmb{g}}_{j,\theta}\right)=\left(d_{j1}'',\dots,d_{jn}''\right)^{\top}.$$ On the one hand, from Lemma \ref{lema_novo} we have that
\begin{eqnarray}
\max_{0\leq j \leq p_n}\sup_{\theta \in \Theta_n^\ast}\max_{1\leq i\leq n}|d_{ji}'|=O_p\left(\frac{v_n}{hf(h)}\right).\label{num1}
\end{eqnarray}
On the other hand, from the uniform convergence of $\widehat{g}_{j,\theta}(\chi)$ to $g_{j,\theta_0}(\chi)$ (see Lemma \ref{lemaA3}) together with the fact that 
\begin{equation}
\max_{0\leq j\leq n}\max_{1\leq i\leq n} |g_{j,\theta_0}(\mathcal{X}_i)|=O(1)\label{num2}
%\nonumber
\end{equation}
(see Assumption (\ref{mom_1})), we obtain that
\begin{equation}
\max_{0\leq j \leq p_n}\sup_{\theta \in \Theta_n^\ast}\max_{1\leq i\leq n} |\widehat{g}_{j,\theta}(\mathcal{X}_i)|=O_p(1) \label{num3};
\end{equation}
so,
\begin{eqnarray}
\max_{0\leq j \leq p_n}\sup_{\theta \in \Theta_n^\ast}\max_{1\leq i\leq n}|d_{ji}''|=O_p\left(1\right). \label{num4}
\end{eqnarray}
Taking into account (\ref{num1}) and (\ref{num4}), we have that
\begin{eqnarray}
\left|B_1\right|=\left|\sum_{i=1}^nd_{0i}'d_{0i}''\right|\leq n\max_{1\leq i\leq n}\left|d_{0i}'\right|\max_{1\leq i\leq n}\left|d_{0i}''\right|=O_p\left(n\frac{v_n}{h f(h)}\right) \mbox{ uniformly on } \theta \in \Theta_n^\ast.\nonumber \\\label{b1}
\label{Op1}
\end{eqnarray}
From (\ref{num1}) and (\ref{num2}) we obtain that
\begin{eqnarray}
\left|B_2\right|=2\left|\sum_{i=1}^ng_{0,\theta_0}\left(\mathcal{X}_i\right)d_{0i}'\right|\leq 2n \max_{1\leq i\leq n}\left|g_{0,\theta_0}\left(\mathcal{X}_i\right)\right| \max_{1\leq i\leq n}\left|d_{0i}'\right|=O_p\left(n\frac{v_n}{hf(h)}\right)  \label{b2}
\end{eqnarray}
uniformly on $\theta \in \Theta_n^\ast$.
From Lemma \ref{lemaA1} and expression (\ref{num1}) we obtain that
\begin{eqnarray}
\left|B_3\right|&=&2\left|\sum_{i=1}^{n} d_{0i}'\varepsilon_i\right|=O_p\left(n^{1/2+1/r_{\varepsilon}} \frac{v_n}{hf(h)} \log n\right)  \mbox{ uniformly on } \theta \in \Theta_n^\ast.\label{b3}
\end{eqnarray}
If one takes Assumption (\ref{null_par3}) into account, then similar reasonings as those used to obtain (\ref{b1})-(\ref{b3}) give that, uniformly on $\theta \in \Theta_n^\ast$,
\begin{equation}
B_4=O_p\left(n s_n\frac{v_n}{hf(h)}\right), \
%.\label{Op3}
B_5=O_p\left(ns_n\frac{v_n}{hf(h)}\right), \
%.\label{Op5}
B_7=O_p\left(n^{1/2+1/r_{\varepsilon}}s_n\frac{v_n}{hf(h)}\log n\right)
\label{b457}
\end{equation}
and
\begin{equation}
B_8=O_p\left(\frac{n s_n^2v_n}{hf(h)}\right), \ B_9=O_p\left(ns_n^2\frac{v_n}{hf(h)}\right).
\label{b89}
\end{equation}
The term $B_6$ can be re-written in the following manner:
\begin{equation}
B_6= \left(\widehat{\pmb{g}}_{0,\theta}-\widehat{\pmb{g}}_{0,\theta_0}\right)^{\top}\widehat{\pmb{G}}_{\theta S_n}\pmb{\beta}_{0S_n}+\widehat{\pmb{g}}_{0,\theta_0}^{\top}\left(\widehat{\pmb{G}}_{\theta }-\widehat{\pmb{G}}_{\theta_0}\right)_{S_n}\pmb{\beta}_{0S_n}=B_{61}+B_{62}.
\nonumber %\label{num8}
\end{equation}
If one considers (\ref{num3}) instead of (\ref{num2}), then similar reasonings as those used to obtain the orders of $B_4$ and $B_5$ (see (\ref{b457})) give
$$
B_{61}=O_p\left(ns_n\frac{v_n}{hf(h)}\right) \mbox{ and } B_{62}=O_p\left(ns_n\frac{v_n}{hf(h)}\right) \mbox{ uniformly on } \theta \in \Theta_n^\ast
$$
respectively; so, we have that
\begin{equation}
B_{6}=O_p\left(ns_n\frac{v_n}{hf(h)}\right) \mbox{ uniformly on } \theta \in \Theta_n^\ast. \label{b6}
\end{equation}
It is  noteworthy that, as consequence of our assumptions, all the $O_p(\cdot)$ in (\ref{Op1})-(\ref{b6}) are $O_p(n u_n^2)$. Therefore, we have proved that
\begin{equation}
B=O_p(n u_n^2) \mbox{ uniformly on } \theta \in \Theta_n^\ast.
\label{b-order}
\end{equation}
The term $A_{15}$ (see (\ref{expr5})) can be studied in a similar way as (A6) in \cite{aneiros_2015}, but considering our Lemma \ref{lemaA3} instead of Lemma A.3 of \cite{aneiros_2015}. Specifically, denoting
\begin{equation}
r_n=\frac{\log p_n\psi_{\mathcal{C}}\left(1/n\right)}{nf(h)}, \label{rn}
\end{equation}
it can be obtained that:
\begin{equation}
\left(\pmb{g}_{0,\theta_0}-\widehat{\pmb{g}}_{0,\theta}\right)^{\top}\left(\pmb{G}_{\theta_0}-\widehat{\pmb{G}}_{\theta}\right)_{S_n}\pmb{u}=O_p\left(ns_n^{1/2}(h^{2\alpha}+r_n)\right)\left\|\pmb{u}\right\|  \mbox{ uniformly on } \theta \in \Theta_n^\ast \label{Op6}
\end{equation}
(consequence of Lemma \ref{lemaA3} and Cauchy-Schwarz inequality; take into account that, because it is assumed that $ns_nv_n=O(hf(h))$, it verifies that $r_n^{\ast}=r_n$, where $v_n$ was defined in (\ref{Theta}) while $r_n^{\ast}$ was defined in (\ref{rn.ast}) and used in Lemma \ref{lemaA3}),
\begin{equation}
\left(\pmb{g}_{0,\theta_0}-\widehat{\pmb{g}}_{0,\theta}\right)^{\top}\pmb{\eta}_{\theta_0S_n}\pmb{u}=O_p\left(n^{1/2}s_n^{1/2}(h^\alpha+r_n^{1/2})\log n\right)\left\|\pmb{u}\right\|  \mbox{ uniformly on } \theta \in \Theta_n^\ast \label{Op7}
\end{equation}
(consequence of Cauchy-Schwarz inequality, Lemma \ref{lemaA1} and Lemma \ref{lemaA3}),
\begin{equation}
\pmb{\beta}_{0S_n}^{\top}\left(\pmb{G}_{\theta_0}-\widehat{\pmb{G}}_{\theta}\right)_{S_n}^{\top}\left(\pmb{G}_{\theta_0}-\widehat{\pmb{G}}_{\theta}\right)_{S_n}\pmb{u}=O_p\left(ns_n^{3/2}(h^{2\alpha}+r_n)\right)\left\|\pmb{u}\right\|   \label{Op8}
\end{equation}
uniformly on $\theta \in \Theta_n^\ast$ (consequence of Cauchy-Schwarz inequality and Lemma \ref{lemaA3}),
\begin{equation}
\pmb{\beta}_{0S_n}^{\top}\left(\pmb{G}_{\theta_0}-\widehat{\pmb{G}}_{\theta}\right)_{S_n}^{\top}\pmb{\eta}_{\theta_0S_n}\pmb{u}=O_p\left(n^{1/2}s_n^{3/2}(h^\alpha+r_n^{1/2})\log n\right)\left\|\pmb{u}\right\|   \label{Op9}
\end{equation}
uniformly on $\theta \in \Theta_n^\ast$ (consequence of Cauchy-Schwarz inequality, Lemma \ref{lemaA2} and Lemma \ref{lemaA3}) and
\begin{equation}
\pmb{\varepsilon}^{\top}\left(\pmb{G}_{\theta_0}-\widehat{\pmb{G}}_{\theta}\right)_{S_n}\pmb{u}
=O_p\left(n^{1/2+1/r_{\varepsilon}}s_n^{1/2}(h^\alpha+r_n^{1/2})\log n\right)\left\|\pmb{u}\right\|   \label{Op10}
\end{equation}
uniformly on $\theta \in \Theta_n^\ast$ (consequence of Cauchy-Schwarz inequality, Lemma \ref{lemaA1} and Lemma \ref{lemaA3}). Then, from Lemma \ref{lemaA8} toghether with the fact that all the orders $O_p(\cdot)$ involved in (\ref{Op6})-(\ref{Op10}) are $O_p(n u_n)$, we obtain that
\begin{equation}
A_{15}=O_p(n u_n)\left\|\pmb{u}\right\|  \mbox{ uniformly on } \theta \in \Theta_n^\ast. \label{Op11}
\end{equation}
Now we focus on the term $A_{16}$ (see (\ref{descomposicion_L})). Using Lemma \ref{lemaA5} (considering $s_n$, $\pmb{X}_{\theta S_n}$ and $\pmb{B}_{\theta_0 S_n\times S_n}$ instead of $p_n$, $\pmb{X}_{\theta}$ and $\pmb{B}_{\theta_0}$, respectively) we have that
\begin{equation}
A_{16}=nu_n^2\left(\pmb{u}^{\top}\pmb{B}_{\theta_0S_n\times S_n}\pmb{u}+o_p(1)\right)  \label{main}
\end{equation}
uniformly over $\{\pmb{u}\in\mathbb{R}^{p_n}, ||\pmb{u}||=C\}$ and over $\theta\in\Theta_n^\ast$. From (\ref{descomposicion_L}), (\ref{b}), (\ref{b-order}), (\ref{Op11}) and (\ref{main}), we obtain that
\begin{equation}
2A_1
= O_p(nu_n^2)+O_p\left(n u_n^2\right)||\pmb{u}||-nu_n^2\left(\pmb{u}^{\top}\pmb{B}_{\theta_0S_n\times S_n}\pmb{u}+o_p(1)\right),
\label{L_final}
\end{equation}
where all the orders of convergence are uniform in $||\pmb{u}||=C$ and over $\theta\in\Theta_n^\ast$. Focusing now on $A_2$ (see (\ref{a})), let us note that $A_2$ is only linked to linear part of the model (\ref{modelo}). Therefore, from result (A14) in Aneiros et al (2015) we have that
\begin{equation}
A_2=O\left(n u_ns_n^{1/2}\delta_n\right)||\pmb{u}||+O\left(nu_n^2\rho_n\right)||\pmb{u}||^2.\label{P_final}
\end{equation}
Finally, from expressions (\ref{Q*_inicial}), (\ref{L_final}) and (\ref{P_final}), we obtain that
\begin{eqnarray}
\mathcal{Q}^*(\pmb{\beta}_{0S_n},\theta_0)-\mathcal{Q}^*\left(\pmb{\beta}_{0S_n}+u_n\pmb{u},\theta_0+v_n v\right)&=&O_p(nu_n^2)+O_p\left(n u_n^2\right)||\pmb{u}|| \nonumber\\ &+& O\left(n u_ns_n^{1/2}\delta_n\right)||\pmb{u}||
+ O\left(nu_n^2\rho_n\right)||\pmb{u}||^2 \nonumber\\ &-&nu_n^2\left(\pmb{u}^{\top}\pmb{B}_{\theta_0S_n\times S_n}\pmb{u}+o_p(1)\right).%\nonumber\\
\label{Q_final}
\end{eqnarray}
Therefore, taking into account that $s_n^{1/2}\delta_n=O(u_n)$ and $\rho_n\rightarrow 0$ as $n\rightarrow\infty$, together with Assumption (\ref{mom_4}), it is possible to choose a sufficiently large $C$ in such a way that the last term in (\ref{Q_final}) dominates the other terms uniformly on $||\pmb{u}||=C$. This fact completes the proof of (\ref{rates}).

Now, we will obtain the result in (\ref{minimo}). Because of $\left(\widehat{\pmb{\beta}}_{0S_n},\widehat{\theta}_0\right)$ is a local minimizer of  $\mathcal{Q}^*(\pmb{\beta}_{S_n},\theta)$ verifying (\ref{rates}), to prove (\ref{minimo}) it suffices to obtain that:
\begin{equation}
\mathcal{Q}\left(\left(\widehat{\pmb{\beta}}_{0S_n}^{\top},\pmb{0}_{p_n-s_n}^{\top}\right)^{\top},\widehat{\theta}_0\right)=\min_{||\pmb{\beta}_{S_n^c}||\leq Cu_n}\mathcal{Q}\left(\left(\widehat{\pmb{\beta}}_{0S_n}^{\top},\pmb{\beta}_{S_n^c}^{\top}\right)^{\top},\widehat{\theta}_0\right),
\end{equation}
where $\pmb{\beta}_{S_n^c}=\left(\beta_{s_n+1},\dots,\beta_{p_n}\right)^{\top}$. For that, we will show that both
\begin{equation}
\left.\frac{\partial\mathcal{Q}\left(\pmb{\beta},\widehat{\theta}_0\right)}{\partial \beta_j}\right|_{\pmb{\beta}=\pmb{\beta}^{j\vartheta}}>0  \textrm{ for $0 < \vartheta< C u_n$ }
\label{positivo}
\end{equation}
and
\begin{equation}
\left.\frac{\partial\mathcal{Q}\left(\pmb{\beta},\widehat{\theta}_0\right)}{\partial \beta_j}\right|_{\pmb{\beta}=\pmb{\beta}^{j\vartheta}}<0   \textrm{ for $-Cu_n < \vartheta<0$ }
\label{negativo}
\end{equation}
hold, where $j\in\{s_n+1,\dots,p_n\}$ and $\pmb{\beta}^{j\vartheta}$ denotes a vector with dimension $p_n$ obtained from $\left(\widehat{\pmb{\beta}}_{0S_n}^{\top},\pmb{\beta}_{S_n^c}^{\top}\right)^{\top}$ by changing their $j$th component by $\vartheta$. Simple calculations give, for $s_n+1 \leq j\leq p_n$,
%\begin{equation}
\begin{eqnarray}
\left.\frac{\partial\mathcal{Q}\left(\pmb{\beta},\widehat{\theta}_0\right)}{\partial\beta_j}\right|_{\pmb{\beta}=\pmb{\beta}^{j\vartheta}}&=&-\left(\widetilde{\pmb{X}}_{\widehat{\theta}_0}\right)_j^{\top}\left(\widetilde{\pmb{Y}}_{\widehat{\theta}_0}-\widetilde{\pmb{X}}_{\widehat{\theta}_0}\pmb{\beta}_0\right)+\left(\widetilde{\pmb{X}}_{\widehat{\theta}_0}\right)_j^{\top}\widetilde{\pmb{X}}_{\widehat{\theta}_0}\left(\pmb{\beta}^{j\vartheta}-\pmb{\beta}_0\right)\nonumber \\&+&n\mathcal{P}'_{\lambda_{jn}}\left(\left|\vartheta\right|\right)sgn(\vartheta),
\label{der_Q}
\end{eqnarray}
%end{equation}
where $\left(\widetilde{\pmb{X}}_{\widehat{\theta}_0}\right)_j$ denotes de $j$th column of $\widetilde{\pmb{X}}_{\widehat{\theta}_0}$. Therefore, to prove (\ref{positivo}) and (\ref{negativo}) it suffices to show that the sign of (\ref{der_Q}) is determined by $sgn(\vartheta)$. On the one hand, we have that
\begin{eqnarray}
\left(\widetilde{\pmb{X}}_{\widehat{\theta}_0}\right)_j^{\top}\left(\widetilde{\pmb{Y}}_{\widehat{\theta}_0}-\widetilde{\pmb{X}}_{\widehat{\theta}_0}\pmb{\beta}_0\right)&=&\left(\pmb{g}_{j,\theta_0}-\widehat{\pmb{g}}_{j,\widehat{\theta}_0}\right)^{\top}\left(\pmb{g}_{0,\theta_0}-\widehat{\pmb{g}}_{0,\widehat{\theta}_0}\right)+\left(\pmb{g}_{j,\theta_0}-\widehat{\pmb{g}}_{j,\widehat{\theta}_0}\right)^{\top}\pmb{\varepsilon}\nonumber\\
&-&\left(\pmb{g}_{j,\theta_0}-\widehat{\pmb{g}}_{j,\widehat{\theta}_0}\right)^{\top}\left(\pmb{G}_{\theta_0}-\widehat{\pmb{G}}_{\widehat{\theta}_0}\right)\pmb{\beta}_0+\pmb{\eta}_{j,\theta_0 }^{\top}\left(\pmb{g}_{0,\theta_0}-\widehat{\pmb{g}}_{0,\widehat{\theta}_0}\right)\nonumber\\&+&\pmb{\eta}_{j,\theta_0 }^{\top}\pmb{\varepsilon}-\pmb{\eta}_{j,\theta_0 }^{\top}\left(\pmb{G}_{\theta_0}-\widehat{\pmb{G}}_{\widehat{\theta}_0}\right)\pmb{\beta}_0.
\label{expr7}
\end{eqnarray}
The  quantities in expression (\ref{expr7}) can be bounded in the same way as those in (A20) in \cite{aneiros_2015}, but using our Lemma \ref{lemaA3} instead of Lemma A.3 in \cite{aneiros_2015} (note that $\widehat{\theta}_0 \subset \Theta_n$). 
Therefore, it can be obtained that
\begin{equation}
\left(\widetilde{\pmb{X}}_{\widehat{\theta}_0}\right)_j^{\top}\left(\widetilde{\pmb{Y}}_{\widehat{\theta}_0}-\widetilde{\pmb{X}}_{\widehat{\theta}_0}\pmb{\beta}_0\right)=O_p\left(n^{1/2+1/r_{\varepsilon}}\log n\right).
\label{expr8}
\end{equation}
On the other hand, we have that 
\begin{eqnarray}
\left|\left(\widetilde{\pmb{X}}_{\widehat{\theta}_0}\right)_j^{\top}\widetilde{\pmb{X}}_{\widehat{\theta}_0}\left(\pmb{\beta}^{j\vartheta}-\pmb{\beta}_0\right)\right|&\leq& \left\lvert\left\lvert\left(\widetilde{\pmb{X}}_{\widehat{\theta}_0}\right)_j\right\rvert\right\rvert \left\lvert\left\lvert\widetilde{\pmb{X}}_{\widehat{\theta}_0}\right\rvert\right\rvert \left\lvert\left\lvert\pmb{\beta}^{\vartheta j}-\pmb{\beta}_0\right\rvert\right\rvert\nonumber\\&=&\left\lvert\left\lvert\left(\widetilde{\pmb{X}}_{\widehat{\theta}_0}\right)_j\right\rvert\right\rvert \left\lvert\left\lvert\Delta_{\max}\left(\widetilde{\pmb{X}}_{\widehat{\theta}_0}^{\top}\widetilde{\pmb{X}}_{\widehat{\theta}_0}\right)\right\rvert\right\rvert \left\lvert\left\lvert\pmb{\beta}^{j \vartheta}-\pmb{\beta}_0\right\rvert\right\rvert.\label{t2.1}
\end{eqnarray}
From Lemma \ref{lemaA5} together with Assumption (\ref{mom_2}) we obtain that
\begin{eqnarray}
\left\lvert\left\lvert\left(\widetilde{\pmb{X}}_{\widehat{\theta}_0}\right)_j\right\rvert\right\rvert=O_p\left(n^{1/2}\right) \label{t2.2}
\end{eqnarray}
uniformly over $s_n+1\leq j\leq p_n$, while using a similar reasoning of that employed in proof of Lemma \ref{lemaA5} we have that
\begin{equation}
\widetilde{\pmb{X}}_{\widehat{\theta}_0}^{\top}\widetilde{\pmb{X}}_{\widehat{\theta}_0}=n\pmb{B}_{\theta_0}+o_p(n). \label{t2.3}
\end{equation}
Therefore, from (\ref{t2.1}), (\ref{t2.2}) and (\ref{t2.3}) together with the fact that $\left\lvert\left\lvert\pmb{\beta}^{j \vartheta}-\pmb{\beta}_0\right\rvert\right\rvert=O_p(u_n)$, we obtain that
\begin{equation}
\left|\left(\widetilde{\pmb{X}}_{\widehat{\theta}_0}\right)_j^{\top}\widetilde{\pmb{X}}_{\widehat{\theta}_0}\left(\pmb{\beta}^{j\vartheta}-\pmb{\beta}_0\right)\right|=O_p\left(nu_n\Delta_{\max}^{1/2}(\pmb{B}_{\theta_0})\right).
\label{expr9}
\end{equation}
Finally, using (\ref{der_Q}), (\ref{expr8}) and (\ref{expr9}) we obtain that
%\begin{equation}
\begin{eqnarray}
\left.\frac{\partial\mathcal{Q}\left(\pmb{\beta},\widehat{\theta}_0\right)}{\partial\beta_j}\right|_{\pmb{\beta}={\pmb{\beta}}^{j\vartheta}}&=&n\lambda_{jn}\left(O_p\left(n^{-1/2+1/r_{\varepsilon}}\lambda_{jn}^{-1}\log n\right)+O_p\left(\lambda_{jn}^{-1}u_n\Delta_{\max}^{1/2}(\pmb{B}_{\theta_0})\right) \right.\nonumber \\&+&\left.\lambda_{jn}^{-1}\mathcal{P}'_{\lambda_{jn}}\left(\left|\vartheta\right|\right)sgn(\vartheta)\right)\nonumber
\end{eqnarray}
%\end{equation}
Thus, taking into account our assumptions, we have proved that the sign of $\partial\mathcal{Q}\left(\pmb{\beta},\widehat{\theta}_0\right)/\partial\beta_j|_{\pmb{\beta}=\pmb{\beta}^{j\vartheta}}$ is completely determined
by that of $\vartheta$. Therefore equations (\ref{positivo}) and (\ref{negativo}) are checked and, as a consequence, the proof of (\ref{minimo}) is completed.

Because we have proven both (\ref{rates}) and (\ref{minimo}), the proof of our Theorem \ref{teorema1} concludes.
$\blacksquare $

\subsection{Proof of Theorem \ref{teorema2}}
Taking our Theorem \ref{teorema1} into account, similar steps as those used to prove Theorem 3.2(a) in \cite{aneiros_2015} can be followed to obtain the proof of our Theorem \ref{teorema2}. $\blacksquare $

\subsection{Proof of Theorem \ref{teorema3}}
Because $\left(\widehat{\pmb{\beta}}_{0S_n},\widehat{\theta}_0\right)$ is a local minimum of $Q^*(\pmb{\beta}_{S_n},\theta)$, for each $j\in S_n$ it is verified that:
\begin{equation}
\left.\frac{\partial\mathcal{Q}^*\left(\pmb{\beta}_{S_n},\theta\right)}{\partial\beta_j}\right|_{\left(\pmb{\beta}_{S_n},\theta\right)=\left(\widehat{\pmb{\beta}}_{0S_n},\widehat{\theta}_0\right)}=0.
\label{der_0}
\end{equation} 
After some Taylor expansion and using assumptions (\ref{u3}), (\ref{null_par2}), the fact that $u_n/\min_{j\in S_n}\{\lambda_{jn}\}=o(1)$ and Theorem \ref{teorema1}, we obtain that
\begin{eqnarray}
\left.\frac{\partial\mathcal{Q}^*\left(\pmb{\beta}_{S_n},\theta\right)}{\partial\beta_j}\right|_{(\pmb{\beta}_{S_n},\theta)=\left(\widehat{\pmb{\beta}}_{0S_n},\widehat{\theta}_0\right)}&=&-\left(\widetilde{\pmb{X}}_{\widehat{\theta}_0}\right)_j^{\top}\left(\widetilde{\pmb{Y}}_{\widehat{\theta}_0}-\widetilde{\pmb{X}}_{\widehat{\theta}_0{S_n}}\pmb{\beta}_{0S_n}\right) \nonumber\\ &+&\left(\widetilde{\pmb{X}}_{\widehat{\theta}_0}\right)_j^{\top}\widetilde{\pmb{X}}_{\widehat{\theta}_0{S_n}}\left(\widehat{\pmb{\beta}}_{0S_n}-\pmb{\beta}_{0S_n}\right)
\nonumber\\&+&n\mathcal{P}'_{\lambda_{jn}}\left(\left|\beta_{0j}\right|\right)\textrm{sgn}(\beta_{0j})+n\mathcal{P}''_{\lambda_{jn}}\left(\left|\beta_{0j}\right|\right)\left(\widehat{\beta}_{0j}-\beta_{0j}\right)\nonumber\\&+&O_p(nu_n^2).\nonumber
\label{der_Q*_minimo}
\end{eqnarray}
Then, in virtue of (\ref{der_0}), it can be written
\begin{eqnarray}
\pmb{0}&=&-\widetilde{\pmb{X}}_{\widehat{\theta}_0S_n}^{\top}\left(\widetilde{\pmb{Y}}_{\widehat{\theta}_0}-\widetilde{\pmb{X}}_{\widehat{\theta}_0{S_n}}\pmb{\beta}_{0S_n}\right)+\widetilde{\pmb{X}}_{\widehat{\theta}_0S_n}^{\top}\widetilde{\pmb{X}}_{\widehat{\theta}_0{S_n}}\left(\widehat{\pmb{\beta}}_{0S_n}-\pmb{\beta}_{0S_n}\right)\nonumber\\
&+&n\pmb{c}_{S_n}+n\pmb{V}_{S_n\times S_n}\left(\widehat{\pmb{\beta}}_{0S_n}-\pmb{\beta}_{0S_n}\right)+O_p\left(s_n^{1/2}nu_n^2\right).
\label{der_0_des}
\end{eqnarray}
Now, from (\ref{der_0_des}) and Lemma \ref{lemaA7} we have that
\begin{align}
&n^{1/2}\left(\pmb{B}_{\theta_0S_n\times S_n}+\pmb{V}_{S_n\times S_n}+o_p(1)\right)\left(\widehat{\pmb{\beta}}_{0S_n}-\pmb{\beta}_{0S_n}
+\left(\pmb{B}_{\theta_0S_n\times S_n}+\pmb{V}_{S_n\times S_n}+o_p(1)\right)^{-1}\pmb{c}_{S_n}\right)\nonumber\\
&=n^{-1/2}\widetilde{\pmb{X}}_{\widehat{\theta}_0S_n}^{\top}\left(\widetilde{\pmb{Y}}_{\widehat{\theta}_0}-\widetilde{\pmb{X}}_{\widehat{\theta}_0{S_n}}\pmb{\beta}_{0S_n}\right)+O_p\left(s_n^{1/2}n^{1/2}u_n^2\right).
\label{der_0_lema}
\end{align}
We should note that the first term on the right-hand side of the equality (\ref{der_0_lema}) matches the term $A_{15}$ in (\ref{descomposicion_L}) (when $\theta=\widehat{\theta}_0$ is considered in (\ref{descomposicion_L})) after multiplying it by $n^{-1/2}$ and removing $\pmb{u}$. Note also that the orders in (\ref{Op6})-(\ref{Op10}) are still true if both the vector $\pmb{u}$ is removed and $p_n$ in $r_n$ is changed by $s_n$ (note that $p_n$ comes from Lemma \ref{lemaA3}, where the maximum is taken over $p_n$ elements; in the particular case of expressions (\ref{Op6})-(\ref{Op10}), the corresponding number of elements is $s_n$). Therefore, taking into account the decomposition (\ref{expr5}) of $A_{15}$, denoting by $\gamma_n$ the maximum of the orders in Equations (\ref{Op6})-(\ref{Op10}) when $p_n$ in $r_n$ is changed by $s_n$,  and multiplying each side of (\ref{der_0_lema}) by $\pmb{A}_n\sigma_{\varepsilon}^{-1}\pmb{B}_{\theta_0S_n\times S_n}^{-1/2}$, we obtain that
\begin{align}
&n^{1/2}\pmb{A}_n\sigma_{\varepsilon}^{-1}\pmb{B}_{\theta_0S_n\times S_n}^{-1/2}\left(\pmb{B}_{\theta_0S_n\times S_n}+\pmb{V}_{S_n\times S_n}\right)\left(\widehat{\pmb{\beta}}_{0S_n}-\pmb{\beta}_{0S_n}
+\left(\pmb{B}_{\theta_0S_n\times S_n}+\pmb{V}_{S_n\times S_n}\right)^{-1}\pmb{c}_{S_n}\right)\nonumber\\
&=n^{-1/2}\pmb{A}_n\sigma_{\varepsilon}^{-1}\pmb{B}_{\theta_0S_n\times S_n}^{-1/2}\pmb{\eta}_{\theta_0S_n}^{\top}\pmb{\varepsilon}+\pmb{A}_n\pmb{B}_{\theta_0S_n\times S_n}^{-1/2}O_p\left(n^{-1/2}\gamma_n+s_n^{1/2}n^{1/2}u_n^2\right).\label{der_0_A15}
\end{align}
Note that  from assumptions in Theorem \ref{teorema3}, $(n^{-1/2}\gamma_n+s_n^{1/2}n^{1/2}u_n^2)=o(1)$ holds. Now using that $\pmb{A}_n\pmb{A}_n^{\top}\rightarrow\pmb{A}$ together with Assumption (\ref{mom_4}), we have that
\begin{eqnarray}||\pmb{A}_n\pmb{B}_{\theta_0S_n\times S_n}^{-1/2}||^2\leq ||\pmb{A}_n||^2||\pmb{B}_{\theta_0S_n\times S_n}^{-1/2}||^2=\Delta_{max}\left(\pmb{A}\right)(1+o(1))\frac{1}{\Delta_{min}\left(\pmb{B}_{\theta_0S_n\times S_n}\right)}=O(1).\nonumber
\end{eqnarray}
Therefore, expression (\ref{der_0_A15}) can be simplified in the following way:
\begin{align}
&n^{1/2}\pmb{A}_n\sigma_{\varepsilon}^{-1}\pmb{B}_{\theta_0S_n\times S_n}^{-1/2}\left(\pmb{B}_{\theta_0S_n\times S_n}+\pmb{V}_{S_n\times S_n}\right)\left(\widehat{\pmb{\beta}}_{0S_n}-\pmb{\beta}_{0S_n}
+\left(\pmb{B}_{\theta_0S_n\times S_n}+\pmb{V}_{S_n\times S_n}+o_p(1)\right)^{-1}\pmb{c}_{S_n}\right)\nonumber\\
&=n^{-1/2}\pmb{A}_n\sigma_{\varepsilon}^{-1}\pmb{B}_{\theta_0S_n\times S_n}^{-1/2}\pmb{\eta}_{\theta_0S_n}^{\top}\pmb{\varepsilon}+o_p(1).\label{der_final}\nonumber
\end{align}
As consequence, the result will be proved if we show that
\begin{equation}
n^{-1/2}\pmb{A}_n\sigma_{\varepsilon}^{-1}\pmb{B}_{\theta_0S_n\times S_n}^{-1/2}\pmb{\eta}_{\theta_0S_n}^{\top}\pmb{\varepsilon}=\sum_{i=1}^n\pmb{Z}_{ni}\overset{\textnormal{d}}{\longrightarrow}\textrm{N}(\pmb{0},\pmb{A}),
\label{nor_zni}
\end{equation}
where we use the notation $\pmb{Z}_{ni}=n^{-1/2}\pmb{A}_n\sigma_{\varepsilon}^{-1}\pmb{B}_{\theta_0S_n\times S_n}^{-1/2}\pmb{\eta}_{i,\theta_0S_n}\varepsilon_i$.

On the one hand, following exactly the same development used in \cite{aneiros_2015}, we obtain that the iid sequence of $q$-dimensional random vectors $\{\pmb{Z}_{ni}\}$ satisfies the conditions of the Lindeberg--Feller central limit  theorem. On the other hand,
since assumptions (\ref{centred_error}) and (\ref{var_indep_2}) are verified,  the reasoning used in \cite{aneiros_2015} gives us
\begin{align}
&\mathbb{E}\left(\sum_{i=1}^n\pmb{Z}_{ni}\right)=\pmb{0},\nonumber\\
&\textrm{Var}\left(\sum_{i=1}^n\pmb{Z}_{ni}\right)=\pmb{A}_n^{\top}\pmb{A}_n\rightarrow\pmb{A}.
\end{align}
Therefore, the result (\ref{nor_zni}) is verified, which completes the proof. $\blacksquare$

\subsection{Proof of Theorem \ref{teorema4}}
We have that
\begin{eqnarray}
\left|\widehat{m}_{\theta}\left(\chi\right)-m_{\theta_0}\left(\chi\right)\right|&=&\left|\sum_{i=1}^n w_{n,h,\theta}\left(\chi,\mathcal{X}_i\right)\left(m_{\theta_0}\left(\mathcal{X}_i\right)+\varepsilon_i\right)-m_{\theta_0}\left(\chi\right)+\sum_{j \in S_n}\widehat{g}_{j,\theta}(\chi)\left(\beta_{0j}-\widehat{\beta}_{0j}\right)\right|\nonumber\\
&\leq&\left|\sum_{i=1}^n w_{n,h,\theta}\left(\chi,\mathcal{X}_i\right)\left(m_{\theta_0}\left(\mathcal{X}_i\right)+\varepsilon_i\right)-m_{\theta_0}\left(\chi\right)\right|\nonumber\\
&+&\left|\sum_{j \in S_n}\left(\widehat{g}_{j,\theta}(\chi)-g_{j,\theta_0}(\chi)\right)\left(\beta_{0j}-\widehat{\beta}_{0j}\right)\right|\nonumber\\
&+&\left|\sum_{j \in S_n}g_{j,\theta_0}(\chi)\left(\beta_{0j}-\widehat{\beta}_{0j}\right)\right|\nonumber\\
&\leq&\left|\sum_{i=1}^n w_{n,h,\theta}\left(\chi,\mathcal{X}_i\right)\left(m_{\theta_0}\left(\mathcal{X}_i\right)+\varepsilon_i\right)-m_{\theta_0}\left(\chi\right)\right|\nonumber\\
&+& s_n^{1/2} \sup_{\chi\in\mathcal{C},j\in S_n,\theta\in \Theta_n}\left|\widehat{g}_{j,\theta}(\chi)-g_{j,\theta_0}(\chi)\right|\left|\left|\widehat{\pmb{\beta}}_0-\pmb{\beta}_0\right|\right|\nonumber\\
&+& s_n^{1/2}\sup_{\chi\in\mathcal{C},j\in S_n,\theta\in \Theta_n}\left|g_{j,\theta_0}(\chi)\right|\left|\left|\widehat{\pmb{\beta}}_0-\pmb{\beta}_0\right|\right|.\label{descomp}
\end{eqnarray}
Applying, to the decompositions considered in the proofs of Lemmas \ref{Lema1}-\ref{lemaA3} (see Section \ref{apendice_lemas}), the techniques used in \cite{ferraty_2010} to prove their Theorem 2, we obtain that
\begin{equation}
\sup_{\chi\in \mathcal{C},\theta\in\Theta_n}\left|\sum_{i=1}^n w_{n,h,\theta}\left(\chi,\mathcal{X}_i\right)\left(m_{\theta_0}\left(\mathcal{X}_i\right)+\varepsilon_i\right)-m_{\theta_0}\left(\chi\right)\right|=O_p\left(h^{\alpha}+\sqrt{\frac{\psi_{\mathcal{C}}\left(1/n\right)}{nf(h)}}\right)\label{sup1}
\end{equation}
(remember that, as a consequence of our assumptions on $v_n$, we have that $\psi_{\Theta_n}\left(1/n\right)=0$). In addition, Lemma \ref{lemaA3} give us
\begin{equation}
\sup_{\chi\in\mathcal{C},j\in S_n,\theta\in \Theta_n}\left|\widehat{g}_{j,\theta}(\chi)-g_{j,\theta_0}(\chi)\right|=O_p\left(h^{\alpha}+\sqrt{r_n}\right),\label{sup2}
\end{equation}
where $r_n$ was defined in (\ref{rn}). Therefore, since $\sup_{\chi\in\mathcal{C},j\in S_n}\left|g_{j,\theta_0}(\chi)\right|=O(1)$, using Theorem \ref{teorema1} and expressions (\ref{descomp}), (\ref{sup1}) and (\ref{sup2}), we can finally obtain
\begin{eqnarray}
\sup_{\theta\in\Theta_n}\sup_{\chi\in\mathcal{C}}\left|\widehat{m}_{\theta}\left(\chi\right)-m_{\theta_0}\left(\chi\right)\right|&=&O_p\left(h^{\alpha}+\sqrt{\frac{\psi_{\mathcal{C}}\left(1/n\right)}{nf(h)}}\right)\nonumber\\&+&O_p\left(s_n^{1/2}u_n\left(h^{\alpha}+\sqrt{r_n}\right)\right)+O_p\left(s_n^{1/2}u_n\right)\nonumber\\
&=&O_p\left(h^{\alpha}+\sqrt{\frac{\psi_{\mathcal{C}}\left(1/n\right)}{nf(h)}}\right)+O_p\left(s_n^{1/2}u_n\right). \ \blacksquare \nonumber
\end{eqnarray}

\subsection{Proof of Corollary \ref{cor}}
Trivial. $\blacksquare $

\subsection{Proof of Corollary \ref{final-cor}}
From the first part and second one in Condition A in Corollary \ref{final-cor}, one obtains that $\psi_{\mathcal{C}}\left(1/n\right) \approx \log n$ (see Example 4 in \citealt*{ferraty_2010}, page 338) and $f(h) \approx h$ (see Lemma 13.6 in \citealt*{ferratyvieu_2006}), respectively. Therefore, taking into account that $h \approx C(\log n / n)^{1/(2\alpha+1)}$, one has that
\begin{equation}
\label{end1}
h^{\alpha}+\sqrt{\frac{\psi_{\mathcal{C}}\left(1/n\right)}{nf(h)}}=O\left( \left(\frac{\log n}{n}\right)^{\alpha/(2\alpha + 1)}\right).
\end{equation}
In addition, taking into account that $u_n=O(\sqrt{s_n}n^{-1/2})$ (see Condition D in Corollary \ref{final-cor}), together with the fact that $s_n\approx cn^\gamma$ with $0<2\gamma\leq 1-2\alpha/(2\alpha+1)$, one obtains that 
\begin{equation}
\label{end2}
\sqrt{s_n}u_n \approx c' n^{\gamma-1/2}=O\left( \left(\frac{\log n}{n}\right)^{\alpha/(2\alpha + 1)}\right).
\end{equation}
(\ref{end1}) and (\ref{end2}) conclude the proof. $\blacksquare $

\section{Technical lemmas}
\label{apendice_lemas}
This section shows some known lemmas used to prove the results in this paper. In addition, novel lemmas, as well as their proofs, are presented. Their interest is not restricted to the proof of our theorems but they could be useful in other contexts. In fact, some assumptions used in our novel lemmas are more general than the corresponding ones imposed in our theorems.

\subsection{General assumptions}
\label{general-assump}
Let us present some additional assumptions to be used in some of the lemmas in Section \ref{results-lemas}.

\begin{description}

\item[\textit{Condition on the set of directions and the associated topology.}] The set of directions, $\Theta_n$ (see (\ref{Theta})), satisfies
\begin{equation}
	\Theta_n \subset \bigcup_{j=1}^{N_{\Theta_n,\epsilon}}B(\theta_{\epsilon,j},\epsilon), \mbox{ where } \epsilon=1/n
	\label{cover-Theta}
\end{equation}
and $N_{\Theta_n,\epsilon}$ is the minimal number of open balls in $(\mathcal{H},d(\cdot,\cdot))$ of radius $\epsilon$ which are necessary to cover $\Theta_n$.% \subset \mathcal{H}$.

\item[\textit{Conditions on the entropies and the balls in (\ref{cover-C}).}] Let $\psi_{\Theta_n}(\epsilon)$ denotes the Kolmogorov entropy of $(\Theta_n, d(\cdot,\cdot))$ (that is, $\psi_{\Theta_n}(\epsilon)=\log(N_{\Theta_n,\epsilon})$). It is assumed that:
\begin{equation}
	\exists \beta>1 \textrm{ such that } p_n\exp\left\{\left(1-\beta\log p_n\right)\left(\psi_{\mathcal{C}}\left(\frac{1}{n}\right)+\psi_{\Theta_n}\left(\frac{1}{n}\right)\right)\right\}\rightarrow 0 \textrm{ as } n\rightarrow \infty ,
	\label{fun_spaces}
\end{equation}
and
\begin{equation}
	\sup_{\theta \in \Theta_n} \max_{k \in \{1,\ldots,N_{\mathcal{C},1/n}^\theta\}}d(\chi_k^{\theta},\chi_{k^\ast}^{\theta^\ast})=O(1/n) \mbox{ (for notation, see (\ref{asterisco}))}.
	\label{chi-chi}
\end{equation}
\item[\textit{Condition linking the entropies and the small-ball probabilities.}] There exists a constant $C_{13}>0$ such that, for $n$ large enough,
\begin{equation}
	\psi_{\mathcal{C}}\left(\frac{1}{n}\right)+\psi_{\Theta_n}\left(\frac{1}{n}\right)\leq\frac{C_{13} n f(h)}{\alpha_n\log p_n}, \mbox{ where } \alpha_n\rightarrow\infty \mbox{ as } n\rightarrow\infty.
	\label{entropy}
\end{equation}
(The function $f(\cdot)$ was defined in (\ref{h33}) and (\ref{small_ball}))
\end{description}

\begin{remark}
\label{rem-app}
In the theorems presented in this paper, the condition imposed on $v_n$ (i.e., $ns_nv_n=O(hf(h))$) implies that
\begin{equation}
	\Theta_n \subset B(\theta_0,1/n). \label{cover-Theta_0}
\end{equation}
Nevertheless, that does not necessarily happen in the novel lemmas proposed in Section \ref{results-lemas} (they are applicable in more general scenarios). For this motive, the more general assumption (\ref{cover-Theta}), as well as the new assumptions (\ref{fun_spaces}),  (\ref{chi-chi}) and (\ref{entropy}), are introduced here (note that assumptions (\ref{fun_spaces_0}), (\ref{chi-chi_0}) and (\ref{entropy_0}), used in our theorems, are particular cases of (\ref{fun_spaces}), (\ref{chi-chi}) and (\ref{entropy}), respectively: it suffices to consider $N_{\Theta_n,\epsilon}=1$ and $\theta_{\epsilon,1}=\theta_0$ in (\ref{cover-Theta})). Assumptions (\ref{cover-Theta}), (\ref{fun_spaces}) and (\ref{entropy}) (together with (\ref{cover-C}) and (\ref{small_ball}) in Section \ref{assumptions}) are common when one needs to obtain uniform orders over $\mathcal{C}$ and $\Theta_n$ (see, for instance, \citealt*{wang_2016}). Finally, the condition (\ref{chi-chi}) is really specific to the functional setting addressed here and, therefore, requires a deeper discussion. It is a technical assumption that links the topologies of $(\mathcal{C}, d_\theta(\cdot,\cdot))$ and $(\Theta_n,d(\cdot,\cdot))$, allowing to bound the difference $d_\theta(\chi_k^{\theta},\cdot)-d_{\theta^\ast}(\chi_{k^\ast}^{\theta^\ast},\cdot)$ by means of bounds based on the topology of $(\Theta_n,d(\cdot,\cdot))$ (for details, see the proof of Lemma \ref{Lema1}). In fact, condition (\ref{chi-chi}) could be changed by the more general (but maybe harder to interpret) one
$$\sup_{\chi\in\mathcal{C}}\sup_{\theta \in \Theta_n} \max_{k \in \left\{1,\ldots,N_{\mathcal{C},1/n}^\theta\right\}}\left|d_\theta\left(\chi_k^{\theta},\chi\right) - d_{\theta^\ast}\left(\chi_{k^\ast}^{\theta^\ast},\chi\right) \right|=O\left(1/n\right).$$ It is worth noting that condition (\ref{chi-chi}) is satisfied if, for instance, the following condition holds:
\begin{equation}
	d(\chi_k^{\theta},\chi_{k^\ast}^{\theta^\ast}) \leq Cd(\theta,\theta^\ast), \ \mbox{uniformly on } {\theta \in \Theta_n} \mbox{ and }  k \in \left\{1,\ldots,N_{\mathcal{C},1/n}^\theta\right\}, \label{alternativa}
\end{equation}
where $C$ denotes a positive constant. Actually, condition (\ref{alternativa}) can be seen as a smoothness condition: roughly speaking, it imposes ``smooth changes'' between the coverings (\ref{cover-C}) induced by the topologies of $(\mathcal{C}, d_\theta(\cdot,\cdot))$ ($\theta \in \Theta$) when the indexes $\theta$ are close (to be more precise, see the definition of both $k^\ast$ and $\theta^\ast$ in (\ref{asterisco})).
\end{remark}

\subsection{Additional notation}
The following notation:
$$
k_{(\theta,k,\epsilon)}^{\ast}=\arg\min_{k' \in \left\{1,\ldots,N_{\mathcal{C},\epsilon}^{\theta_{j_{(\theta,\epsilon)}}}\right\}}d(\chi_{\epsilon,k}^{\theta},\chi_{\epsilon,k'}^{\theta_{\epsilon,j_{(\theta,\epsilon)}}}), \mbox{ where }
j_{(\theta,\epsilon)}=\arg\min_{j \in \{1,\ldots,N_{\Theta_n,\epsilon}\}}d(\theta,\theta_{\epsilon,j}),$$
and
\begin{equation}
r_n^{\ast}=\frac{\log p_n\left(\psi_{\mathcal{C}}\left(1/n\right)+\psi_{\Theta_n}\left(1/n\right)\right)}{nf(h)},
\label{rn.ast}
\end{equation}
generalizes the previous notation $k_{(\theta,k,\epsilon)}^{0}$ (see (\ref{varias})) and $r_n$ (see (\ref{rn})), respectively, to the more general setting considered in Section \ref{apendice_lemas}. In the sake of brevity, we will denote
\begin{equation}
\chi_{k}^{\theta}=\chi_{1/n,k}^{\theta}, \ \theta^\ast=\theta_{1/n,j_{(\theta,1/n)}} \mbox{ and } k^\ast=k_{(\theta,k,1/n)}^{\ast}.
\label{asterisco}
\end{equation}
Finally, we introduce the statistics
$$
\widehat{F}_{\theta}(\chi)=\frac{\sum_{i=1}^{n}K\left( d_{\theta}\left( \chi,\mathcal{X}_{i}\right) /h\right) }{n\mathbb{E}%
\left( K\left( d_{\theta}\left( \chi,\mathcal{X}\right) /h\right) \right)} \mbox{ and }
\widehat{g}_{j,\theta}^{\ast}(\chi)= \widehat{F}_{\theta}(\chi)\widehat{g}_{j,\theta}(\chi) \
\left(j=0,1,\ldots ,p_{n}\right),
$$
which will be used in the proofs of some of our lemmas.

\subsection{Results}
\label{results-lemas}
\begin{lemma}(Lemma 3 in \citealt*{aneiros_2008})
Let $\{V_i\}_{i=1}^n$ be a zero-mean, stationary, independent and real process verifying that
$\exists r > 4$ such that $\max_{1\leq i\leq n}\mathbb{E}|V_i|^r = O(1)$. Assume that $\{a_{ij}, i, j = 1,\dots, n\}$ is a sequence of positive numbers such that
$\max_{1\leq i,j\leq n} |a_{ij}| = O(a_n)$. Then,
\begin{equation*}
	\max_{1\leq j\leq n}\left|\sum_{j=1}^n a_{ij} V_i\right|=O_p\left(a_nn^{1/2+1/r} \log n\right).
\end{equation*}
\label{lemaA1}
The conclusion of this lemma remains unchanged when the $a_{ij}$ are random variables satisfying the conditions earlier in probability.
\end{lemma}
\begin{lemma}(Lemma A.2 in \citealt*{aneiros_2015})
Let $\{V_{ijk}\}_{i=1}^n$ $\left(1\leq j\leq u_n, 1\leq k\leq v_n\right)$ be independent random variables with zero mean and $\forall r\geq2$,
$E|V_{ijk}|^r\leq C_V(r!/2)$, where $0 < C_V <\infty$ is a constant. Assume that $\{a_{ijk}, 1\leq i\leq n,1\leq j\leq u_n,1\leq k\leq v_n\}$ is a set of positive numbers such that
$\max_{\substack{1\leq i\leq n\\ 1\leq j\leq u_n\\1\leq k\leq v_n}} |a_{ijk}| = O(a_n)$. If $u_nv_nn^{-\log n}\rightarrow0$ as $n\rightarrow\infty$ then:
\begin{equation*}
	\max_{1\leq j\leq u_n}\max_{1\leq k\leq v_n}\left|\sum_{i=1}^n a_{ijk} V_{ijk}\right|=O_p\left(a_nn^{1/2}\log n\right).
\end{equation*}
\label{lemaA2}
The conclusion of this lemma remains unchanged when the $a_{ijk}$ are random variables satisfying the conditions earlier in probability.
\end{lemma}
\begin{lemma}
\label{Lema1} Under assumptions (\ref{cover-C}), (\ref{h33}), (\ref{small_ball}), (\ref{kernel}), (\ref{chi}) and (\ref{cover-Theta})-(\ref{entropy}), if in addition $\sup_{\theta \in \Theta_n} \left<\theta,\theta\right>^{1/2}=O(1)$ and $\mathcal{X}_i$ are i.i.d, we have that there exists a positive constant, $C$, such that, for all $\epsilon >0$ and $n$ large enough,  %, for any $\epsilon>0$ and $0<\delta<1$%
$$
P\left(\sup_{\theta \in \Theta_n} \sup_{\chi\in \mathcal{C}} \left|\widehat{F}_{\theta}(\chi) - 1\right|>\epsilon \sqrt{%
	r_{n}^{\ast}}\right) \leq C\left(p_{n}^{-C\epsilon ^{2}} + \left(N_{\Theta_n,1/n}N_{\mathcal{C},1/n}\right)^{1-C\epsilon
	^{2}\log p_n} \right).
$$
\end{lemma}

\noindent \textit{Proof of Lemma \ref{Lema1}.} One can write
\begin{equation}
\sup_{\theta \in \Theta_n} \sup_{\chi\in \mathcal{C}} \left|\widehat{F}_{\theta}(\chi) - 1\right| \leq F_1 + F_2 + F_3,
\label{L1.1}
\end{equation}
where we have denoted
$$
F_1=\sup_{\theta \in \Theta_n} \sup_{\chi\in \mathcal{C}} \left|\widehat{F}_{\theta}(\chi) - \widehat{F}_{\theta}\left(\chi^{\theta}_{k_{(\theta,\chi,1/n)}}\right)\right|,$$
$$F_2=\sup_{\theta \in \Theta_n} \sup_{\chi\in \mathcal{C}} \left|\widehat{F}_{\theta}\left(\chi^{\theta}_{k_{(\theta,\chi,1/n)}}\right) - \mathbb{E}\left(\widehat{F}_{\theta}\left(\chi^{\theta}_{k_{(\theta,\chi,1/n)}}\right)\right)\right|
$$
and
$$
F_3=\sup_{\theta \in \Theta_n} \sup_{\chi\in \mathcal{C}} \left|\mathbb{E}\left(\widehat{F}_{\theta}\left(\chi^{\theta}_{k_{(\theta,\chi,1/n)}}\right)\right)-\mathbb{E}\left(\widehat{F}_{\theta}(\chi)\right)\right|,
$$
with
\begin{equation}
k_{(\theta,\chi,1/n)}=\arg\min_{k \in \{1,\ldots,N_{\mathcal{C},1/n}^{\theta}\}}d_{\theta}(\chi,\chi^{\theta}_{1/n,k})
\label{k-}
\end{equation}
(see also notation (\ref{asterisco})). Taking assumptions (\ref{h33}) and (\ref{small_ball}) into account, and following the same lines as \cite{ferraty_2010} to prove their Lemma 8, one can obtain that there exists a positive constant, $C$, such that, for all $\epsilon >0$ and $n$ large enough
\begin{equation}
P\left(F_1>\epsilon \sqrt{%
	r_{n}^{\ast}}\right) \leq Cp_{n}^{-C\epsilon ^{2}} \ \mbox{ and } F_3=O\left(\sqrt{r_{n}^{\ast}}\right).
\label{L1.2}
\end{equation}
Focusing now on $F_2$, we have that
\begin{equation}
F_2 \leq F_{21} + F_{22} + F_{23},
\label{L1.3}
\end{equation}
where we have denoted
$$
F_{21}=\sup_{\theta \in \Theta_n} \max_{k \in \left\{1,\ldots,N_{\mathcal{C},1/n}^\theta\right\}} \left|\widehat{F}_{\theta}(\chi_k^{\theta}) - \widehat{F}_{\theta^\ast}\left(\chi^{\theta^\ast}_{k^\ast}\right)\right|,$$
$$F_{22}=\sup_{\theta \in \Theta_n} \max_{k \in \left\{1,\ldots,N_{\mathcal{C},1/n}^\theta\right\}} \left|\widehat{F}_{\theta^\ast}\left(\chi^{\theta^\ast}_{k^\ast}\right) - \mathbb{E}\left(\widehat{F}_{\theta^\ast}\left(\chi^{\theta^\ast}_{k^\ast}\right)\right)\right|
$$
and
$$
F_{23}=\sup_{\theta \in \Theta_n} \max_{k \in \left\{1,\ldots,N_{\mathcal{C},1/n}^\theta\right\}}\left|\mathbb{E}\left(\widehat{F}_{\theta^\ast}\left(\chi^{\theta^\ast}_{k^\ast}\right)\right)-\mathbb{E}\left(\widehat{F}_{\theta}\left(\chi_k^{\theta}\right)\right)\right|.
$$
(For notation, see (\ref{asterisco})) First, we consider the terms $F_{21}$ and $F_{23}$. Taking into account that
$$
\left|d_\theta\left(\chi_k^{\theta},\mathcal{X}_i\right) - d_{\theta^\ast}\left(\chi_{k^\ast}^{\theta^\ast},\mathcal{X}_i\right) \right| \leq d\left(\chi_k^{\theta},\chi_{k^\ast}^{\theta^\ast}\right)\left<\theta,\theta\right>^{1/2}+
d\left(\chi_{k^\ast}^{\theta^\ast},\mathcal{X}_i\right)d\left(\theta,\theta^\ast\right),
$$
and using assumptions (\ref{chi}) and (\ref{chi-chi}), together with the condition $\sup_{\theta \in \Theta_n} \left<\theta,\theta\right>^{1/2}=O(1)$, we obtain that
$$
\sup_{\theta \in \Theta_n} \max_{k \in \left\{1,\ldots,N_{\mathcal{C},1/n}^\theta\right\}}\left|d_\theta\left(\chi_k^{\theta},\mathcal{X}_i\right) - d_{\theta^\ast}\left(\chi_{k^\ast}^{\theta^\ast},\mathcal{X}_i\right) \right|=O\left(1/n\right).
%\label{L1.4}
$$
Therefore, similar steps as those used to obtain (\ref{L1.2}) can be followed to get
\begin{equation}
P\left(F_{21}>\epsilon \sqrt{%
	r_{n}^{\ast}}\right) \leq Cp_{n}^{-C\epsilon ^{2}} \ \mbox{ and } F_{23}=O\left(\sqrt{r_{n}^{\ast}}\right).
\label{L1.4}
\end{equation}
Finally, we study the term $F_{22}$. We have that
$$
F_{22}=\max_{j \in \left\{1,\ldots,N_{\Theta_n,1/n}\right\}} \max_{k \in \left\{1,\ldots,N_{\mathcal{C},1/n}^{\theta_j}\right\}} \left|\widehat{F}_{\theta_j}\left(\chi^{\theta_j}_{k}\right) - \mathbb{E}\left(\widehat{F}_{\theta_j}\left(\chi^{\theta_j}_{k}\right)\right)\right|.
$$
Therefore, similar steps as those used in \cite{ferraty_2010} (page 345) to obtain their Lemma 8 can be followed to get
\begin{equation}
P\left(F_{22} > \epsilon\sqrt{r_n^{\ast}}\right) \leq C \left(N_{\Theta_n,1/n} N_{\mathcal{C},1/n}\right)^{1-C\epsilon^2\log p_n}.
\label{L1.5}
\end{equation}
(\ref{L1.1})-(\ref{L1.5}) complete the proof. $\blacksquare $

\begin{lemma}
\label{Lema2} Under the assumptions of Lemma \ref{Lema1}, if in addition assumptions (\ref{centred_error}), (\ref{var_indep_1}), (\ref{var_indep_2}) and (\ref{mom_1}) hold, then there exists a positive constant, $C$, such that, for all $\epsilon >0$ and $n$ large enough
$$
P\left(\sup_{\theta \in \Theta_n} \sup_{\chi\in \mathcal{C}} \left|\widehat{g}_{j,\theta}^{\ast}(\chi)-\mathbb{E}\left(\widehat{g}_{j,\theta}^{\ast}(\chi)\right)\right|>\epsilon \sqrt{%
	r_{n}^{\ast}}\right) \leq C\left(p_{n}^{-C\epsilon ^{2}} + \left(N_{\Theta_n,1/n}N_{\mathcal{C},1/n}\right)^{1-C\epsilon
	^{2}\log p_n} \right),
$$
uniformly on $j=0,1,\ldots,p_n$.
\end{lemma}
\noindent \textit{Proof of Lemma \ref{Lema2}.} This proof can be easily obtained combining the techniques considered in the proof of Lemma 11 in \cite{ferraty_2010} with the decompositions (adapted to the new setting) used in the proof of our Lemma \ref{Lema1}. $\blacksquare$
\begin{lemma}
\label{Lema3} Under assumptions (\ref{Theta}), (\ref{small_ball}), (\ref{kernel}), (\ref{smooth}) and (\ref{chi}), if in addition $v_n$ in (\ref{Theta}) verifies $v_n=O(h)$, we have that%
$$
\sup_{\theta \in \Theta_n} \sup_{\chi\in \mathcal{C}} \left|\mathbb{E}\left(\widehat{g}_{j,\theta}^{\ast}(\chi)\right)-g_{j,\theta_0}(\chi)\right|=O\left(h^{\alpha}\right),
$$
uniformly on $j=0,1,\ldots,p_n$.
\end{lemma}

\noindent \textit{Proof of Lemma \ref{Lema3}.} Firstly, we note that, if $d_{\theta}(\mathcal{X},\chi)<h$ holds, then, from the fact that $v_n=O(h)$ together with Assumption (\ref{chi}), we have that
\begin{eqnarray*}
d_{\theta_0}(\mathcal{X},\chi)&\leq& \left|d_{\theta_0}(\mathcal{X},\chi)-d_{\theta}(\mathcal{X},\chi)\right|+d_{\theta}(\mathcal{X},\chi)=\left|\left<\mathcal{X}-\chi,\theta_0-\theta\right>\right|+d_{\theta}(\mathcal{X},\chi)\nonumber\\
&\leq & \left<\mathcal{X}-\chi, \mathcal{X}-\chi\right>^{1/2}  \left<\theta_0-\theta,\theta_0-\theta\right>^{1/2}+d_{\theta}(\mathcal{X},\chi) \leq Ch. %\label{dist0}
\end{eqnarray*}
This fact, together with both assumptions (\ref{small_ball}) and (\ref{smooth}), allows to follow the same steps as in the proof of Lemma 10 in \cite{ferraty_2010} and to get the thesis of our Lemma \ref{Lema3}. $\blacksquare $

\begin{lemma}
\label{lemaA3}
Under assumptions (\ref{centred_error}), (\ref{cover-C}), (\ref{Theta}), (\ref{h33}), (\ref{small_ball}), (\ref{kernel})-(\ref{mom_1}) and (\ref{cover-Theta})-(\ref{entropy}), if in addition $v_n$ in (\ref{Theta}) verifies $v_n=O(h)$ and $p_n\rightarrow\infty$ as $n\rightarrow\infty$, then
$$
\max_{0\leq j \leq p_n} \sup_{\theta \in \Theta_n} \sup_{\chi\in \mathcal{C}}\left\{ \left|\widehat{g}_{j,\theta}(\chi)-g_{j, \theta_0}(\chi)\right|\right\} =O_{p}\left( h^{\alpha }+\sqrt{r_n^{\ast}}\right).
$$
\end{lemma}	
\noindent \textit{Proof of Lemma \ref{lemaA3}.} It verifies that
\begin{equation}
\widehat{k}_{j,\theta}(\chi)=\left(\widehat{g}_{j,\theta}^{\ast}(\chi)-\mathbb{E}\left(\widehat{g}_{j,\theta}^{\ast}(\chi)\right)\right) + 
\left(\mathbb{E}(\widehat{g}_{j,\theta}^{\ast}(\chi))-g_{j,\theta_0}(\chi)\right) + \left(1-\widehat{F}_{\theta}(\chi)\right)g_{j,\theta_0}(\chi), \label{L4.1}
\end{equation}
where we have denoted
\begin{equation}
\widehat{k}_{j,\theta}(\chi)=\widehat{F}_{\theta}(\chi)\left(\widehat{g}_{j,\theta}(\chi)-g_{j,\theta_0}(\chi)\right).
\label{L4.0}
\end{equation}
Therefore, from Lemmas \ref{Lema1}-\ref{Lema3} together with (\ref{L4.1}), we obtain that there exists a positive constant, $C$, such that, for all $\epsilon >0$ and $n$ large enough,
\begin{equation}
P\left(\sup_{\theta \in \Theta_n} \sup_{\chi\in \mathcal{C}}\left\vert \widehat{k}_{j}(\chi)\right\vert >\epsilon
\left(h^{\alpha }+\sqrt{r_{n}^{\ast}}\right)\right) \leq C\left(p_{n}^{-C\epsilon ^{2}} + \left(N_{\Theta_n,1/n}N_{\mathcal{C},1/n}\right)^{1-C\epsilon
	^{2}\log p_n} \right),
\label{L4.2}
\end{equation}%
uniformly on $j=0,1,\ldots,p_n$. In addition, taking Lemma \ref{Lema1} into account together with the facts that $r_{n}^{\ast}\rightarrow 0$ and $%
p_{n}\rightarrow \infty $ as $n\rightarrow \infty$, we obtain that, for any $0<\delta<1$ and $n$ large enough,
\begin{eqnarray}
P\left( \inf_{\theta \in \Theta_n} \inf_{\chi\in \mathcal{C}}\widehat{F}_{\theta}(\chi)\geq \delta \right) &\geq &1-P\left( \sup_{\theta \in \Theta_n} \sup_{\chi\in \mathcal{C}}\left\vert \widehat{F}_{\theta}(\chi)-1\right\vert
>1-\delta \right)  \notag \\
&\geq &1-P\left( \sup_{\theta \in \Theta_n} \sup_{\chi\in \mathcal{C}}\left\vert \widehat{F}_{\theta}(\chi)-1\right\vert >\delta
\epsilon \sqrt{r_{n}}\right)  \notag \\
&\geq &1-C\left(p_{n}^{-C\delta ^{2}\epsilon ^{2}} + \left(N_{\Theta_n,1/n}N_{\mathcal{C},1/n}\right)^{1-C\delta
	^{2}\epsilon ^{2}\log p_n} \right)\geq 1/2. \nonumber\\ \label{L4.3}
\end{eqnarray}%
Now, from (\ref{L4.2}) and (\ref{L4.3}) have that that 
\begin{eqnarray}
P\left( \sup_{\theta \in \Theta_n} \sup_{\chi\in \mathcal{C}}\frac{\left\vert \widehat{k}_{j}(\chi)\right\vert }{\widehat{F}_{\theta}(\chi)}%
>\epsilon (h^{\alpha }+\sqrt{r_{n}^{\ast}})\right) &\leq &P\left( \sup_{\theta \in \Theta_n} \sup_{\chi\in \mathcal{C}}\frac{\left\vert \widehat{k}_{j}(\chi)\right\vert }{\widehat{F}_{\theta}(\chi)}>\epsilon (h^{\alpha
}+\sqrt{r_{n}^{\ast}}) \mid \inf_{\theta \in \Theta_n} \inf_{\chi\in \mathcal{C}}d(\chi)\geq \delta \right)  \notag \\
&&+P\left( \inf_{\theta \in \Theta_n} \inf_{\chi\in \mathcal{C}}\widehat{F}_{\theta}(\chi)<\delta \right)  \notag \\
&\leq&C\left(p_{n}^{-C\delta ^{2}\epsilon ^{2}} + \left(N_{\Theta_n,1/n}N_{\mathcal{C},1/n}\right)^{1-C\delta
	^{2}\epsilon ^{2}\log p_n} \right).   \label{L4.4}
\end{eqnarray}%
Finally, using (\ref{L4.4}) we obtain that,
\begin{eqnarray*}
P\left( \sup_{0\leq j\leq p_{n}}\sup_{\theta \in \Theta_n} \sup_{\chi\in \mathcal{C}}\frac{\vert \widehat{k}_{j}(\chi)\vert }{\widehat{F}_{\theta}(\chi)}>\epsilon (h^{\alpha }+\sqrt{r_{n}^{\ast}})\right) &\leq
&\sum_{j=0}^{p_{n}}P\left( \sup_{\theta \in \Theta_n} \sup_{\chi\in \mathcal{C}}\frac{\vert
	k_{j}(\chi)\vert }{d(\chi)}>\epsilon (h^{\alpha }+\sqrt{r_{n}^{\ast}})\right) \\
&\leq &Cp_{n}\left(p_{n}^{-C\delta ^{2}\epsilon ^{2}} + (N_{\Theta_n,1/n}N_{\mathcal{C},1/n})^{1-C\delta \label{L4.5}
	^{2}\epsilon ^{2}\log p_n} \right).
\end{eqnarray*}%
The proof is completed taking the notation (\ref{L4.0}) into account and choosing $\delta$ and $\epsilon$ in (\ref{L4.5}) such that $%
C\delta ^{2}\epsilon ^{2}=\beta$ (see (\ref{fun_spaces})). $\blacksquare$
\begin{lemma}
Under assumptions of Lemma \ref{lemaA3}, if in addition $p_n/n^{\log n}\rightarrow0$ as $n\rightarrow\infty$, then
$$
\max_{0\leq j \leq p_n} \sup_{\theta \in \Theta_n} \sup_{\chi\in \mathcal{C}}\left\{ \left|\widehat{g}_{j,\theta}(\chi)-\widehat{g}_{j, \theta_0}(\chi)\right|\right\} =O_{p}\left(\frac{v_n}{h f(h)} \right) .
$$
\label{lema_novo}
\end{lemma}	

\noindent \textit{Proof of Lemma \ref{lema_novo}.} Let us denote 
\begin{equation*}
\widehat{F}_{\theta}^*(\chi)=\frac{\sum_{i=1}^n K(d_{\theta}(\mathcal{X}_i,\chi)/h)}{n f(h)}
\mbox{ and }
\widehat{g}_{j,\theta}^{**}(\chi)=\widehat{F}_{\theta}^*(\chi)\widehat{g}_{j,\theta}(\chi).
\end{equation*}
It is easy to obtain the decomposition
\begin{eqnarray}
\widehat{g}_{j,\theta}(\chi)-\widehat{g}_{j,\theta_0}(\chi)&=&\frac{1}{\widehat{F}_{\theta}^*(\chi)}\widehat{g}_{j,\theta}^{**}(\chi)-\frac{1}{\widehat{F}_{\theta_0}^*(\chi)}\widehat{g}_{j,\theta_0}^{**}(\chi)\nonumber\\
&=&\frac{1}{\widehat{F}_{\theta}^*(\chi)}\left(\widehat{g}_{j,\theta}^{**}(\chi)-\widehat{g}_{j,\theta_0}^{**}(\chi)\right)+\widehat{g}_{j,\theta_0}^{**}(\chi)\left(\frac{1}{\widehat{F}_{\theta}^*(\chi)}-\frac{1}{\widehat{F}_{\theta_0}^*(\chi)}\right).\label{decomp_g}
\end{eqnarray}
Now, we are going to analyze the terms in (\ref{decomp_g}). Let us denote
$$
Z_{i0}=Y_i \mbox{  and  } Z_{ij}=X_{ij} \ (1 \leq j \leq p_n).
$$
We have that
\begin{eqnarray}
\left|\widehat{g}_{j,\theta}^{**}(\chi)-\widehat{g}_{j,\theta_0}^{**}(\chi)\right|
&\leq&\frac{\sum_{i=1}^n\left|d_{\theta}\left(\mathcal{X}_i,\chi\right)-d_{\theta_0}\left(\mathcal{X}_i,\chi\right)\right||Z_{ij}|}{nhf(h)}\leq \frac{v_n}{h f(h)}\frac{1}{n}\sum_{i=1}^n|Z_{ij}|\nonumber\\
&=&O_p\left(\frac{v_n}{h f(h)}\right),\label{dif_g_ast_ast}
\end{eqnarray}
uniformly over $0 \leq j \leq p_n$, $\theta\in\Theta_n$ and $\chi\in\mathcal{C}$. (Note that the first inequality in (\ref{dif_g_ast_ast}) is a consequence of Assumption (\ref{kernel}) while assumptions (\ref{Theta}) and (\ref{chi}) give the second one. Finally, the equality comes from Assumption (\ref{mom_1}) together with Lemma \ref{lemaA2} applied to the centred variables $\{|Z_{ij}|-\mathbb{E}(|Z_{ij}|\}_i)$) In a similar way (considering $Z_{ij}=1$ in (\ref{dif_g_ast_ast})), one obtains
\begin{eqnarray}
\left|\widehat{F}_{\theta}^*(\chi)-\widehat{F}_{\theta_0}^*(\chi)\right|=O_p\left(\frac{v_n}{nh f(h)}\right),\label{dif_F_ast}
\end{eqnarray}
uniformly over $\theta\in\Theta_n$ and $\chi\in\mathcal{C}$. 

Now, we focus on $\widehat{g}_{0,\theta_0}^{**}(\chi)$. From assumptions (\ref{h33}), (\ref{small_ball}) and (\ref{kernel}), we have that (see \citealt*{ferratyvieu_2006})
\begin{equation*}
Cf(h)\leq \mathbb{E}\left(K\left(d_{\theta}(\mathcal{X}_i,\chi)/h \right)\right)\leq C'f(h),
\end{equation*}
where $C$ and $C'$ denote positive constants. Therefore, there exist positive constants, $C^*$ and $C^{'*}$, such that
\begin{equation}
C^*\widehat{F}_{\theta}(\chi)\leq \widehat{F}_{\theta}^*(\chi)\leq C^{*'}\widehat{F}_{\theta}(\chi).\label{cota_F_ast1}
\end{equation}
On the one hand, from (\ref{cota_F_ast1}) together with Lemma \ref{Lema1} we obtain that %there exist positive constants, $C^{\ast}$ and $C^{*'}$, such that
\begin{equation}
C^*\left(1+o_p(1)\right)\leq \widehat{F}_{\theta}^*(\chi)\leq C^{*'}\left(1+o_p(1)\right),\label{cota_F_ast2}
\end{equation}
uniformly over $\theta\in\Theta_n$ and $\chi\in\mathcal{C}$. On the other hand, from the uniform convergence of $\widehat{g}_{j,\theta}(\chi)$ to $g_{j,\theta_0}(\chi)$ (see Lemma \ref{lemaA3}) together with the fact that 
\begin{equation}
\max_{0\leq j\leq n}\max_{1\leq i\leq n} |g_{j,\theta_0}(\mathcal{X}_i)|=O(1)
\nonumber
\end{equation}
(see Assumption (\ref{mom_1})), we obtain that
\begin{equation}
\max_{0\leq j \leq p_n}\sup_{\theta \in \Theta_n^\ast}\max_{1\leq i\leq n} |\widehat{g}_{j,\theta}(\mathcal{X}_i)|=O_p(1). \label{cota_g_hat}
\end{equation}
As a consequence of (\ref{cota_F_ast2}) and (\ref{cota_g_hat}), we have that
\begin{equation}
\widehat{g}_{j,\theta}^{**}(\chi)=\widehat{F}_{\theta}^*(\chi)\widehat{g}_{j,\theta}(\chi)=O_p(1),\label{g_ast_ast}
\end{equation}
uniformly over $0 \leq j \leq p_n$, $\theta\in\Theta_n$ and $\chi\in\mathcal{C}$. 

Finally, (\ref{decomp_g}), (\ref{dif_g_ast_ast}), (\ref{dif_F_ast}), (\ref{cota_F_ast2}) and (\ref{g_ast_ast}) give the result of the lemma. $\blacksquare$

\begin{lemma}(Lemma A.4 in \citealt*{aneiros_2015})
\label{lemaA4}
Let us assume that $\pmb{\eta}_{i,\theta_0}^{\top}$ $(i=1,\dots,n)$ are iid random vectors. If, in addition, $\mathbb{E}\left(\eta_{\theta_0,1j}^4\right)<C$ uniformly on $1\leq j\leq p_n$, then
\begin{equation*}
	\pmb{u}^{\top}\left(\pmb{\eta}_{\theta_0}^{\top}\pmb{\eta}_{\theta_0}-n\pmb{B}_{\theta_0}\right)\pmb{u}=O_p\left(n^{1/2}p_n\right),  \textrm{ uniformly over } \{\pmb{u}\in\mathbb{R}^{p_n}, ||\pmb{u}||=M\}.
\end{equation*}
\end{lemma}

\begin{lemma}
\label{lemaA5}
Let us assume that $\pmb{\eta}_{i,\theta_0}$ $(i=1,\dots,n)$ are iid random vectors. Under assumptions (\ref{cover-C}), (\ref{Theta}), (\ref{h33}), (\ref{small_ball}), (\ref{kernel}), (\ref{smooth}), (\ref{chi})-(\ref{mom_2}) and (\ref{cover-Theta})-(\ref{entropy}) ($g_{0,\theta_0}$ and $Y$ not included in assumptions (\ref{smooth}) and (\ref{mom_1}), respectively), if in addition $v_n$ in (\ref{Theta}) verifies $v_n=O(h)$, $p_n\rightarrow\infty$,  $p_n=o(n^{1/2})$, $nh^{4\alpha}=O(1)$ and  $$\log^2 p_n= O\left(n\left(\frac{f(h)}{\psi_{\mathcal{C}}\left(1/n\right)+\psi_{\Theta_n}\left(1/n\right)}\right)^2\right)$$
as $n\rightarrow\infty$, then we have that
\begin{equation*}
	\pmb{u}^{\top}\left(\widetilde{\pmb{X}}_{\theta}^{\top}\widetilde{\pmb{X}}_{\theta}-n\pmb{B}_{\theta_0}\right)\pmb{u}=o_p(n), \textrm{ uniformly over } \{\pmb{u}\in\mathbb{R}^{p_n}, ||\pmb{u}||=M\} \textrm{ and over } \theta\in\Theta_n.
\end{equation*}
\end{lemma}
\noindent \textit{Proof of Lemma \ref{lemaA5}.} To prove this result, the outline used in proof of Lemma A.5 in \cite{aneiros_2015} can be exactly followed, but now  our Lemma \ref{lemaA3} it is needed to conclude instead of Lemma A.3 in \cite{aneiros_2015}. $\blacksquare$
\begin{lemma}(Lemma A.6 in \citealt*{aneiros_2015})
\label{lemaA6}
Let us assume that $\pmb{\eta}_{i,\theta_0 S_n}$ $(i=1,\dots,n)$ are iid random vectors. If in addition $s_n^2/n=o(1)$ and $\max_{1\leq j\leq s_n}\mathbb{E}(\eta_{1j,\theta_0}^4)=O(1)$, then 
\begin{equation*}
	\left\lvert\left\lvert n^{-1}\pmb{\eta}_{\theta_0S_n}^{\top}\pmb{\eta}_{\theta_0S_n}-\pmb{B}_{\theta_0S_n\times S_n}\right\rvert\right\rvert=o_p(1),
\end{equation*}
\end{lemma}

\begin{lemma}
\label{lemaA7}
Let us assume that $\pmb{\eta}_{i,\theta_0 S_n}$ $(i=1,\dots,n)$ are iid random vectors. If, in addition, assumptions (\ref{cover-C}), (\ref{Theta}), (\ref{h33}), (\ref{small_ball}), (\ref{kernel}), (\ref{smooth}), (\ref{chi}), (\ref{mom_1}) and (\ref{cover-Theta})-(\ref{entropy}) hold (but using $s_n$ instead of $p_n$, and $g_{0,\theta_0}$ and $Y$ not included in assumptions (\ref{smooth}) and (\ref{mom_1}), respectively), and $v_n$ in (\ref{Theta}) verifies $v_n=O(h)$, $p_n\rightarrow\infty$ and $$\max\left\{h,s_nh^{\alpha},s_n^2/n,s_n^2\log s_n/\left(n\left(\frac{f(h)}{\psi_{\mathcal{C}}\left(1/n\right)+\psi_{\Theta_n}\left(1/n\right)}\right)\right)\right\}=o(1),$$ then 
\begin{equation*}
	n^{-1}\widetilde{\pmb{X}}_{\theta S_n}^{\top}\widetilde{\pmb{X}}_{\theta S_n}=\pmb{B}_{\theta_0S_n\times S_n}+o_p(1), \textrm{ uniformly over } \theta\in\Theta_n.
\end{equation*}

\end{lemma}
\noindent \textit{Proof of Lemma \ref{lemaA7}.} The scheme of proof of Lemma A.7 in \cite{aneiros_2015} can be exactly followed, taking into account that Lemma A.3 and Lemma A.6 in \cite{aneiros_2015} should be replaced by Lemma \ref{lemaA3} and Lemma \ref{lemaA6}, respectively, of this paper. $\blacksquare$

%%%%%%%%%%%%%%%%%%%%%%%%%%%%%%%%%
\begin{lemma}(Lemma A.8 in \citealt*{aneiros_2015})
\label{lemaA8}
Let us assume that $\left(\pmb{\eta}_{i,\theta_0S_n}^{\top},\varepsilon_i\right)$ $(i=1,\dots,n)$ are i.i.d random vectors with mean zero, and $\{\pmb{\eta}_{i,\theta_0S_n}\}$  and $\{\varepsilon_i\}$ are independent. If, in addition, $\mathbb{E}(\varepsilon_i)=0$, $\mathbb{E}(\varepsilon_i^2)<C$,  $\mathbb{E}(\eta_{1j,\theta_0})<C$ uniformly on $1\leq j\leq s_n$, then
\begin{equation*}
	\pmb{\varepsilon}^{\top}\pmb{\eta}_{\theta_0S_n}^{\top}\pmb{u}=O_p\left(n^{1/2}s_n^{1/2}\right), \textrm{ uniformly over } \{\pmb{u}\in\mathbb{R}^{p_n}, ||\pmb{u}||=M\}.
\end{equation*}

\end{lemma}
	
\end{document}